# Complete Embedded Minimal Surfaces of Finite Total Curvature

David Hoffman*  
MSRI  
1000 Centennial Drive  
Berkeley, CA 94720

Hermann Karcher†  
Math. Institut  
Universitat Bonn  
D-5300 Bonn, Germany

March 16, 1995

# Contents

*Supported by research grant DE-FG02-86ER25015 of the Applied Mathematical Science subprogram of the Office of Energy Research, U.S. Department of Energy, and National Science Foundation, Division of Mathematical Sciences research grants DMS-9101903 and DMS-9312087 at the University of Massachusetts.

†Partially supported by Sonderforschungsbereich SFB256 at Bonn.

# 1 Introduction

We will survey what is known about minimal surfaces $S \subset \mathbb{R}^3$, which are complete, embedded, and have finite total curvature: $\int_S |K| dA < \infty$. The only classically known examples of such surfaces were the plane and the catenoid. The discovery by Costa [14, 15], early in the last decade, of a new example that proved to be embedded sparked a great deal of research in this area. Many new examples have been found, even families of them, as will be described below. The central question has been transformed from whether or not there are any examples except surfaces of rotation to one of understanding the structure of the space of examples.

Up to this point, *every new example* of a complete, embedded minimal surface of finite total curvature has been discovered first by using the global version of the Enneper-Riemann-Weierstrass representation, which is essentially due to Osserman [58, 59]. This involves knowledge of the compact Riemann structure of the minimal surface, as well as its Gauss map and other geometric-analytic data. One of our goals is to show how this is done in the case that has been most completely analyzed, namely surfaces with genus $\geq 1$ and three topological ends. An important quality of this construction is that the Riemann surface and the meromorphic data are constructed simultaneously under the assumption of symmetry. Moreover, once this is done, parameters must be found in order to have a well-defined finite-total-curvature surface. This parameter search is typically done by computer using a combination of relatively simple numerical routines and relatively complex graphics tools ([8, 41]). In many cases a full theoretical analysis, as is done here in Section 4 for the three-ended surfaces of Theorem 3.3 has yet to be carried out. Moreover, solution of the period problem does not at all guarantee that the surfaces are embedded. In Section 5 we present examples of Callahan, Hoffman, Karcher, Meeks and Wohlgemuth that lie in one parameter families containing both embedded and immersed surfaces. In fact the period problem and embeddedness are totally independent issues. There are examples of Weierstrass data meeting all necessary conditions (Proposition 2.4) for embeddedness, for which the period problem is not solvable (a genus-one example with two catenoid ends does not exist but Weierstrass data for such a surface –even a very symmetric one– does) and others for which the period problem is solvable for a family of surfaces, which are embedded outside of a compact set of $R^3$, but are not embedded.

The survey is organized as follows. In Section 2 we present the basic tools of the subject, the most important of which is the Weierstrass-Enneper representation. In Section 2.2 we describe the construction of Chen and Gackstatter [10], which produces a genus-one surface with the symmetries and end behavior of Enneper's surface. To our knowledge this complete minimal surface of finite total curvature was the first one explicitly constructed by first specifying a geometric property — in this case, end behavior — and then deriving the necessary Weierstrass data. It is not, of course, embedded but its construction, as presented here, has most of the features of the construction of the three-ended examples in Section 4. In Section 2.3, the hypothesis of embeddedness is used to derive relationships between the geometry of the surface and its analytic representation. Propositions 2.3, 2.4 and 2.5 gather together all necessary conditions including the relationship between flux, logarithmic growth rates of the ends and residues of the complex differential of the height function.

In Section 3, we present the few global rigidity theorems that are known. (Theorems 3.1 and 3.4, due to Lopez-Ros [51], Schoen [67] and Costa [16].) We present a proof (in Section 3.1) of the Lopez-Ros theorem, which states that a complete minimal surface of genus zero and finite total curvature must be the plane or the catenoid. Our proof follows that of Perez-Ros [60]. We also state the existence result, Theorem 3.3, for the three-ended, complete minimal surfaces with genus $k-1$ and $k$ vertical planes of symmetry ([30]). The details of the construction of these surfaces are presented in Section 4. We include here the estimation of the parameters that solve the period problem when $k > 2$. The values of the parameters that close the periods determine the logarithmic growth rates of the ends of these surfaces. For the surface to be embedded, they must lie in a certain range, which they (happily) do. This is done in Sections 4.5–4.9.

In Section 5, we survey other known examples and discuss what little is known about the structure of the space of complete embedded minimal surfaces of finite total curvature. Section 5.2 presents some conjectures about this space.

Finite total curvature implies finite topology, even without the additional assumption of embeddedness. The converse is not true; the helicoid is simply connected, nonflat and periodic, so its curvature is infinite, while its topology is finite. Up until recently, the helicoid was the only known embedded minimal surface with finite topology and infinite total curvature. In 1992, we discovered, with Fusheng Wei, a complete embedded minimal surface of genus one with one end — asymptotic to the helicoid — that has infinite total curvature [32, 33]. The details of this construction are outside the scope of this work. However, the extent to which finite topology implies finite total curvature is discussed in Section 6.

Section 7 discusses the index of stability of a complete minimal surface. The basic results [22, 25, 26], the equivalence of finite index and finite total curvature, and the fact that the index is completely defined by the Gauss map, are discussed. There is, as yet, no known relationship between embeddedness and properties of the index. This final section is therefore, strictly speaking, misplaced in this survey. However the ideas and techniques may, in the long run, prove useful in the study of embeddedness of minimal surfaces.

**Acknowledgements.**

The authors would like to thank Jim Hoffman and Ed Thayer at the Center for Geometry, Analysis, Numerics and Graphics, at the University of Massachusetts, Amherst for the preparation of the computer graphics illustrations in this article. Peter Norman, and Pascal Romon read earlier versions and gave us useful critical comments. Harold Rosenberg helped clarify some important points in Section 2.3. We are grateful to M. Kotani and S. Nayatani for useful comments on Section 7. Robert Bryant and Michael Wolf gave us many useful comments. The authors enjoyed the hospitality of MSRI where the final work on this paper was completed.

## 2  Basic theory and the global Weierstrass representation

Let $\mathcal{S}, \langle\rangle$ be an oriented surface with a Riemannian metric, and let $\nabla$ denote its Riemannian connection. For a smooth function $f\colon \mathcal{S} \to \mathrm{R}$, we will denote its differential by $df$ and its gradient by $\mathrm{grad}\, f$. They are related by $Vf = df(V) = \langle \mathrm{grad}\, f, V\rangle, V \in T\mathcal{S}$, which can be thought of as the defining equation for $\mathrm{grad}\, f$. If $w$ is a one form, its covariant derivative $\nabla w$ is defined by the relation

$$(\nabla_U w)V = Uw(V) - w(\nabla_U V). \tag{2.1}$$

The divergence operator is $div = tr\nabla$ and the Laplacian of a smooth function $f$ is given by $\Delta f := div$ grad $f$. Note that

$$\langle \nabla_U \text{ grad } f, V\rangle = U\langle \text{grad } f, V\rangle - (\nabla_U V)f = (\nabla_U df)(V). \tag{2.2}$$

Hence

$$\Delta f = \Sigma(\nabla_{e_i} df)(e_i). \tag{2.3}$$

The operator $R_{\frac{\pi}{2}}$ on vector fields, which rotates each $T_p M$ by $\frac{\pi}{2}$ in the positive direction, is covariantly parallel: i.e. $\nabla Rot_{\frac{\pi}{2}} = 0$. This simply means that for any vector fields $U, V$, $Rot_{\frac{\pi}{2}}(\nabla_U V) = \nabla_U Rot_{\frac{\pi}{2}} V$, which can be easily checked. This implies that for a vector field $V$

$divV = 0 \Leftrightarrow \nabla Rot_{\frac{\pi}{2}} V$ is symmetric,

(i.e. $\langle \nabla_U Rot_{\frac{\pi}{2}} V, \widetilde{U}\rangle = \langle \nabla_{\widetilde{U}} Rot_{\frac{\pi}{2}} V, U\rangle$ for all vector fields $U, \widetilde{U}$).

Associated to a vector field $W$ is the dual one form $\omega$ defined by $\omega(U) = \langle U, W\rangle$. One can verify by elementary means that

$d\omega = 0 \Leftrightarrow \nabla W$ is symmetric (i.e. $\langle \nabla_U W, \widetilde{U}\rangle = \langle \nabla_{\widetilde{U}} W, U\rangle$ for all vector fields $U, \widetilde{U}$.)

$\Leftrightarrow W$ is locally the gradient of some differentiable function.

Combining these two equivalences and applying them to $W = \text{grad } f$, where $f : \mathcal{S} \to \text{R}$ is some smooth function, we have

$\Delta f = div(\text{grad } f) = 0 \Leftrightarrow Rot_{\frac{\pi}{2}} (\text{grad } f)$ is (locally) the gradient of some function $f^*$.

Smooth functions $f : \mathcal{S} \to \text{R}$ satisfying $\Delta f = 0$ are called *harmonic functions*. For a harmonic function $f$, the locally defined functions $f^*$, where grad $f^* = Rot_{\frac{\pi}{2}} (\text{grad } f)$ are called *conjugate* harmonic functions associated to $f$. Note that $df^* = -df \circ Rot_{\frac{\pi}{2}}$. Because grad $f$ and grad $f^*$ are orthogonal and have the same length, the locally defined mapping given by $(f, f^*)$ is conformal wherever grad $f \neq 0$. Thus, $f + if^*$ is a local conformal coordinate in a neighborhood of any point $p \in \mathcal{S}$ where $df_p \neq 0$.

Suppose our surface $\mathcal{S}$ is given as an oriented immersed surface in $\text{R}^3$. Let $X$ denote its position vector; we think of $X$ as an immersion $X : \mathcal{S} \to \text{R}^3$. $\mathcal{S}$ inherits a Riemannian metric $\langle \cdot, \cdot \rangle$ from $\text{R}^3$ and its Riemannian connection for vector fields can be considered to be the projection of differentiation in $\text{R}^3$; if $U, V$ are vector fields on $\mathcal{S}$, $\nabla_U V$ is defined by the relation

$$DX(\nabla_U V) = [U(DX(V))]^T, \tag{2.4}$$

where $[]^T$ denotes projection onto $DX_p(T_p\mathcal{S})$. For $p \in \mathcal{S}$ let $N(p)$ denote the unit normal to $X(\mathcal{S})$ at $X(p)$. The endomorphism $Rot_{\frac{\pi}{2}}$ on $T_p \mathcal{S}$ is just the pullback by $DX$ of counter-clockwise rotation around $N(p)$ in $DX[T_p\mathcal{S}]$. We think of $N$ as the Gauss mapping $N:\mathcal{S} \to S^2$ and denote its differential by $DN$. By identification of $T_p\mathcal{S}$ with $DX_p(T_p\mathcal{S})$ and $DX_p(T_p\mathcal{S})$ with $T_{N(p)}S^2$ by rigid translation in $\mathrm{R}^3$, we may think of $DN_p$ as an endomorphism of $T_p\mathcal{S}$. As such it is called the *shape operator*, and denoted by $S$. The eigenvalues of $S$ are called the *principal curvatures*, and its eigenvectors are called *principal directions*. The average of the principal curvatures is the *mean curvature* and is denoted by $H$.

An important formula with strong consequences for minimal surfaces is:

$$\Delta X = -2HN. \tag{2.5}$$

By $\Delta X$, we mean $(\Delta x_1, \Delta x_2, \Delta x_3)$ where $(x_1, x_2, x_3)$ are the Euclidean coordinates in $\mathrm{R}^3$. We will verify (2.5) as follows. Writing $DX = (dx_1, dx_2, dx_3)$ and applying (2.1) coordinatewise, we have

$$(\nabla_U DX)(V) = U(DX(V)) - DX(\nabla_U V)$$
$$= U(DX(V)) - [U(DX(V))]^T \quad \text{(by 2.4)}$$
$$= [U(DX(V))]^N \quad \text{(where } []^N \text{ denotes projection onto } DX(T\mathcal{S})^\perp)$$
$$= (U(DX(V)) \cdot N)N.$$

The scalar factor of $N$ in the last expression is ($\pm$, depending on your background) the Second Fundamental Form. Because $DX(V) \cdot N = 0$, $U(DX(V)) \cdot N = -DX(V) \cdot DN(U) = -\langle S(U), V\rangle$, where $S$ is the shape operator. Thus

$$(\nabla_U DX)(V) = -\langle S(U), V\rangle N.$$

From (2.2) we conclude that $\Delta X = -(trS)N$, which is (2.5), as desired.

**Remark 2.1** *Since the covariant derivative is most often developed on Riemannian manifolds without a given immersion, we indicate how the above definition is used for practical computations. For example, if $W$ is a vector field in $R^3$ we can consider its tangential component*

$$W^T := W - \langle W, N\rangle N$$

*along the immersed surface $X(S)$. This also gives us a vector field $V$ on the domain surface $S$ such that $DX(V) = W^T$. Strictly speaking it is the covariant derivative of $V$ that was defined, but one often allows abbreviations, e.g. writing $W$ rather than $W \circ X$, in computations like:*

$$\nabla_U(W^T) := DX(\nabla_U V)$$
$$= (U(DX(V)))^T$$
$$= (U(W - \langle W, N\rangle N))^T$$
$$= (UW)^T - \langle W, N\rangle DN(U).$$

**Minimal Surfaces**

A minimal surface, $\mathcal{S}$, is one for which $H \equiv 0$; from (2.5) we can conclude that the coordinate functions of a minimal surface are harmonic. We also may conclude from the discussion above that if $x$ is a coordinate function, $x + ix^*$ defines a conformal coordinate chart in a neighborhood of any $p \in \mathcal{S}$ where $dx \neq 0$. Since $X$ is assumed to be an immersion, we conclude that at every $p \in \mathcal{S}$, at least one of the coordinate functions $x_i$ satisfies $dx_{i|p} \neq 0$. Thus minimal surfaces inherit naturally, from their Euclidean coordinates, an atlas of holomorphic functions that define its Riemann surface structure.

Because $trS = H \equiv 0$ on a minimal surface $\mathcal{S}$, the Gauss mapping $N \colon \mathcal{S} \to S^2$ satisfies $\langle DN(U), DN(V) \rangle = \kappa^2 \langle U, V \rangle$, where $\pm \kappa$ are the principal curvatures, and also reverses orientation. We orient $S^2 \subset \mathrm{R}^3$ with the outward-pointing normal (so that *its* Gauss map is the identity and its principal curvatures equal to $+1$) and let $\sigma$ be stereographic projection from the north pole to the complex plane considered to be the $(x_1, x_2)$-plane, also positively oriented. With these choices, $\sigma$ is orientation-*reversing* and $g := \sigma \circ N \colon \mathcal{S} \to \mathrm{C} \cup \{\infty\}$ orientation-*preserving* and conformal, whenever $DN \neq 0$. Thus $g$ is a *meromorphic function* on $\mathcal{S}$, considered now to be a Riemann surface.

We have shown that a minimal surface has a natural Riemann surface structure, with respect to which the Gauss map is meromorphic. Now we will us this information to realize a minimal surfaces as the real parts of holomorphic curves. Recall that each coordinate function $x$ gives rise, locally, to a conjugate harmonic function $x^*$ and globally to a holomorphic differential form

$$dx + idx^* = dx - i(dx \circ Rot_{\frac{\pi}{2}}).$$

Associated to our minimal immersion $X = (x_1, x_2, x_3) \colon S \to \mathrm{R}^3$ we also have the "conjugate" minimal immersion $X^* = (x_1^*, x_2^*, x_3^*)$ (which is well defined only on some covering $\widetilde{S}$ of $S$). The relation $dX^* = -dX \circ Rot_{\frac{\pi}{2}}$ shows that $X$ and $X^*$ are locally immersions and that $\Psi := X + iX^*$ satisfies $d\Psi \circ Rot_{\frac{\pi}{2}} = i \cdot d\Psi$. In particular $\Psi$ is holomorphic and

$$X = Re\Psi = Re \int d\Psi = Re \int (dX + idX^*).$$

With respect to any local holomorphic parameter $z = u_1 + iu_2$ we have the usual connection with real differentiation:

$$d\Psi = \Psi' dz$$
$$= \left(\frac{\partial X}{\partial u_1} + i\frac{\partial X^*}{\partial u_1}\right)dz = \left(\frac{\partial X}{\partial u_1} - i\frac{\partial X}{\partial u_2}\right)dz.$$

It is convenient to define for a real-valued differentiable function $w$

$$w_z = \frac{\partial w}{\partial z} := \frac{1}{2}\left(\frac{\partial w}{\partial u_1} - i\frac{\partial w}{\partial u_2}\right),$$

and extend this function linearly to complex-valued functions:

$$(v + iw)_z = v_z + iw_z.$$

If $f = v + iw$ is holomorphic, $f_z = f'$; the desirability of this identity is the reason for the factor of $\frac{1}{2}$. The Cauchy-Riemann equation for $f$ is $v_z = iw_z$, and extending this operation coordinate-wise to $\Psi = X + iX^*$, we have

$$\Psi' = \Psi_z = (X_z + iX^*_z)$$
$$= 2X_z = \frac{\partial X}{\partial u_1} - i\frac{\partial X}{\partial u_2}$$
$$= 2iX^*_z.$$

We find an additional property by computing

$$(\Psi'^2) = 4X_z^2 = \left|\frac{\partial X}{\partial u_1}\right|^2 - \left|\frac{\partial X}{\partial u_2}\right|^2 - 2i\frac{\partial X}{\partial u_1} \cdot \frac{\partial X}{\partial u_2} = 0,$$

because all holomorphic reparametrizations are conformal. The image $\Psi(\widetilde{S})$ in $\mathbf{C}^3$ is therefore called a holomorphic null-curve. Conversely, we see that the real part of a holomorphic null-curve is *conformally* parametrized and, of course, harmonic with respect to the given paramterization. Since conformal changes of the metric do not change the harmonicity of functions, the real part of a holomorphic null-curve is also harmonic with respect to the induced Riemannian metric and hence minimal.

Having established a firm connection with complex analysis, we now tie in some more geometry by rewriting $d\Psi$ in terms of the Gauss map. This leads to a representation formulated by Enneper, Riemann and Weierstrass. Recall that the Gauss map became a meromorphic function after we distinguished a "vertical" direction for the stereographic projection; we take this to be the $x_3$-axis. Now we make full use of the fact that our surface is a complex curve: While quotients of one-forms on a surface are not functions, the quotients $\frac{d\Psi_1}{d\Psi_3}, \frac{d\Psi_2}{d\Psi_3}$ of holomorphic one-forms on the complex curve are functions, meromorphic functions on the underlying Riemann surface. Moreover, because of $X_z \cdot X_z = 0$, these functions are quadratically related:

$$\left(\frac{d\Psi_1}{d\Psi_3}\right)^2 + \left(\frac{d\Psi_2}{d\Psi_3}\right)^2 = -1.$$

Finally, since

$$\frac{\partial X}{\partial u_1} = Re(d\Psi(\frac{\partial}{\partial u_1})) \quad \text{and}$$
$$\frac{\partial X}{\partial u_2} = Re(d\Psi(\frac{\partial}{\partial u_2}))$$

are orthogonal tangent vectors to the minimal surface, we must be able to compute the normal, i.e., the Gauss map $g$, from $\frac{d\Psi_1}{d\Psi_3}, \frac{d\Psi_2}{d\Psi_3}$. With a small trick, this procedure can be reversed. Put

$$f := \frac{-(d\Psi_1 + id\Psi_2)}{d\Psi_3}$$
$$= \frac{d\Psi_3}{d\Psi_1 - id\Psi_2} \quad (\text{because } (\frac{d\Psi_1}{d\Psi_3})^2 + (\frac{d\Psi_2}{d\Psi_3})^2 = -1).$$

Then we have
$$d\Psi = (\frac{1}{2}(\frac{1}{f} - f), \frac{i}{2}(\frac{1}{f} + f), 1)d\Psi_3,$$

which now depends only on one meromorphic function $f$ and one holomorphic differential $d\Psi_3$. From now on, we will write $d\Psi_3$ as $dh$. Recall

$$dh := d\Psi_3 = dx_3 - idx_3 \circ Rot_{\frac{\pi}{2}}.$$

We know already that $f$ must be closely related to the Gauss map. In fact, we find that inverse stereographic projection of $f$, namely

$$(2Ref, 2Imf, |f|^2 - 1)/(|f|^2 + 1)$$

is orthogonal to the tangent vectors $Re(\frac{d\Psi}{d\Psi_3}), Im(\frac{d\Psi}{d\Phi_3})$ since

$$\left\langle (\frac{1}{2}(\frac{1}{f} - f), \frac{i}{2}(\frac{1}{f} + f), 1), (2Ref, 2Imf, |f|^2 - 1) \right\rangle_{\mathrm{C}} = 0 \in \mathrm{C}.$$

This identifies $f$ as the stereographic projection of the Gauss map of the minimal surface. We now state the "Weierstrass representation theorem" with emphasis on its global formulation due to Osserman [59].

**Theorem 2.1** *Suppose $S$ is a minimal surface in $R^3$, $M$ its Riemann surface, $g = \sigma \circ N$ the stereographic projection of its Gauss map. Then $S$ may be represented (up to a translation) by the conformal immersion*

$$X(p) = Re \int \Phi, \text{ where} \tag{2.6}$$

$$\Phi = (\phi_1, \phi_2, \phi_3) = \left((g^{-1} - g)\frac{dh}{2}, i(g^{-1} + g)\frac{dh}{2}, dh\right). \tag{2.7}$$

*Conversely, let $M$ be a Riemann surface, $g: M \to \mathrm{C} \cup \{\infty\}$ a meromorphic function and $dh$ a holomorphic one-form on $M$. Then (2.6) and (2.7) define a conformal minimal mapping of some covering of $M$ into $R^3$, which is regular provided the poles and zeros of $g$ coincide with the zeros of $dh$. The mapping $X$ is well-defined on $M$ if and only if no component of $\Phi$ in (??) has a real period. That is*

$$Period_\alpha(\Phi) =: Re \oint_\alpha \Phi = 0 \tag{2.8}$$

*for all closed curves $\alpha$ on $M$.*

We now compute the basic geometric quantities of $M$ in terms of $g$ and $dh$. The metric induced on $M$ by $X$ can be expressed as

$$ds^2 = \frac{1}{2}|\Phi|^2 = \frac{1}{4}(|g| + |g|^{-1})^2|dh|^2, \tag{2.9}$$

The metric will be complete provided $\int_\delta ds = \infty$ for every divergent curve $\delta$ on $M$. The Gauss curvature of this metric may be computed as the ratio of the area stretching of $dN$ to that of $dX$. Thus:
$$K = \frac{-4|dg|^2}{(|g|+|g|^{-1})^2|dh|^2} \cdot \frac{4}{(1+|g|^2)^2}.$$
The last term is the conformal stretching of stereographic projection. The minus sign comes from the fact that stereographic projection reverses orientation. We may rewrite this as
$$K = \frac{-16}{(|g|+|g|^{-1})^4} \left|\frac{dg/g}{dh}\right|^2. \tag{2.10}$$
To express the shape operator $S = DX \cdot DN$ in terms of $g$ and $dh$, we compute $X_{zz}dz^2$ and take the real part of the dot product with $N$:
$$X_{zz}dz^2 = \frac{1}{4}(X_{11} - X_{22} - 2iX_{12}) \cdot (du_1^2 - du_2^2 + 2idu_1 du_2).$$
Since in general $X_{ij} \cdot N = -X_i \cdot N_j = -\langle S(\frac{\partial}{\partial u_i}), \frac{\partial}{\partial u_j}\rangle := -s_{ij}$ and $X_{11} \cdot N = -X_{22} \cdot N$ because we are using conformal coordinates on a minimal surface:
$$-2Re((X_{zz} \cdot N)dz^2) = (s_{11}du_1^2 + s_{22}du_2^2 + 2s_{12}du_1 du_2). \tag{2.11}$$

The right-hand side is the quadratic form associated with the shape operator $S$, expressed in terms of the conformal coordinates $(u_1, u_2)$, where $z = u_1 + iu_2$. The left-hand-side of (2.11) can be written, using (2.7), as $-Re(\Phi_z dz \cdot N)$. Differentiating (2.7) gives
$$\Phi_z dz = \frac{-1}{2}((g^{-1}+g), i(g^{-1}-g), 0)\frac{dg}{g}dh,$$
and, using the fact that
$$N = \sigma^{-1} \circ g = \frac{2}{1+|g|^2}(Re\, g, Im\, g, \frac{|g|^2-1}{2}),$$
we have from (2.11) that
$$\Sigma s_{ij}du_i du_j = Re\{\frac{dg}{g}dh\}.$$
For a tangent vector $V = v_1\frac{\partial}{\partial u_1} + v_2\frac{\partial}{\partial u_2}$, we may express the Second Fundamental Form (2.11) as the real part of a holomorphic quadratic differential:
$$\langle S(V), V\rangle = Re\{\frac{dg}{g}(V) \cdot dh(V)\}. \tag{2.12}$$

This formula allows us to conclude that a curve $c$ on $X(M)$ is
$$\text{asymptotic} \iff \frac{dg}{g}(\dot{c}) \cdot dh(\dot{c}) \in i\mathrm{R}.$$
$$\text{principal} \iff \frac{dg}{g}(\dot{c}) \cdot dh(\dot{c}) \in \mathrm{R}. \tag{2.13}$$

For a minimal surface $X$ as in (2.7) and (2.8) with the Weierstrass data $\{g, dh\}$ on $M$, we make the following definition.

**Definition 2.1** *The* associate surfaces *to $X$ are the minimal immersions $X_\theta$ given by the Weierstrass data $\{g, e^{i\theta} dh\}$, $0 \leq \theta \leq \frac{\pi}{2}$, using (2.7) and (2.8). The* conjugate surface $X^*$ *is equal to $-X_{\frac{\pi}{2}}$, where $X_{\frac{\pi}{2}}$ is the associate surface with data $\{g, i\,dh\}$.*

**Remark 2.1.1** *The Weierstrass data $\{e^{i\theta} g, dh\}$ produces the same minimal surface as $\{g, dh\}$ rotated by an angle $\theta$ around the vertical axis.*

One sees immediately that $X_\theta$ may not be well-defined on $M$ because the period condition (2.8) is not necessarily satisfied. We have

$$X^* = -\,Im \int \Phi, \text{ and}$$

$$X_\theta = (\cos\theta) X - (\sin\theta) X^*$$

locally. Clearly $X_\theta$ is well-defined on $M$, for all $\theta$, $0 < \theta \leq \frac{\pi}{2}$, if and only if $X^*$ is well-defined on $M$. Also $X^{**} = -X$. By definition, all members of the associate family share the same Gauss map and from (2.8) it follows that they are locally isometric.

**The Schwarz Reflection Principle for Minimal Surfaces**

A curve on any surface $M$ in $\mathrm{R}^3$ is a straight line if and only if its geodesic curvature and its normal curvature vanish. A principal curve on $M$ that is not a straight line has the property that it is a geodesic if and only if it lies in a plane orthogonal to the surface. Thus (2.13) can be used on geodesics to identify lines and planar principal curves. From (2.13) we see that principal geodesics of $X$ correspond to straight lines of $X^*$, and vice versa. The Gauss map along a principal geodesic has values parallel to the plane of this curve, and the same values along the conjugate straight line – so that the line is orthogonal to the plane. The Schwarz reflection principle for minimal surfaces states that minimal surfaces with such curves must possess Euclidean symmetries.

*If a minimal surface contains a line segment $L$, then it is symmetric under rotation by $\pi$ about $L$. (If a minimal surface is bounded by a line segment $L$, it may be extended by rotation by $\pi$ about $L$ to a smooth minimal surface containing $L$ in its interior.)*

*If a nonplanar minimal surface contains a principal geodesic— necessarily a planar curve— then it is symmetric under reflection in the plane of that curve. (If a minimal surface meets a plane orthogonally on its boundary the surface may be extended by reflection through the plane to a smooth minimal surface with this curve in its interior.)*

We will prove this below. Finding a geodesic on a minimal surface satisfying one of the conditions of (2.13) is equivalent to finding a straight line or a planar geodesic, which with the Reflection Principle implies the existence of a Euclidean symmetry.

Suppose a minimal surface $S$ contains a line segment $L$; without loss of generality we may assume that $L$ is a portion of the $x_3$-axis containing $\vec{0}$. The surface $S^*$ conjugate to $S$ contains a planar geodesic in a horizontal plane. Because the coordinate functions of $S$ are harmonic $z := x_3 + ix_3^*$ defines a conformal parameter in a neighborhood of $L$. To avoid confusion, we write $z = u_1 + iu_2$ and we may translate so that $(x_1, x_2, x_3)(0) = \vec{0}$. Then

$$x_3(u_1) = u_1, \ x_3^*(u_1) = 0$$

$$x_1(u_1) = x_2(u_1) = 0.$$

Thus $x_3 + ix_3^*$ maps a line segment on the real axis to the real axis, while $x_1 + ix_1^*$ and $x_2 + ix_2^*$ map this same segment to the imaginary axis.

The Schwarz Reflection Principle for complex analytic functions states that if an analytic function $f$, defined in a neighborhood of a segment of the real axis, maps the segment of the real axis into a line $\ell \subset \mathrm{C}$ then $f(\bar{z}) = \rho \circ f(z)$, where $\rho$ is reflection in $\ell$. (Moreover if $f$ is defined as a one-sided neighborhood of the real segment, the extension of $f$ by $f(\bar{z}) := \rho \circ f(z)$ defines an extension of $f$ analytic on a neighborhood of the real line segment).

Since reflection in the real (resp. imaginary) axis is $z \to \bar{z}$ (resp. $z \to -\bar{z}$), we may conclude that

$$x_3 + ix_3^*(\bar{z}) = \overline{x_3 + ix_3^*} \quad \text{by definition, and}$$

$$x_k + ix_k^*(\bar{z}) = -(\overline{x_k + ix_k^*}) \quad k = 1, 2.$$

Hence

$$X(\bar{z}) := (x_1, x_2, x_3)(\bar{z}) = (-x_1, -x_2, x_3)(z)$$

$$X^*(\bar{z}) := (x_1^*, x_2^*, x_3^*)(\bar{z}) = (x_1^*, x_2^*, -x_3^*)(z).$$

That is, $S$ is symmetric under rotation about the $x_3$-axis, and $S^*$ is symmetric under reflection in the $(x_1, x_2)$-plane, as claimed. (The statements in parentheses of the Schwarz Reflection Principle for minimal surfaces follow from the statements in parentheses for the classical reflection principle.) This completes the proof of the Schwarz Reflection Principle.

**The Catenoid and the Helicoid** As a simple illustration, consider the catenoid, which can be represented as $M = \mathrm{C} - \{0\}$ with Weierstrass data $\{g, dh\} = \{z, \frac{dz}{z}\}$. See Figure 2.0.

The mapping $\rho_\theta(z) \to e^{i2\theta}\bar{z}$ is a reflection in the line $L_\theta$ given by $c(t) = te^{i\theta}$. Using (2.9) we see that $\rho_\theta$ is an isometry, and therefore its fixed point set $X(L_\theta)$ is a geodesic. The mapping (2.6) is well defined because there is only one homology cycle on which to check (2.8), namely that generated by the cycle $|z| = 1$, and:

$$\mathrm{Period}_{|z|=1}(\Phi) = Re \int_{|z|=1} \left( \left( \frac{dz}{z^2} - dz \right), i \left( \frac{dz}{z^2} + dz \right), \frac{dz}{z} \right)$$

$$= Re(0, 0, 2\pi i) = \vec{0}.$$

Along rays $c(t) = te^{i\theta}, t > 0, \dot{c}(t) = e^{i\theta}$ and

$$\frac{dg}{g}(e^{i\theta}) \cdot dh(e^{i\theta}) = \left( \frac{e^{i\theta}}{te^{i\theta}} \right)^2 = t^{-2} \in \mathrm{R}.$$

Figure 2.0 The catenoid and helicoid.

The catenoid is the image of $C-\{0\}$ under the Weierstrass mapping with $g = z$ and $dh = cdz/z$, $c = 1$. The image of $\gamma := \{|z| = 1\}$ is the waist-circle of the catenoid (left). When $c = i$, the Weierstrass mapping is multivalued; its third coordinate has a vertical period of length $2\pi$ on $\gamma$ (center). The multivalued mapping produces the helicoid, which can also be produced by reflection across the boundary lines of the image. For $c = e^{i\theta}, \theta \neq 0$, mod $\pi$, the image of $\gamma$ is a helix and the full image surface is not embedded. For any $c, |c| = 1$, the full image surface is complete and has infinite total curvature unless $c = 1$, in which case we have the catenoid with total curvature equal to $-4\pi$.

Hence each ray is a planar geodesic in a vertical plane containing the vector $(\cos\theta, \sin\theta, 0)$. Any two reflections must generate a rotation about a common axis and we can conclude, as we know, that the catenoid is a surface of rotation. Since

$$x_3(z) = Re \int_{z_0}^{z} \frac{dz}{z} = \ln \frac{|z|}{|z_0|},$$

we see that the surface diverges to $\infty$ (*resp.* $-\infty$) in the vertical direction as $|z| \to \infty$ (*resp.* $0$).

If we consider the conjugate surface $X^*$ of the catenoid, we have Weierstrass data $\{z, idz/z\}$ on $C - \{0\}$. But now

$$\text{Period}_{|z|=1}\{\Phi^*\} = Re \int_{|z|=1} i\Phi = Re(0, 0, i\ln z) = (0, 0, 2\pi).$$

Hence $X^*$ is *not* well-defined on $C - \{0\}$. It gives rise to a multi-valued mapping, the image of which we know to be fibred by horizontal lines and invariant under a vertical translation. This surface is, of course, the helicoid. Note that since the helicoid and catenoid are locally conjugate, they are locally isometric. But they are not congruent because the catenoid contains no lines. The Gauss map of the catenoid covers $S^2$ exactly once, omitting two points. Hence its total curvature is $-4\pi$. The helicoid is nonflat and periodic; hence it has infinite total curvature.

In the case of finite total curvature, $\int_M K dA > -\infty$, much more is known about the Weierstrass representation.

## 2.1 Finite total curvature

The representation Theorem 2.1 can be significantly sharpened when the minimal surface in question is complete and has finite total curvature.

**Theorem 2.2 (Osserman, [48, 59])** *Let $X: M \to R^3$ be a complete conformal minimal immersion with finite total curvature (not necessarily embedded). Then:*

i) *$M$ is conformally diffeomorphic to $\overline{M}_k - \{p_1, \ldots, p_r\}$ where $\overline{M}_k$ is a closed Riemann surface of genus $k$ and $p_1 \ldots, p_r$ are points in $\overline{M}_k$, $r \geq 1$;*

ii) *$X$ is proper;*

iii) *The Gauss map $N: M \to S^2$, which is meromorphic on $M$, extends to a meromorphic function on $\overline{M}_k$; the holomorphic one-form $dh$ extends to a meromorphic one-form on $\overline{M}_k$.*

iv) The total curvature is an integer multiple of $4\pi$ and satisfies

$$\int_M KdA \leq -4\pi(k+r-1) = 2\pi(\chi(M)-r), \tag{2.14}$$

where $k$ and $r$ are integers defined in Statement i), and $\chi(M)$ is the Euler characteristic of $M$.

Theorem 2.2 is a prescription for constructing finite total curvature examples. We will give an example of how to use it in Section 2.2.

In what follows, we will assume that $X\colon M \to \mathbf{R}^3$ is a complete, regular, conformal minimal immersion that has finite total curvature, and that $M = \overline{M}_k - \{p_1,\ldots,p_r\}$, where $\overline{M}_k$ is a compact, connected Riemann surface of genus $k$.

**Definition 2.2** Let $D_j$ be a punctured neighborhood of $p_j \in \overline{M}_k$. We will refer to $X(D_j) = E_j$ as an end of $M$, and define, for fixed $R > 0$,

$$S_{R,j} = \{q \in S^2 | Rq \in E_j\} \tag{2.15}$$

where $S^2$ denotes the unit sphere in $\mathbf{R}^3$.

**Theorem 2.3 (Gackstatter [24] Jorge-Meeks [42])** Let $S_{R,j}$ be as defined in (2.15).

i) For each $j = 1,\ldots,r$, $S_{R,j}$ converges smoothly, as $R \to \infty$, to a great circle, covered an integral number of times;

ii) Let $d_j$ be the multiplicity of the great circle $\lim_{R\to\infty} S_{R,j}$. Then

$$\int_M KdA = 2\pi\left(2(1-k) - r - \sum_{j=1}^r d_j\right) \tag{2.16}$$

$$= 2\pi\left(\chi(M) - \sum_{j=1}^r d_j\right).$$

Note that (2.16) gives an interpretation of the difference between the right- and left-hand sides of (2.14).

**Definition 2.3** The *total spinning* of $S$ is the integer quantity

$$n(S) = \sum_{j=1}^r d_j.$$

It is immediate from the definitions of $n(M)$ and $d_j$ that $n(M) \geq r$, with equality if and only if each end has multiplicity one. This may be interpreted in terms of the end at $p_j$. We will say that $M$ has an "embedded end at $p_j$" provided $X$ embeds some punctured neighborhood of $p_j$. The end at $p_j$ is embedded if and only if $d_j = 1$, as follows from the definition of $d_j$ in Theorem 2.3 ii).

**Remark 2.2** *For later use in Section 7 we will need to estimate the degree, d, of the Gauss map of a complete minimal surface M in terms of the genus. We have from (2.16 )*

$$d = (k + \frac{r}{2} - 1 + \frac{n(M)}{2}).$$

*If $r \geq 2$ then clearly $d \geq k+1$. If $r = 1$, then $n(M)$ is just $d_1$ and must be odd, and $d_1 = 1$, if and only if the end of the surface is embedded. Since $\frac{r}{2} + \frac{n(M)}{2} \geq 1$ with equality if and only if $r = 1$ and $n(M) = d_1 = 1$, it follows that $d \geq k+1$ on a complete minimal surface, unless it has precisely one end, and that end is embedded. However such an exceptional minimal surface must be a plane, as we will show below. Hence on a complete nonplanar minimal surface*

$$d \geq k+1.$$

Several arguments can be used to show that a complete immersed minimal surface $M$ of finite total curvature with a single embedded end is a plane. The simplest uses the halfspace theorem (Theorem 2.4 below in Section 2.3.3), which says that a complete, properly immersed minimal surface that lies in a halfspace must be a flat plane. Because we have finite total curvature and an embedded end, we can use Proposition 2.3 iv) and Remark 2.4 i) to conclude the same thing as follows. The end in question is asymptotic to the end of a plane or the end of a catenoid. In either case, the surface $M$ lies in a halfspace. If the end is catenoidal it is clear that the plane defining the smallest halfspace containing $M$ is tangent to $M$, violating the Maximum Principle.(See Section 2.3.3.) If the end of $M$ is planar, then $M$ lies in some slab. But again the boundary planes of the smallest slab are tangent to $M$, a contradiction unless the slab has zero thickness; i.e., $M$ is a plane.

## 2.2 The example of Chen-Gackstatter

Chen and Gackstatter [10, 11] were the first to construct a complete, finite-total-curvature example by explicitly solving the period problem. In retrospect, and from our point of view, what they did was put a handle in Enneper's surface, making it a genus-one surface with one end, while preserving all of the symmetries and end behavior. See Figure 2.2. This was the first known complete minimal torus with finite total curvature and one end. We will present this construction, following Karcher [43], Barbosa and Colares [2] and Thayer [71, 70].

Enneper's surface is given by the data

$$g = z \qquad dh = zdz \qquad (2.17)$$

on C, the sphere minus the point at infinity. There are no nontrivial closed cycles, so there is no period problem. Moreover, the metric induced by using (2.17) in (2.7) and (2.8) is given as in (2.9) by $ds = \frac{1}{2}(|z|^2 + 1)dz$ and is invariant under the maps $z \rightarrow e^{i2\theta}\overline{z}$, for any fixed $\theta$. The line $L_\theta$ through $\vec{0}$ given by $e^{i\theta}t, \infty < t < \infty$, is the fixed-point set of this mapping and so is mapped to a geodesic on Enneper's surface. The normal curvature (multiplied by the length squared of the tangent vector) of this geodesic is given, according to (2.12) by $Re\{e^{2i\theta}\}$. By (2.13), $L_\theta$ is mapped to a straight line on Enneper's surface if $\theta = \pm\frac{\pi}{4}$, and to a planar geodesic if and only if $\theta = 0, \frac{\pi}{2}$. Assuming we start the integration in (2.6) at 0, this implies that Enneper's surface has the $(x_1, x_3)$- and the $(x_2, x_3)$-plane as reflective planes of symmetry, and contains the lines $x_1 \pm x_2 = x_3 = 0$.

Figure 2.2 Enneper's surface (left) and the Chen-Gackstatter Surface (right).

The planar symmetry lines cut the surface into four congruent pieces. The symmetry lines themselves cross at $\vec{0}$ and the end point at infinity. In addition, the two straight lines cross at $\vec{0}$ and at the end, and the surface is invariant under rotation about these lines.

We will construct a genus-one surface $S$ with the same symmetries and end behavior as Enneper's surface. From (2.16) we may conclude that the winding $d_1$ at the end of Enneper's surface is equal to three. The two reflections generate a rotation, $R$, of order two about the vertical axis. By (2.16), the total curvature of this surface must be $-8\pi$, which is equivalent to saying that the degree of $g$ is two. By the symmetry assumptions, the Gauss map is vertical at the end. By symmetry again, and the fact that the degree of $g$ is two, all the vertical normals are fixed points of both reflections, which means they lie on the vertical axis. We place the intersection of the two lines at $\vec{0}$. Because rotation about the lines on the surface reverses orientation, any axis point, $p$, of the surface (where the normal must be vertical) other than $\vec{0}$, is paired with $-p$, where $g$ has the same value. There must be at least one such point $p$ on the positive $x_3$-axis and we orient the surface so that $g(p) = \infty$. Since $g(\vec{0})$ is vertical, it follows that $g(0) = 0$ and $g = 0$ at the end. There can be no other points on the surface on the vertical axis where the Gauss map is vertical. We have

$$g(e) = 0, \ g(\vec{0}) = 0, \ g(p) = g(-p) = \infty, \ p = (0, 0, \alpha),$$

where the end is at $e$ on the torus $\overline{S}$.

Choose a triangulation of $\overline{S}$ whose vertices include $0, \pm p$, and $e$. This triangulation projects to a triangulation of $\overline{S}/R$ with half the number of faces and half the number of edges. But the number of vertices is $\frac{V}{2} + 2$, where $V$ is the number of vertices in the triangulation of $\overline{S}$. Since $\overline{S}$ is a torus

$$\chi(\overline{S}/R) = \frac{\chi(\overline{S})}{2} + 2 = 2\,;$$

so $\overline{S}/R$ is a sphere.

Let $z \colon \overline{S} \to \overline{S}/R = S^2$ denote projection onto the quotient. The meromorphic function $z$ has degree two and is branched at the points $\{\vec{0}, \pm p, e\} \subset \overline{S}$. Without loss of generality, we may assume that $z(\vec{0}) = 0$, $z(p) = 1$ and $z(e) = \infty$. Reflection $K$ in either one of the vertical planes of symmetry induces the same order-two automorphism of $S^2$. This automorphism fixes $0, 1, \infty$ and $z(-p)$. We may conclude that this automorphism is complex conjugation $z \to \overline{z}$ and therefore $z(-p)$ is real. Rotation about either line in $S$ induces an orientation reversing automorphism of $S^2$. This automorphism fixes $0$ and $\infty$ while interchanging $\pm p$. Thus rotation induces $z \to -\overline{z}$ on $S^2$. Moreover, $z(-p) = -1$ and the two lines on $S$ are projected by $z$ onto the imaginary axis.

We now have the torus $\overline{S}$ presented as a double covering of $S^2$, branched at $0, \infty$ and $\pm 1$. Therefore, the surface $\overline{S}$ is the Riemann surface

$$w^2 = z(z-1)(z+1)\,.$$

This is the square torus on which $z$ is the Weierstrass $\mathcal{P}$-function, normalized to have a double pole at $e = (\infty, \infty)$, a double zero at $\vec{0} = (0,0)$, (and branch points at $\pm p = (\pm 1, 0)$). Introducing $\gamma = z/w$, we have

$$\gamma^2 = \frac{z}{(z-1)(z+1)}\,. \tag{2.18}$$

The function $\gamma$ has degree equal to two, and has simple zeros at $\vec{0}$ and $e$, and simple poles at $\pm p$. Therefore $g = \rho\gamma$ for some nonzero complex constant $\rho$. The Gauss map is either purely real or purely imaginary along curves of $S$ projecting onto the real axis in the $z$-plane. This is because these curves lie in the assumed symmetry planes $x_2 = 0$ and $x_1 = 0$. Because $\gamma$ is real or imaginary over the real z-axis, it follows that $\rho$ is either real or imaginary. After a rotation of $\mathrm{R}^3$, if necessary, we may assume $\rho$ is real and positive. From Theorem 2.1 we know that the differential $dh$ has simple zeros at $0$ and $\pm p$ and no poles at any other finite points. At the end $e$, it must have a triple pole (producing end behavior like that of Enneper's surface). The one-form $dz$ has precisely these zeros and poles. Hence, $dh = c\,dz$ for some nonzero complex constant $c$, which by scaling we may assume to be unitary. We will now show that $c$ is real. To see this, observe that along lifts of the real axis $\frac{d\gamma}{\gamma}(\frac{\partial}{\partial t})$ is real as is $dz(\frac{\partial}{\partial t})$. But these curves are geodesic curvature lines so, by (2.13), $c$ must be real.

We may therefore choose $c = 1$ and conclude that

$$g = \rho\gamma \qquad dh = dz, \qquad (2.19)$$

where $z$ is the projection defined above and $\rho > 0$ is to be determined.

All that remains is to choose $\rho > 0$ so that (2.8) is satisfied. Clearly, the third component has no periods because $z$ is well-defined on the square torus. Because $e$ lies on the intersection of two vertical symmetry planes, the real period of $X$ must be orthogonal to both symmetry planes, and hence zero, for any closed curve around $e$.

To determine the value of $\rho$ that makes $X$ in (2.6) well-defined, we exploit the symmetry of the situation. Consider the curve on the square torus that is the lift of the interval $[-1, 0]$ in the $z$-plane, on which $\gamma$ is positive. (It is always real). This curve must be in the vertical plane $x_2 = 0$, assuming we start the integration from $(0, 0)$. We require that

$$Re \int_0^{-1} (g^{-1} - g) dh = 0,$$

or

$$0 = Re \int_0^{-1} (\frac{1}{\rho\gamma} - \rho\gamma) dz = \rho^{-1} \int_{-1}^0 \sqrt{\frac{(t^2-1)}{t}} dt - \rho \int_{-1}^0 \sqrt{\frac{t}{(t^2-1)}} dt.$$

Both of the integrals are convergent and positive. Choosing $\rho$ to be the positive square root of the ratio of these integrals insures that this way of evaluating $-p = X((-1, 0))$ places it on the $x_3$-axis. Choice of negative square roots in the integrand will result in the same value of $\rho$, placing $X((1, 0))$ on the vertical axis. Thus $X$ is well-defined on the torus and produces the genus-one Chen-Gackstatter example with the symmetry of Enneper's surface.

**Remark 2.3** *i) The only complete minimal surfaces in $R^3$ with total curvature $-4\pi$ are the Catenoid and Enneper's surface (Osserman [58]). Recently Lopez [49] proved that the unique genus-one complete minimal surface with total curvature $-8\pi$ is the Chen-Gackstatter example. This result was also contained in the thesis of D. Bloss [3].*

*ii) The Chen-Gackstatter genus-one construction has been generalized in two different ways. First, Karcher realized that the kth-order Enneper surface given by the Weierstrass data*

$$g(z) = z^k, dh = g(z)dz, z \in C,$$

*could also be used as a model for the end behavior and symmetry of higher genus examples. Using the Lopez-Ros parameter (see Section 3.1, Equation (3.5)) he was able to show that the periods could be killed for all $k > 0$ [43]. Second, one can add more handles by stacking them. Chen and Gackstatter did this themselves in the genus-two case [11]. Computationally, higher genus examples were found by E. Thayer, who used a computer to solve the period problem (2.8). Thayer also computed higher genus examples (up to 35) with higher winding orders [71]. N. do Espírito-Santo established existence for stacked handles in the genus-three case [19]. (See Remark 4.2.)*

## 2.3 Embeddedness and finite total curvature: necessary conditions

We continue the general discussion of Section 2.1. Suppose $X: M \to \mathbb{R}^3$ has an embedded end at some $p \in \overline{M}_k - M$. Let $D - \{p\}$ be a punctured neighborhood of $p$ on which $X$ is one-to-one. We may assume, after a rotation if necessary, that the Gauss map, which extends to $p$, takes on the value $(0,0,1)$ at $p$. The two simplest examples are the plane and (either end of) the catenoid. Written as graphs over the $(x_1, x_2)$-plane, these ends have bounded and logarithmic growth respectively. This is in fact all that can happen.

**Proposition 2.1** ([59] [67]) *Suppose $X: M \to R^3$ has a complete embedded end $X: D - \{p\} \to R^3$ of finite total curvature with the value of the extended Gauss map equal to $(0,0,\pm 1)$ at $p$. Then, outside of a compact set, $X(D - \{p\})$ is the graph (over the exterior of a bounded domain in the $(x_1, x_2)$-plane) with the following asymptotic behavior:*

$$x_3(x_1, x_2) = \alpha \log \rho + \beta + \rho^{-2}(\gamma_1 x_1 + \gamma_2 x_2) + \mathcal{O}(\rho^{-2}) \tag{2.20}$$

*for $\rho = (x_1^2 + x_2^2)^{1/2}$. Moreover the first two components $\phi_1, \phi_2$ of $\Phi$ in (2.7) have poles of order two at $p$ and have no residues, while the third component $dh = \phi_3$ is either regular (which happens if and only if $\alpha = 0$ in (2.20)) or has a simple pole.*

**Proof.** We will prove this proposition in the case that the Gauss map at the embedded end is one-to-one. In this case we may assume, without loss of generality, that $p = 0$ and $g(z) = z$ on $D$. (In general we may assume $g(z) = z^k$ on $D$. The reader is invited to carry out a similar computation, showing that $\alpha = 0$ when $k > 1$.)

From (2.7) and (2.8) we have

$$x_1(z) = Re \int_{z_0}^z \phi_1 = Re \int_{z_0}^z \frac{1}{2}(g^{-1} - g)dh = Re \int_{z_0}^z \frac{1}{2}(z^{-1} - z)dh$$

$$x_2(z) = Re \int_{z_0}^z \phi_2 = Re \int_{z_0}^z \frac{i}{2}(g^{-1} + g)dh = Re \int_{z_0}^z \frac{i}{2}(z^{-1} + z)dh$$

$$x_3(z) = Re \int_{z_0}^z dh$$

It is straightforward to compute

$$2(x_1 - ix_2) = \int_{z_0}^z g^{-1}dh - \overline{\int_{z_0}^z gdh} = \int_{z_0}^z z^{-1}dh - \overline{\int_{z_0}^z zdh}.$$

We will now look more closely at $dh$. Writing

$$dh = (\sum_{j=1}^k c_{-j}z^{-j} + c_0 + zw_1(z))dz,$$

$w_1(z)$ holomorphic, we observe that $x_1 - ix_2$ is not well-defined if $c_0 \neq 0$. The first integral in the expression for $x_1 - ix_2$ would produce a term $c_0 \ln z$, which is not a well-defined function of $z$. For similar reasons (using the second integral) $c_{-2} = 0$. Knowing that $c_0 = 0$ allows us to conclude that $dh$ has a pole at 0; otherwise the metric (2.9)

$$\frac{1}{2}(|g| + |g|^{-1})|dh| = \frac{1}{2}(|z| + |z|^{-1})|dh| = \frac{1}{2}(|z|^2 + 1)\frac{|dh|}{|z|}$$

is not complete at $|z| = 0$. We claim now that $dh$ has a simple pole at $z = 0$. If not, let $c_{-k}z^{-k}, k \geq 3$, be the highest order term in the expansion. Then the dominant term in $\int z\,dh$ is $-(c_{-k}/k)z^{-k}, k \geq 3$. This means that on $|z| = \varepsilon$ small, $x_1 - ix_2(z)$ has self intersections. However, since we are assuming the end is embedded, this contradicts the fact that the normal is close to vertical near 0, so the projection of the image of a small circle $|z| = \varepsilon$ is one-to-one. Thus, $dh = (c_{-1}z^{-1} + zw_1(z))dz$, $w_1(z)$ holomorphic, as claimed. Because $x_3 = Re \int_{z_0}^{z} dh$ is well-defined, $c_{-1}$ must be real. Moreover, $\phi_1$ and $\phi_2$ have double poles without residues at $z = 0$.

We now can write

$$x_1 - ix_2 = -c_{-1}/z + c + zw_2(z) + \overline{z}\overline{w_3(z)}$$

$$x_3 = const + c_{-1} \ln |z| + \mathcal{O}(|z|^2),$$

where $w_i(z)$ holomorphic and $c$ is a constant. Because the end is complete, $c_1 \neq 0$. We compute

$$\rho^2 = |c_{-1}|^2/|z|^2 - 2Re(\overline{c}_{-1}c/\overline{z}) + const + \mathcal{O}(|z|)$$

and $Re(x_1 - ix_2) \cdot \overline{c}$, which we express as follows:

$$Re(c_{-1}/z)\overline{c} = -Re(x_1 - ix_2) \cdot \overline{c} + |c|^2 + \mathcal{O}(|z|).$$

From the expression for $\rho^2$, one sees that $|z| = \mathcal{O}(\rho^{-1})$. Using this fact, and substitution of the second equation into the expression for $\rho^2$, gives after rearrangement:

$$|c_{-1}|^2/|z|^2 = \rho^2 + 2Re(x_1 - ix_2 \cdot \overline{c}) + const + \mathcal{O}(\rho^{-1})$$

$$= \rho^2(1 + \rho^{-2})(2Re(x_1 - ix_2 \cdot c) + \mathcal{O}(\rho^2)).$$

Taking the logarithm of both sides and using the fact that $\log(1 + a) = a + \mathcal{O}(a^2)$, we get

$$\ln |z| = \ln |c_{-1}| - \ln \rho - \frac{Re(x_1 - ix_2) \cdot c}{\rho^2} + \mathcal{O}(\rho^{-2}).$$

Substitution in the equation above for $x_3(z)$ yields:

$$x_3(z) = const - c_{-1} \ln \rho - \rho^{-2}c_{-1}Re(x_1 - ix_2) \cdot c + \mathcal{O}(\rho^{-2}), \qquad (2.20)'$$

which gives (2.20).

**Definition 2.4** *Consider an embedded end of a complete minimal surface of finite total curvature. The end is said to be a* flat *or* planar *end if* $\alpha = 0$ *in equation (2.20), and is a* catenoid *end* otherwise. *The constant* $\alpha$ *is called the* logarithmic growth *of the end.*

**Remark 2.4** *From (2.20)':*

*i) It is evident that an embedded finite total curvature end is asymptotic to an actual catenoid end ($h = \alpha \log r + \beta$) or to a plane ($h \equiv \beta$);*

*ii) It also follows from equations (2.20)' that the logarithmic growth of a complete embedded end, parameterized on a punctured disk $D - \{p\}$, is equal to $-c_1$, where $c_1$ is the residue of dh at p. If the end is parameterized by the exterior of a disk ($p = \infty$) then the logarithmic growth is equal to the residue of dh at infinity;*

*iii) In the proof of Proposition 2.1 we showed that if the extended Gauss map is unbranched at an embedded end of finite total curvature, that end is a catenoid end. Since the only other geometric possibility is to have a flat end, it follows that branching of the extended Gauss map at an embedded end of finite total curvature is equivalent to the flatness of that end.*

**Proposition 2.2 (Gackstatter [24] Jorge-Meeks [42])**

*i) All of the ends of M are embedded $\Leftrightarrow n(M) = r \Leftrightarrow$ equality holds in (2.14):*

$$\int K dA = -4\pi(k + r - 1). \tag{2.21}$$

*ii) If no two ends of M intersect then, after a rotation, the Gauss map N satisfies $N(p_i) = (0, 0, \pm 1), i = 1, \ldots, r$.*

**Remark 2.5** *If $X \colon M \to R^3$ is an embedding, then the conditions of Proposition 2.2 must hold. However, these conditions are not sufficient to imply that X is an embedding. See the discussion and examples presented in Section 5.*

Proof *i)* follows immediately from Theorem 2.3 *ii)*, formula (2.16) and the observation made after Definition 2.3 that $d_j = 1$ when an end is embedded. Statement *ii)* follows from Remark 2.4, because two ends that are of planar- or catenoid-type must intersect if their limit normals are not parallel.

Because formula (2.21) is important and we did not present the proof of Theorem 2.3 we will give a direct proof of it here. Formula (2.21) is clearly true for the flat plane ($k = 0, r = 1, K \equiv 0$). Hence we may assume $K < 0$ except at isolated points. $M$ is by assumption a compact Riemann surface of genus $k$, from which $r$ points, $p_1, \ldots p_r$, have been removed. Choose closed curves $\widetilde{\beta}_j$ that bound small disks $D_j$ centered at $p_j$. Since the ends are assumed to be embedded, Remark 2.4,i) implies that they are asymptotic to half-catenoids or planes. This means that a sequence of curves homotopic to $\widetilde{\beta}_j$, whose lengths approach a minimum for all such curves, cannot diverge. Because $K < 0$ except at isolated points of $M$, a *unique* length-minimizing geodesic $\beta_j$ exists in the homotopy class of $\widetilde{\beta}_j$ ([27]). Each geodesic bounds an end representative that is topologically an annulus. By intersecting a planar end with a large sphere or a catenoid end with a plane orthogonal to the limit normal and sufficiently far out, we may approximate the annular end $E_j$ by a finite annulus bounded on one side by the geodesic $\beta_j$ and on the other by a closed curve $c$. Moreover the closed curve $c$ can be chosen, according to Proposition 2.1 and formula (2.20), to be of the form $c_\rho(\theta)$, $0 \leq \theta \leq 2\pi$, with

$$c_\rho(\theta) = (\rho\cos\theta, \rho\sin\theta, \alpha\log\rho + \rho^{-2}(\gamma_1\cos\theta + \gamma_2\sin\theta) + \mathcal{O}(\rho^{-2})).$$

Here $\rho > 0$ can be made as large as desired, and $\gamma_1, \gamma_2$ are fixed. Increasing $\rho$ corresponds to pushing $c_\rho(\theta)$ farther out on the end; as $\rho \to \infty$, $c_\rho(\theta)$ diverges. It is straightforward to see from the equation for $c_\rho(\theta)$ that its total curvature converges quickly to $2\pi$. Using the Gauss-Bonnet formula, we have

$$\int_{E_j} KdA \cong 2\pi\chi(E_j) - \int_{\beta_j} k_g ds - \int_c \kappa ds = -\int_c \kappa ds,$$

since the Euler characteristic of an annulus is zero and $\beta_j$ is a geodesic. But since $c$ is asymptotically a circle, we have shown that the total curvature of the end bounded by $\beta_j$ is exactly $-2\pi$.

We next show that the $\beta_j$ are pairwise disjoint. If two touch but do not cross, then of course they must coincide, which implies that $S$ is the union of two ends with a common boundary. The surface $S$ is an annulus ($g=0, r=2$) with total curvature $-4\pi$, so (2.21) is satisfied in this case. (By a result of Osserman [58], the only such surface is the catenoid. See also Lemma 3.1 in Section 3.1.) If two $\beta_j$ cross, then a segment of one together with a segment of the other bounds a region $\Sigma$ of the surface. Since $\Sigma$ is simultaneously inside of two ends it is easy to see that $\Sigma$ is simply connected. This leads to a contradiction of the fact that each $\alpha_j$ is the unique geodesic in its homotopy class. (Alternatively, we could use the Gauss-Bonnet formula again to get a contradiction; the Euler characteristic of $\Sigma$ is 1 and its boundary consists of two geodesics meeting at two vertices whose exterior angles are less than $\pi$ in absolute value:

$$2\pi = 2\pi\chi(\Sigma) = \int_\Sigma KdA + \Sigma \text{ exterior angles } < \int_\Sigma KdA + 2\pi.$$

This clearly contradicts the fact that $K < 0$ almost everywhere on $\Sigma$.) Hence the $\beta_j$ are disjoint.

Now remove the $r$ ends from $S$ by cutting along the geodesics $\alpha_j$ $j = 1, \ldots r$. This leaves a compact surface $\widehat{M}$ of genus $k$, bounded by $r$ closed geodesics. A third (and final) application of Gauss-Bonnet gives

$$2\pi(2 - 2k - r) = 2\pi\chi(\widehat{M}) = \int_{\widehat{M}} KdA.$$

Combining this with the previously established fact that each end has total curvature $-2\pi$ gives
$$\int_M KdA = \int_{\widehat{M}} KdA + r(-2\pi) = -4\pi(k+r-1),$$
which is (2.21). See [9] for more details.

If $S$ is embedded and complete with finite total curvature, then we know by Proposition 2.1 that outside of a sufficiently large compact set of $\mathrm{R}^3$, $S$ is asymptotic to a finite number of half-catenoids and planes, which may be assumed to have the same vertical limit normal (Proposition 2.2 *ii)*). If $S$ is not a plane it must have *at least two* catenoid ends, one with positive, the other with negative logarithmic growth. This follows from the Halfspace Theorem (Theorem 2.4, which will be proved in Section 2.3.3). It states that a complete properly immersed minimal surface without boundary in $\mathrm{R}^3$, which is not a flat plane, cannot lie in any halfspace. If $S$ has only one catenoid end, it would be in a halfspace yet not be a plane, a contradiction.

Since we are assuming that $S$ is connected, the fact that $S$ is properly embedded implies that $\mathrm{R}^3 - S$ consists of precisely two components. Outside of a sufficiently large compact set, the ends of $S$ are stacked and thus ordered from top to bottom. This also means that the limit normals $(0, 0 \pm 1)$ alternate from one end to the next and that the logarithmic growth rates are also ordered: If the ends correspond to points $\{p_1 \ldots p_r\} \subset \overline{M}$, and $\alpha_j$ is the logarithmic growth rate of the $j$th end:

$$\alpha_1 \leq \alpha_2 \cdots \leq \alpha_r, \qquad \alpha_1 \alpha_r < 0. \tag{2.22}$$

### 2.3.1 Flux

On a surface $S$, with boundary, let $\nu$ denote the outward-pointing unit-normal vector field on $\partial S$. For any vector field $W$ that is $C^1$ with compact support on $S$, the divergence theorem states that
$$\int_S (div_S W)dA = \int_{\partial S} \langle W, \nu \rangle ds.$$

If $S \subset \mathrm{R}^3$ and $W$ is a vector field on $\mathrm{R}^3$, we may restrict $W$ to $S$ and consider its tangential component $W^T := W - \langle W, N \rangle N$. Suppose $W$ is a Killing field and $S$ a minimal surface. Then we have $div_S W^T = 0$, which can be seen as follows. For any tangent vector $U$, and any vector field: $W \in \mathrm{R}^3$

$$\nabla_U W^T = [UW^T]^T = [UW - U(\langle W, N \rangle N)]^T$$
$$= [UW]^T - \langle W, N \rangle DN(U).$$

Since $W$ is Killing, $\langle (UW)^T, U \rangle = \langle UW, U \rangle = 0$. Since $S$ is minimal, $tr_S DN = 0$. Hence $div_S W^T = tr_S \nabla W^T = tr_S[DW]^T - \langle W, N \rangle tr_S DN = 0$, as claimed.

Figure 2.3.1 Flux

Flux on a curve is computed by integrating the outward-pointing unit normal vector on a cycle. The total flux of the boundary of a minimal surface is zero.

The divergence theorem yields

$$\int_{\partial S} \langle W, \nu \rangle ds = 0\,,$$

since $\langle W^T, \nu \rangle = \langle W, \nu \rangle$. Thus the one-form $\omega(V) := \langle V, W \rangle$ is a closed form on $S$ and for closed curves $\gamma$, $\int_\gamma \omega$ is an homology invariant.

A constant vector field $T$ is a Killing field associated to translation in that direction. Since

$$\int_\gamma \langle T, \nu \rangle ds = \langle T, \int_\gamma \nu ds \rangle$$

and the left-hand side is a homology invariant for all constant vector fields $T$, so is $\int_\gamma \nu ds$.

**Definition 2.5** *The flux along $\gamma$ is the vector quantity*

$$Flux\,([\gamma]) := \int_\gamma \nu(s) ds\,.$$

We remark that the flux can also be defined by the same equation for any closed curve $\gamma$ contained in $S$, using $\nu(s) = Rot_{\frac{\pi}{2}}(d\gamma/ds)$. This coincides with our definition of flux when $\gamma$ is a boundary curve of $S$ and is easily seen to be homology-invariant.

**Proposition 2.3** *Suppose $X \colon M \to R^3$ is a complete conformal minimal immersion given by $X = \operatorname{Re} \int \Phi$ as in (2.6), (2.7). Suppose $\gamma = X(\widetilde{\gamma}), \widetilde{\gamma} \subset M$ is a closed curve. Let $X^*$ denote the conjugate surface (See Definition 2.7), which is a (possibly multi-valued) conformal immersion $X^* \colon M \to R^3$:*

i) *Flux $([\gamma]) = -\operatorname{Period}_{\widetilde{\gamma}} X^*$;*

ii) *If $\widetilde{\gamma}$ is in the homology class of a puncture, p,*

$$\operatorname{Flux}([\gamma]) = -2\pi \ \operatorname{Residue}_p \Phi \,;$$

iii) *If $\widetilde{\gamma}$ is in the homology class of a puncture p, which represents a vertical embedded end of finite total curvature,*

$$\operatorname{Flux}([\gamma]) = (0, 0, 2\pi\alpha) \,, \tag{2.23}$$

*where $\alpha$ is the logarithmic growth of the end. In particular the flux of a homology class representing a finite-total-curvature end is always vertical and is zero if and only if the end is flat;*

iv) *If $X(M)$ is an embedding with finite total curvature and logarithmic growth rates $\alpha_1, \ldots \alpha_r$ at the ends, then*

$$\sum_1^r \alpha_j = 0 \,. \tag{2.24}$$

**Proof.** A constant vector field $T$ is the gradient, in $\mathrm{R}^3$ of a linear function $f$. The tangential component $T^T$ of $T$ restricted to $\mathcal{S} := X(M)$ is the (Riemannian) gradient of $f|_S$, which is harmonic, according to (2.5), because $X$ is a minimal immersion. Moreover, $Rot_{\frac{\pi}{2}} T^T = \operatorname{grad} f^*$, where $f^*$ is the locally defined harmonic conjugate of $f$. (See the discussion between (2.3) and (2.4).) Now we compute

$$T \cdot \operatorname{Flux}([\gamma]) = \int_\gamma \langle T^T, \nu \rangle ds = \int_\gamma \langle Rot_{\frac{\pi}{2}} T^T, Rot_{\frac{\pi}{2}} \nu \rangle ds$$

$$= \int_\gamma \langle \operatorname{grad} f^*, \frac{-d\gamma}{ds} \rangle ds = -\int_\gamma df^* \,.$$

Since this is true for *any* constant vector field $T$ in $\mathrm{R}^3$, it follows that Flux $([\gamma]) = -$ Period $_{\widetilde{\gamma}} X^*$, which is i). From this, ii) follows immediately since

$$\int_{\widetilde{\gamma}} \Phi = \operatorname{Re} \int_{\widetilde{\gamma}} \Phi + i\operatorname{Im} \int_{\widetilde{\gamma}} \Phi = \operatorname{Period}_{\widetilde{\gamma}} X + i\operatorname{Period}_{\widetilde{\gamma}} X^* = -i \ \operatorname{Flux}([\gamma]) \,. \tag{2.25}$$

To establish iii) we use Remark 2.4ii), together with statement ii) above and the fact that the third component of $\Phi$ is $dh$. We have that the third component of Flux $([\gamma])$, where $[\gamma]$ is the homotopy class of a puncture $p$, must equal $-2\pi \operatorname{Residue}_p dh = 2\pi\alpha$, where $\alpha$ is the logarithmic growth of the end at $p$. According to Proposition 2.1, the first two components of $\Phi$ have no residues at $p$.

If $X(M)$ has finite total curvature then according to Theorem 2.2iii), $dh$ is globally defined, with poles only at the ends $p_1 \ldots p_r$. Hence, using statement iii)

$$0 = \sum_{i=1}^{r} \operatorname{Residue}_{p_i} dh = 2\pi (\sum_{i=1}^{r} \alpha_i) \,.$$

This establishes statement iv). □

**Remark 2.6** *If we consider the normal component of the Killing field, $W^N =: w^N N$, where $w^N = \langle W, N \rangle$, a calculation shows that*

$$-\Delta_S w^N + 2K w^N = 0 \,,$$

*where $K$ is the Gauss curvature of $S$. That is, $w^N N$ is a Jacobi-field. From a geometric point of view this is expected, since a Killing field is an infinitesimal isometry, and the Jacobi equation is derived from the formula for the second variation of area of a minimal surface. This is considered in more detail in Section 7.*

### 2.3.2 Torque

Let $R_{\vec{u}}$ be the Killing field associated with counter-clockwise rotation about the axis $\ell_{\vec{u}}$ in the $\vec{u}$ direction. From the identity $(U \wedge V) \cdot W = det(U, V, W)$, for vectors $U, V, W$ in $\mathrm{R}^3$, we have

$$(X \wedge \nu) \cdot \vec{u} = (\vec{u} \wedge X) \cdot \nu = R_{\vec{u}} \cdot \nu \,, \tag{2.26}$$

where $X$ is the position vector of a minimal surface $S$ and $\nu$ is the outward-pointing normal to a component $\gamma$ of $\partial S$. Because $R_{\vec{u}}$ is a Killing field, $\int_\gamma R_{\vec{u}} \cdot \nu$ is a homology invariant and $\int_{\partial S} R_{\vec{u}} \cdot \nu = 0$. This motivates defining the *torque* of a closed curve $\gamma$ on $S$ as the vector-valued quantity.

**Definition 2.6**

$$Torque_0(\gamma) = \int_\gamma X \wedge \nu \,.$$

Torque was introduced by Kusner in [45, 46]. From (2.26) it follows that the component of torque in the $\vec{u}$ direction is $\int_\gamma R_{\vec{u}} \cdot \nu$. If we move the origin from $0$ to $W \in \mathrm{R}^3$ and let $\widehat{X}$ be the position vector measured from $W$, then $\widehat{X} = X - W$ and we can compute

$$\begin{aligned}Torque_W(\gamma) &= Torque_0(\gamma) - W \wedge \int_\gamma \nu \\ &= Torque_0(\gamma) - W \wedge Flux(\gamma) \,.\end{aligned} \tag{2.27}$$

It follows that if $Flux(\gamma) = 0$, the torque of $\gamma$ does not depend upon the base-point of $X$. It also follows that the torque does not change if we move the base point parallel to the flux. If $E$ is an embedded end we may speak of the torque of $E$; it is the torque of a closed curve in the homology class of the ends. From Proposition 2.1 the following geometric proposition can be deduced.

**Proposition 2.4** *If $E$ is a vertical catenoid end the torque of $E$ is a horizontal vector. If $E$ is a vertical flat end at which the degree of the Gauss map is at least three, the torque of $E$ vanishes.*

In the case of a vertical catenoid end, we know from Proposition 2.3iii) that the flux is a vertical vector of the form $(0, 0, 2\pi\alpha)$ where $\alpha$ is the logarithmic growth. We can conclude from (2.27) that there is a unique horizontal vector $W = (w_1, w_2, 0)$ for which the *Torque* of $E$ is vertical. From Proposition 2.4, it follows that there is a unique vertical line on which the torque vanishes.

**Definition 2.7** *For a vertical catenoid end $E$, the (unique) vertical line $\ell_E$ on which the torque vanishes is called the* axis *of $E$.*

The torque of $E$ measured from a base point displaced by $W$ from $\ell_E$ is

$$\begin{aligned} Torque_W(E) &= W \wedge (0, 0, 2\pi\alpha) \\ &= 2\pi\alpha(w_2, -w_1, 0) \,. \end{aligned} \tag{2.28}$$

In the case of an embedded finite total curvature minimal surface with catenoid ends $E_1 \ldots E_n$, with logarithmic growth rates $\alpha_i$ and either no flat ends or all flat ends of order three or greater, the total torque of all the ends must sum to zero. For any base point let $W_i$ be the horizontal displacement of $\ell_{E_i}, i = 1 \ldots n$. Then

$$\begin{aligned} 0 &= \Sigma \, Torque_{-W_i}(E_i) = 2\pi Rot_{\frac{\pi}{2}}(\Sigma_{i=1}^n \alpha_i W_i), \quad \text{or} \\ \Sigma \alpha_i W_i &= 0 \,. \end{aligned} \tag{2.29}$$

**Remark 2.7** *If a surface has only two catenoid ends $E_1, E_2$, and either no flat ends or all flat ends of order three or greater, then placing the base point on $\ell_{E_1}$ forces $W_1 = 0$ and by (2.29) we must have $W_2 = 0$. That is, the catenoid ends $E_1$ and $E_2$ have the* same *axis. Pascal Romon [personal communication] points out that if there are three catenoid ends (and the same hypothesis about flat ends) their axes are coplanar. We may assume, without loss of generality that there are two ends with positive logarithmic growth. If we label these ends $E_1$ and $E_2$ and the third end $E_3$ and choose a base point on $E_3$, then it follows from (2.29) that $W_1 = -\frac{\alpha_2}{\alpha_1} W_2$. Thus the axes $\ell_{E_1}, \ell_{E_2}$ and $\ell_{E_3}$ are coplanar and, because $\frac{\alpha_2}{\alpha_1} < 0$, the axis $\ell_{E_3}$ with negative logarithmic growth lies between the two with positive growth. We note that on all known embedded examples, flat ends of order two do not occur. It is not known if they can occur.*

### 2.3.3 The Halfspace Theorem

We present here a proof of

**Theorem 2.4 (The Halfspace Theorem for Minimal Surfaces, [40])** *A complete, properly immersed, nonplanar minimal surface in $R^3$ is not contained in any halfspace.*

Our proof requires the use of the Maximum Principle for Minimal Surfaces: If $S_1$ and $S_2$ are two connected minimal surfaces with a point $p$ in common, near which $S_1$ lies on one side of $S_2$, then a neighborhood of $p$ in $S_1$ coincides with a neighborhood of $p$ in $S_2$. This implies immediately that the analytic continuations of $S_1$ and $S_2$ coincide. In particular, if $S_1$ and $S_2$ are both complete, $S_1 = S_2$. We refer the reader to [17] for proofs. In the simple case that $S_2$ is a plane, the Maximum Principle states that a minimal surface $S = S_1$ that lies locally on one side of its tangent plane $T_pS = S_2$ must in fact be planar. This is easy to prove from the Weierstrass representation, Theorem 2.1. Without loss of generality, we may assume that $T_pS = \{x_3 = 0\}$ and that $z$ is a local coordinate chart chosen so that $X(0) = p$. The Gauss map of $S$ is vertical at $p$, so that $g(0) = \sigma \circ N(p) = 0, \infty$. Either $g$ is constant, in which case $S$ is flat and is a subset of $T_pS = \{x_3 = 0\}$, or $0$ is an isolated pole or zero of $g$ of order $k \geq 1$. In the latter case,

$$dh = z^k w(z) dz,$$

where $w(z)$ is holomorphic and $c = w(0) \neq 0$. Then

$$x_3(z) = Re \int_0^z dh = Re \int_0^z z^k w(z) dz = Re(\frac{cz^{k+1}}{k+1}) + \mathcal{O}(|z|^{k+2}).$$

Clearly $x_3(z)$ changes sign $2(k + 1)$ times in any small neighborhood of $z = 0$. Hence $S$ does not lie on one side of $T_pS = \{x_3 = 0\}$ near $p$.

**Proof of Theorem 2.4.** By assumption, $S$ is contained in some halfspace $H$. The intersection, $H^*$, of the closed halfspaces containing $S$ and having boundary parallel to $\partial H$ is the smallest closed halfspace containing $S$. We will rotate and translate $S$ so that $H^* = \{x_3 \geq 0\}$, and hence $\partial H^* = \{x_3 = 0\}$. If $S \cap \{x_3 = 0\} \neq \emptyset$, then the simple case of the Maximum Principle implies $S$ is a plane. If $S$ is not a plane $S \cap \{x_3 = 0\} = \emptyset$, and since $S$ is properly immersed, $S$ has no limit points in $\{x_3 = 0\}$. This means each point of $\{x_3 = 0\}$ has positive distance from $S$.

Consider a catenoid with axis parallel to the $x_3$-axis and waist circle in $\{x_3 = 0\}$. Let $C$ denote the half-catenoid in $\{x_3 \leq 0\}$. Because $x_3$ is proper (and unbounded) on $C$ it follows that we can vertically translate $C$ so that the waist circle is in $\{x_3 > 0\}$ but the *solid* half-catenoid inside $C$ is disjoint from $S$. For convenience in the rest of the proof, we move our coordinate system so that $\partial H^* = \{x_3 = -\varepsilon\}$ and the waist-circle of $C$ lies in $\{x_3 = 0\}$. Also, we will now think of $C$ as the boundary for the solid half-catenoid; i.e., we add to $C$ the disk in $\{x_3 = 0\}$ bounded by the waist-circle.

Figure 2.3.3 The vertically translated half-catenoid $C$
in the proof of Theorem 2.4.

Let $tC$ denote the rescaling of $C$ by a factor of $t > 0$. We will say that *$tC$ lies below $S$* if $S$ lies in the closure of the component of $\mathbb{R}^3$ above $tC$. (Similarly we may define what it means for a plane to lie below $S$.) For $t \leq 1$, the waist-circles of $tC$ lie on $C$, so their distance from $S$ is uniformly bounded away from 0. If $S \cap tC \neq \emptyset$, then the Maximum Principle implies $S = tC$, and $S$, being complete, must be a catenoid; this is a contradiction, since $S$ lies in $H^*$. We may conclude that if $tC$ is below $S$, then $tC$ and $S$ have no points in common. Furthermore, since $tC \cap H^*$ is compact and $S \subset H^*$, it follows that every $tC$ lying below $S$ actually has finite positive distance from $S$. In particular the set

$$T := \{t \in (0,1] : tC \text{ is below } S\}$$

is open in $(0,1]$. But if $\{t_n\} \subset T$ and $\lim_{n \to \infty} t_n = \tau > 0$, then $\tau C$ is a half catenoid below $S$; if not, some $t_n C$ must not be below $S$. Hence $T$ is closed, too. This means that $T = (0,1]$, so all $tC$ are below $S$, $0 < t \leq 1$.

Note that the half-catenoids $tC$ converge, as $t \to 0$, to the plane $\{x_3 = 0\}$, which is strictly inside $H^* = \{x_3 \geq -\varepsilon\}$. Since all the $tC$ are below $S$, the limit plane $\{x_3 = 0\}$ is also below $S$. This contradicts the definition of $H^*$ as the intersection of all the closed halfspaces containing $S$ and having boundary parallel to $\partial H$. This contradiction completes the proof. □

**Remark 2.8** *A slight modification of the same proof can be used to prove the following slight improvement of Theorem 2.4: Suppose $S$ is a properly immersed nonplanar minimal surface with compact boundary $\partial S$. If $H$ is a halfspace containing $S$, and $H^*$ a smallest halfspace containing $S$, with $\partial H^*$ parallel to $\partial H$, then $\partial H^* \cap \partial S \neq \emptyset$.*

## 2.4 Summary of the necessary conditions for existence of complete embedded minimal surfaces with finite total curvature

We conclude this section by gathering together the results of our discussion to formulate a sharpened version of Theorem 2.2 for *embedded* complete, nonplanar minimal surfaces of finite total curvature.

**Proposition 2.5** *Let $X: M \to R^3$ be a complete, nonplanar, conformal minimal embedding with finite total curvature. Then:*

i) *$M$ is conformally diffeomorphic to $\overline{M}_k - \{p_1, \ldots, p_r\}$ where $\overline{M}_k$ is a closed Riemann surface of genus $k$ and $p_1, \ldots, p_r$ are points in $\overline{M}_k, r \geq 2$;*

ii) *$X$ is proper. A punctured neighborhood of each $p_i$ is mapped by $X$ onto an end of $S = X(M)$ that is asymptotic to either a plane or a half-catenoid.*

iii) *The Gauss map $N: M \to S^2$ extends to a meromorphic function on $\overline{M}_k$. All the normals on $\{p_1, \ldots, p_r\}$ must be parallel and after a rotation if necessary we may assume that $N(p_i) = (0, 0 \pm 1)$, $i = 1, \ldots r$. The holomorphic one form $dh$, the complex differential of the height function, extends meromorphically to $\overline{M}_k$.*

iv) *The total curvature of $M$ is $-4\pi \cdot$ degree $g = -4\pi(k + r - 1)$, where $g = \sigma \circ N$, and $\sigma: S^2 \to C \cup \{\infty\}$ is stereographic projection.*

v) *The ends are naturally ordered by height from top to bottom. In this ordering, say $p_1 \ldots p_r$, the unit vertical normals at the ends alternate and the logarithmic growth rates $\alpha_j$ are ordered from biggest to smallest. Furthermore,*

$$\sum_{j=1}^{r} \alpha_j = 0. \tag{2.30}$$

*If $S$ is not a plane, $r \geq 2$ and $\alpha_1 \cdot \alpha_r < 0$.*

vi) *At a catenoid end, $g$ has a simple pole or zero, while at a flat end the pole or zero of $g$, has higher order. The height differential $dh$ has a simple pole at a catenoid end. At a planar end where $g$ has a zero or pole of order $m$, $dh$ has a zero of order $|m| - 2 \geq 0$. The logarithmic growth rate is equal to minus the residue of $dh$ at the puncture corresponding to the end.*

To construct a complete, embedded minimal surface of finite total curvature with the Enneper-Weierstrass-Riemann representation ((2.6) (2.7) of Theorem 2.1), Proposition 2.5 dictates necessary conditions on the choice of $\overline{M}_k$, $\{p_1 \ldots p_r\} \subset \overline{M}_k$, $g = \sigma \circ N$ and $dh$.

However, in order for (2.7) to be single-valued on $M_k = \overline{M}_k - \{p_1 \ldots p_r\}$, it is necessary that (2.8) be satisfied:

$$\text{Period}_\alpha(\Phi) = Re \int_\alpha \left( (g^{-1} - g)\frac{dh}{2}, i(g^{-1} + g)\frac{dh}{2}, dh \right) = 0 \qquad (2.31)$$

for all closed curves $\alpha \subset M_k$. (See Figure 3.0.0.) If Statement iii) of Proposition 2.5 is satisfied, the resulting surface must have embedded parallel ends, which implies that outside of some compact set in $\mathrm{R}^3$, $X(M_k)$ is embedded. However, $X$ may not be an embedding. See the examples in Section 5.

## 3 Examples with restricted topology: existence and rigidity

In this section, we present the few uniqueness results that are presently known. They take the following form: if the topological type is restricted or some other geometric property is specified, then the surface is uniquely determined. (See also Theorem 3.4, below.)

**Theorem 3.1** *Let $S$ be a complete, embedded minimal surface of finite total curvature.*

1) *If $S$ has one end, it is a flat plane;*

2) *(The Lopez-Ros punctured sphere Theorem) [51]. If $S$ has genus zero, it is the catenoid or the flat plane;*

3) *(The R. Schoen catenoid characterization) [67]. If $S$ has two ends, it is the catenoid;*

4) *If $S$ has three ends and a symmetry group of order at least $4(k+1)$, where $k$ is the genus, then $S$ is the surface described in Theorem 3.2.*

5) *(Costa's thrice-punctured torus Theorem) [16]. If $S$ has genus one and three ends, it is one of the surfaces $M_{1,x}$ in Theorem 3.3.*

Statement 1) follows from Theorem 2.4 or Proposition 2.5 v). Statement 4) is proved in [36]. We will give a full proof of 2) in Section 3.1.

**Corollary 3.1** *The only complete, embedded minimal surfaces with $\int K dA \geq -8\pi$ are the plane and the catenoid (with total curvature 0 and $-4\pi$ respectively).*

**Proof.** From Proposition 2.2, it follows that the total curvature is an integer multiple of $-4\pi$. If the total curvature is zero, we must have the plane. If it is $-4\pi$, then by a result of Osserman [58] it must be the catenoid. (See also Lemma 3.1 in Section 3.1.)

We need to show that total curvature $-8\pi$ is not possible. In that case we would have, by Proposition 2.2, i), that $k + r = 3$. But the various cases of Theorem 3.1 exclude this possibility. □

Figure 3.0.0 Attempts to produce nonexistent minimal surfaces.

Left: According to the Lopez-Ros punctured sphere theorem (Theorem 3.1,2)), a complete embedded minimal surface with genus zero and finite total curvature is either the plane or the catenoid. Nonetheless, it is possible to write down Weierstrass data that would produce such an example, provided the period problem (2.8) were solvable (which it is not). One can do this with a high degree of symmetry. Under the assumptions that the desired example has three parallel vertical ends –top and bottom catenoidal, and the middle one flat– one vertical plane of symmetry and a horizontal line diverging into the flat end, the resulting Weierstrass data is

$$g = r(z-1)(z+1)^{-1}, \quad dh = dz$$

on $C - \{\pm 1\}$ with $r > 0$. The gap in the image indicates the failure of the Weierstrass-Enneper mapping to be single-valued. The value of $r$ chosen here is 1.

Right: According to the Schoen catenoid characterization (Theorem 3.1,3)) the catenoid is the only complete embedded minimal surface with two ends and finite total curvature. Nevertheless, we can find symmetric Weierstrass data that fails only in that the period problem is not solvable. Under the assumptions that the desired example has genus one, two vertical planes of reflective symmetry and one horizontal plane of reflective symmetry, the Weierstrass data is $\{g, dh\}$, where

$$g^2 = (1-z)(r-z)((1+z)(r+z))^{-1}, \quad dh = (z^2-1)^{-1}dz,$$

$r > 1$, on the rectangular torus determined by the equation relating $g$ and $z$. The image looks like a catenoid through which someone has tried to drill a tunnel from both sides. The tunnels don't meet up as one sees here (for $r = 2.5$). By (Theorem 3.1,3)), it is futile to try other values of $r$.

**Proposition 3.1** *Let $X\colon M \to R^3$ be a connected complete minimal immersion with finite total curvature. Let $n(M)$ be as in Definition 2.3. Then for any $p \in R^3$, the number of points in $X^{-1}(p)$ is at most $n(M) - 1$, with the sole exception of the case when $X(M)$ is the flat plane.*

**Corollary 3.2** *Let $X\colon M \to R^3$ be as in Theorem 3.1. If $n(M) = 2$, $X(M)$ is the catenoid.*

**Proof of Proposition 3.1.** We will use a version of the monotonicity formula to prove this proposition: *If $p$ is a point on a complete immersed minimal surface and $A(r)$ is the area of the minimal surface inside a Euclidean ball of radius $r$, then $A(r)/\pi r^2$ is a nondecreasing function of $r$, which is strictly increasing unless the minimal surface is a collection of planes.*

Suppose that $p \in X(M)$ and $X^{-1}(p)$ contains $n$ points. Assume that $X(M)$ is not a plane. (It cannot be a collection of planes because we assume that $X(M)$ is connected.) Hence, for small $r$, $A(r)/\pi r^2 > n$. But as a consequence of Theorem 2.3, $A(r)/\pi r^2$ converges to $n(M)$ as $r \to \infty$. From monotonicity we conclude that $n < n(M)$. If $X(M)$ is a flat plane $n = 1 = n(M)$. $\square$

**Proof of Corollary 3.2.** It follows from Proposition 3.1 that $X$ is an embedding and since $n(M) = 2$, it is not the plane. It must therefore have two embedded ends. From Theorem 3.1. 3, it must be the catenoid. $\square$

According to Corollary 3.1, the next smallest possible total curvature for an embedded example is $-12\pi$; that is, $k + r = 4$. By Theorem 3.1, this can only happen when $k = 1$ and $r = 3$, and in fact it does.

**Theorem 3.2 ([35, 36])** *For every $k \geq 2$, there exists a complete properly embedded minimal surface of genus $k-1$ with three annular ends. After suitable rotation and translation, the example of genus $k-1$, which we will call $M_k$, has the following properties:*

1. *$M_k$ has one flat end between its top and bottom catenoid ends. The flat end is asymptotic to the $(x_1, x_2)$-plane.*

2. *$M_k$ intersects the $(x_1, x_2)$-plane in $k$ straight lines, which meet at equal angles at the origin. Removal of the $k$ lines disconnects $M_k$. What remains is, topologically, the union of two open annuli;*

3. *The intersection of $M_k$ with any plane parallel (but not equal) to the $(x_1, x_2)$-plane is a single Jordan curve;*

4. *The symmetry group of $M_k$ is the dihedral group with $4k$ elements generated by reflection in $k$ vertical planes of symmetry meeting in the $x_3$-axis, and rotation about one of the lines on the surface in the $(x_1, x_2)$-plane. Removal of the intersection of $M_k$ with the $(x_1, x_2)$-plane and the vertical symmetry planes disconnects $M_k$ into $4k$ congruent pieces, each a graph.*

Figure 3.0.1 Costa's surface (upper left) and a few members of the family
of her deformations, as described in Theorem 3.3.

5. $M_k$ is the unique properly embedded minimal surface of genus $k-1$ with three ends, finite total curvature, and a symmetry group containing $4k$ or more elements.

As will follow from our presentation in Section 4, the Riemann surface $\overline{M}_k$ is given by

$$w^k = -(\frac{1}{2})z^{k-1}(z-1)(z+1). \tag{3.1}$$

The catenoid ends are located at $(z,w) = (\pm 1, 0)$, and the flat end is at $(\infty, \infty)$. The Gauss map and the complex differential of the height function are given by

$$g = \rho w, \quad dh = \frac{dz}{(z-1)(z+1)} \tag{3.2}$$

where $\rho$ is a constant determined by the necessity of satisfying the period condition (2.6) in Proposition 2.5. The symmetry of the surface is used to show that there is only one period condition. (See Figure 3.0.2., left-hand column, for pictures of the $M_k$.)

**Remark 3.1** *The example with $k = 2$ was found by Celso Costa in 1982. He proved it was complete and had embedded ends of the specified type. See Costa [14, 16] and Hoffman [28, 29].*

Figure 3.0.2 Higher-genus embedded minimal surfaces.

The surfaces in the column on the left are the surfaces $M_k, k = 3, 4, 6$ described in Theorem 3.2, with genus 2 (top), 3 (middle) and 5 (bottom). The middle end in each of these surfaces is flat. In each row, the other two surfaces are deformations of the leftmost one. These are the surfaces described in Theorem 3.3. Their middle ends are catenoidal.

The surfaces $M_k$ each lie in a family of embedded minimal surfaces.

**Theorem 3.3 ([30]. See Section 4)** *For every $k \geq 2$, there exists a one parameter, $M_{k,x}, x \geq 1$, of embedded minimal surfaces of genus $k - 1$ and finite total curvature. The surfaces $M_{k,1}$ are precisely the surfaces $M_k$ of Theorem 3.2. The surfaces $M_{k,x}, x > 1$ have all three ends of catenoid type and a symmetry group generated by $k$ vertical planes of reflectional symmetry. The Riemann surface $\overline{M}_{k,x}$ is given by*

$$w^k = -c_x z^{k-1}(z-x)(z+x^{-1}), \quad c_x = (x+x^{-1})^{-1}. \tag{3.3}$$

*The catenoid ends are located at $(z, w) = (x, 0), (-x^{-1}, 0)$ and $(\infty, \infty)$. The Gauss map and differential of the height function are given by*

$$g = \frac{\rho w}{mz+1}, \quad dh = \frac{z(m+z^{-1})dz}{(z-x)(z+x^{-1})}, \tag{3.4}$$

*where $\rho$ and $m$ are constants determined by the period conditions. (When $x = 1$, $m(1) = 0$, and (3.3) and (3.4) give (3.1) and (3.2) for the surfaces $M_k = M_{k,1}$.)*

**Remark 3.2** *The surfaces $M_{k,1}$ of Theorem 3.3 are the surfaces $M_k$ of Theorem 3.2. In Section 4, we find it more convenient to work with the function $u = z/w$ instead of $w$. The formulae (3.3) and (3.4) are presented (e.g. in (4.12)) in terms of $u$ and $z$. See (4.4) where we perform the conversion, into a $(w, z)$ expression, of the Riemann surface equation (3.3).*

In Section 5, we will discuss examples with more than three ends, some other results about the structure of the "space" of embedded finite total curvature examples, and present some questions and conjectures. In Section 4 we describe in detail the ideas behind the construction of the surfaces $M_{k,x}$ of Theorem 3.3. As we pointed out in the introduction, these are the only higher-genus examples that have been fully analyzed.

## 3.1 Complete embedded minimal surfaces of finite total curvature and genus zero: the Lopez-Ros theorem

In this section we prove the second statement of Theorem 3.1, which is Theorem 3.5 below. Recall that for a closed curve $\gamma \subset S$ its flux is defined in Section 2.3.1 as

$$\text{Flux}\,([\gamma]) = \int_\gamma Rot_{\frac{\pi}{2}}(\frac{d\gamma}{ds})ds\,.$$

Flux $([\gamma])$ depends only on the homology class of $\gamma$. We say that $S$ *has vertical flux* provided Flux $([\gamma])$ is a vertical vector for all closed $\gamma \subset S$.

Our presentation follows that of Perez and Ros [60], whose approach yields the following theorem.

**Theorem 3.4** *Let $S$ be a complete, embedded minimal surface with finite total curvature and vertical flux. Then $S$ is the catenoid or the plane.*

As usual, we may realize $S$ by a conformal minimal embedding $X: M \to \mathrm{R}^3$, where $M = \overline{M} - \{p_1 \ldots p_r\}$, and rotate $S$ in $\mathrm{R}^3$ so that $g = 0$, or $\infty$ at the ends $\{p_1 \ldots p_r\}$.

If $S$ is complete and embedded, with genus zero, the only closed curves one has to consider are those associated to the ends $\{p_1 \ldots p_r\}$. In Proposition 2.3ii) we established that the flux of the homology class of $p_j$ is equal to the vertical vector $(0, 0, 2\pi\alpha_j)$, where $\alpha_j$ is the logarithmic growth of the end at $p_j$. Hence we have as an immediate consequence of Theorem 3.4:

**Theorem 3.5 ([51])** *Let $S$ be a complete minimal surface with genus zero and finite total curvature. Then $S$ is the catenoid or the plane.*

The proof of Theorem 3.4 exploits what has come to be called the Lopez-Ros deformation. If $\{g, dh\}$ is the Weierstrass data of $S$, we define on $\overline{M}$ the data

$$g_\lambda = \lambda g \qquad dh_\lambda = dh, \tag{3.5}$$

for any $\lambda > 0$. Notice that the zeros and poles of $\{g_\lambda, dh_\lambda\}$ are the same as those of $\{g, dh\}$. The deformation we consider is given by the Weierstrass representation (2.6), (2.7).

$$X_\lambda := Re \int \Phi_\lambda := Re \int (\frac{1}{2}(\lambda^{-1}g^{-1} - \lambda g), \frac{i}{2}(\lambda^{-1}g^{-1} + \lambda g), 1) dh. \tag{3.6}$$

Note that $X_\lambda$ may be multi-valued. The metric and curvature of $X_\lambda$ are given by (2.9) (2.10) as

$$ds_\lambda = (\lambda^{-1}|g|^{-1} + \lambda|g|)|dh|$$

$$K_\lambda = \frac{-16}{(\lambda^{-1}|g|^{-1} + \lambda|g|)^4} \frac{|dg|}{|dh|}. \tag{3.7}$$

We will prove three propositions about the Lopez-Ros deformation that allow us to prove Theorem 3.4. The propositions are stated here and then the proof of Theorem 3.4 is presented. The proofs of the propositions follow.

**Proposition 3.2** *Let $X: M \to R^3$ be a conformal minimal immersion:*

i) *$X$ is complete if and only if $X_\lambda$ is complete for all $\lambda > 0$;*

ii) *If $X_\lambda$ is single-valued, the total curvature of $X$ is finite if and only if the total curvature of $X_\lambda$ is finite for all $\lambda > 0$;*

iii) *The immersions $X_\lambda$ are single-valued for all $\lambda > 0$ if and only if $X$ has vertical flux.*

**Proposition 3.3** *Let $X: M \to R^3$ be a complete, embedded minimal surface of finite total curvature, for which $X_\lambda$ is single-valued for all $\lambda > 0$. Then $X_\lambda$ is an embedding for all $\lambda > 0$.*

**Proposition 3.4** *Suppose $X: M \to R^3$ is a conformal minimal immersion.*

i) *If for some $p \in M, N(p) = (0, 0 \pm 1)$ (i.e., $g(p) = 0$ or $g(p) = \infty$) then, on every neighborhood of $p$, $X_\lambda$ is not an embedding, for $\lambda$ sufficiently large.*

ii) *Suppose $X: M = \{z| \ 0 < |z| < \varepsilon\} \to R^3$ is a conformal embedding representing a flat end. Then $X_\lambda$ is not an embedding, for $\lambda$ sufficiently large, unless $X$ represents the end of a flat plane.*

**Proof of Theorem 3.4.** Let $X\colon M \to \mathbb{R}^3$ be a conformal minimal embedding with $S = X(M)$. The vertical flux condition guarantees that the Lopez-Ros deformation $X_\lambda$, defined by (3.5) and (3.6), is single-valued for all $\lambda > 0$ (Proposition 3.2 above). Because $X$ is a complete embedding with finite total curvature, the same is true for $X_\lambda, \lambda > 0$ (Propositions 3.2, and 3.3 above). But this means that $S$ has no vertical points (points where $g = 0, \infty$) and no flat ends, unless $S$ is itself a flat plane (Proposition 3.4 above).

Assume that $S$ is not a plane.

Because $S$ is embedded and has finite total curvature, $M = \overline{M} - \{p_1 \ldots p_r\}$, where $\overline{M}$ is compact. The ends at $p \in \{p_1, \ldots p_r\}$ are either planar or catenoid-type. The zeros of $dh$ can occur only at vertical points of $g$ or at flat ends, while $dh$ has simple poles at the catenoid ends. (See Proposition 2.5 in Section 2.3.4.) From the previous paragraph we can conclude that $dh$ has no zeros. (Because $S$ is not a plane, $dh \not\equiv 0$.) Thus $\chi(\overline{M}) = 2$, i.e. $\overline{M}$ is the sphere and $r = 2$. We may conclude from Lemma 3.1 that $S$ is the catenoid. □

**Lemma 3.1** *Suppose $S$ is a complete, embedded minimal surface with total curvature $-4\pi$ and two ends. Then $S$ is a catenoid.*

This follows from (2.21). Proposition 2.5i) and Theorem 3.1,3). It also follows from Osserman's characterization of complete minimal suraces of total curvature $-4\pi$ [58]. We give a direct proof here in order to make this section as self-contained as possible.

**Proof of Lemma 3.1.** Because $g = 0, \infty$ at the ends, and the values must alternate (see Proposition 2.5), we may assume, without loss of generality, that $\overline{M} - \{p_1, p_2\} = \mathbb{C} - \{0\}$ and $g(0) = 0, g(\infty) = \infty$. Since $g$ has no other poles or zeros, it follows that $g(z) = az, a \neq 0$. Without loss of generality, we may assume that $a$ is real, because rotation of $S$ about a vertical axis changes $g$ by a unitary multiplicative factor. Because $dh$ has a simple pole at $0$ and at $\infty$ and no other zeros or poles, $dh = cdz/z$ for some $c \neq 0$. Letting $\zeta = az$ we have $g(\zeta) = \zeta$ and $dh = cd\zeta/\zeta$. Clearly, $dh$ has a real period on $|\zeta| = 1$ if $c$ is not real. Thus our surface is the catenoid as described in Section 2. (See Figure 2.0.) □

We now present proofs of Propositions.

**Proof of Proposition 3.2.** From (3.7) it is evident that for each $\lambda > 0$

$$\underline{c}\, ds \leq ds_\lambda \leq \overline{c}\, ds$$

$$\overline{c}^{\,-4}|K| \leq |K_\lambda| \leq \underline{c}^{\,-4}|K|,$$

where $\overline{c} = \max\{\lambda^{-1}, \lambda\}$ and $\underline{c} = \min\{\lambda^{-1}, \lambda\}$. From these inequalities, statements i) and ii) follow directly. To prove statement iii), we use (2.25) from Section 2.3.1. For any closed curve $\gamma = X(\widetilde{\gamma}), \widetilde{\gamma} \subset M$, we have

$$\int_{\widetilde{\gamma}} \Phi = Re \int_{\widetilde{\gamma}} \Phi + i\, Im \int_{\widetilde{\gamma}} \Phi = \text{Period}_\gamma(X) - i\, \text{Flux}_X([\gamma]) = -i\, \text{Flux}_X([\gamma]).$$

From this identity it is evident that the vertical flux condition is equivalent to the exactness of $\phi_1$ and $\phi_2$. Since $\phi_1 - i\phi_2 = g^{-1}dh$ and $-(\phi_1, +i\phi_2) = gdh$, exactness of $\phi_1$ and $\phi_2$ is equivalent to the exactness of $g^{-1}dh$ and $gdh$. But from (3.6), it is clear that $\text{Period}_\gamma(X_\lambda) = Re \int_\gamma \Phi_\lambda = 0$, for all $\lambda > 0$ and all $[\gamma]$, if and only if $g^{-1}dh$ and $gdh$ are exact. $\square$

**Proof of Proposition 3.3.** Since the Lopez-Ros deformation, (3.5) and (3.6), changes neither zeros nor poles of $g$ and $dh$, the end types of $X_\lambda$ at the ends $p_j$ of $M = \overline{M} - \{p_1 \ldots p_r\}$ are the same for all $\lambda > 0$. That is, flat ends remain flat, and catenoid-type ends remain catenoidal; even more, their logarithmic growth rates, $\alpha = -\text{Residue}_p dh$, are independent of $\lambda$.

Fix $\lambda_0 > 0$ and consider two distinct ends $p_i, p_j$, $i \neq j$. Choose disjoint neighborhoods $D_i, D_j$ of $p_i$ and $p_j$ so that $X(D_i)$ and $X(D_j)$ are end representatives. Suppose $X_{\lambda_0}$ is an embedding. If these ends have different logarithmic growths, i.e., if $\alpha_i \neq \alpha_j$, the distance function from $X(D_i)$ to $X(D_j)$ is not only bounded away from zero but unbounded. This is because the ends are asymptotic to catenoid ends (or to a plane if one of the $\alpha$'s is zero) with different logarithmic growth rates. Since these growth rates are independent of $\lambda$, it follows that $X_\lambda(D_i) \cap X_\lambda(D_j) = \emptyset$ for $\lambda$ sufficiently close to $\lambda_0$. In case $\alpha_i = \alpha_j$ and $X_{\lambda_0}(D_i)$ and $X_{\lambda_0}(D_j)$ are asymptotic to ends with the same growth rate, we appeal to the Maximum Principle at Infinity proved in [52], which states that the distance between these embedded annular ends is bounded away from zero; i.e., they are *not* asymptotic at infinity. Thus, again, $X_\lambda(D_i) \cap X_\lambda(D_j) = \emptyset$ for $\lambda$ sufficiently close to $\lambda_0$. We may conclude that, for $\lambda$ sufficiently close to $\lambda_0$, $X_\lambda$ is an embedding. Stated differently, what we have proved is that

$$L := \{\lambda > 0 \mid X_\lambda \text{ is injective}\}$$

is an open set. It is nonempty since $1 \in L$.

We now wish to show that $L$ is closed. Suppose $\{\lambda_k\}$ is a sequence in $L$ that converges to $\lambda_0$, and suppose $X_{\lambda_0}$ is *not* an embedding; i.e., $X_{\lambda_0}(q_1) = X_{\lambda_0}(q_2)$ for some points $q_1, q_2 \in M, q_1 \neq q_2$. Because $X_{\lambda_k}$ converges, uniformly on compact subsets of $M$, to $X_{\lambda_0}$ it follows that some neighborhood $\mathcal{O}_1$ of $q_1$ is mapped by $X_{\lambda_0}$ to be on one side of the image of some neighborhood $\mathcal{O}_2$ of $q_2$. By the Maximum Principle, stated in Section 2.3.3, $X_{\lambda_0}(\mathcal{O}_1) = X_{\lambda_0}(\mathcal{O}_2)$ and by analyticity, $X_{\lambda_0}: M \to \mathbb{R}^3$ is a finite covering map, whose image $\mathcal{S} := X_{\lambda_0}(M)$ is a complete, embedded minimal surface of finite total curvature in $\mathbb{R}^3$. We will show that $X_{\lambda_0}$ is in fact one-to-one.

Given two disctinct punctures $p_i, p_j$, $i \neq j$, the quantity

$$Re \int_{p_i}^{p_j} dh_\lambda = Re \int_{p_i}^{p_j} dh$$

is a measure of the vertical distance between the two ends. Notice that it is independent of $\lambda$. Since $X$ is an embedding, the absolute value of this integral is infinite if the logarithmic growth rates at $p_i$ and $p_j$ are different. If the logarithmic growth rates are the same, it is still nonzero, a consequence of the Maximum Principle at Infinity quoted above. The mapping $X_{\lambda_0}$ takes a sufficiently small neighborhood of any puncture $p_j$ onto an end of $\mathcal{S}$, and each end of $\mathcal{S}$ is the image of some neighborhood of some puncture point. Hence the number of ends of $\mathcal{S}$ is not greater than the number of ends of $X_{\lambda_0}(M)$, i.e. not greater than $r$. If neighborhoods of two distinct punctures, $p_i, p_j$, $i \neq j$, are mapped to the same end of $\mathcal{S}$, then

$$Re \int_{p_i}^{p_j} dh = Re \int_{p_i}^{p_j} dh_{\lambda_0} = 0,$$

which is impossible. Thus the number of ends of $\mathcal{S}$ equals $r$ and, furthermore, we can find sufficiently small neighborhoods of $p_j$, $1 \leq i \leq r$, so that each neighborhood is mapped by $X_{\lambda_0}$ onto a different end of $\mathcal{S}$. This means that $X_{\lambda_0}$ is one-to-one near the punctures and hence one-to-one everywhere; i.e. $X_{\lambda_0}$ is injective.

□

**Proof of Proposition 3.4.** Suppose $q$ is a point where the Gauss map is vertical. Without loss of generality, we may assume that $g(q) = 0$. Choose conformal coordinates near $q$ so that $q$ corresponds to $z = 0$ and $g(z) = z^k$ for some $k > 0$ on $D_r(0) := \{|z| < r\}$. Because we are assuming regularity, $dh$ must have a zero of order $k$ at $0$, i.e.,

$$g(z) = z^k \qquad dh = z^k(a + zf(z))dz\,,$$

where $a$ is some nonzero complex constant and $f(z)$ is holomorphic on $D_r(0)$. For each $\lambda > 0$, the change of variables $z = \lambda^{-\frac{1}{k}}\zeta$ from $D_{\lambda^{\frac{1}{k}}r}(0)$ to $D_r(0)$ allows us to express $g_\lambda$ and $dh_\lambda$ as follows:

$$g_\lambda(z(\zeta)) = \zeta^k \qquad dh_\lambda = dh = \lambda^{-(1+\frac{1}{k})}\zeta^k(a + \lambda^{-\frac{1}{k}}\zeta f(\lambda^{-\frac{1}{k}}\zeta))d\zeta \tag{3.8}$$

on $D_{\lambda^{\frac{1}{k}}r}(0)$.

Rescale $dh_\lambda$ by a factor of $\lambda^{(1+\frac{1}{k})}$ so that $dh_\lambda = \zeta^k(a + \lambda^{-\frac{1}{k}}\zeta f(\lambda^{-\frac{1}{k}}\zeta))d\zeta$. The rescaled data produce immersions, which we will also call $X_\lambda$, that converge uniformly, on compact subsets of $\mathrm{C} = \lim_{\lambda \to \infty} D_{\lambda^{\frac{1}{k}}r}(0)$, to the immersion produced by the Weierstrass data

$$\widehat{g} = \zeta^k \qquad \widehat{dh} = a\zeta^k d\zeta\,. \tag{3.9}$$

We will call this immersion $\widehat{X} \colon \mathrm{C} \to \mathrm{R}^3$. The limit surface is a rescaled version of the $k$-Enneper surface defined and discussed in Remark 2.3ii). This surface is clearly not *embedded*. (See Figure 3.1.) Since rescaling does not create or destroy self-intersections, it follows that $X_\lambda$ is not injective for $\lambda$ sufficiently large. This completes the proof of the first statement of the proposition.

We now prove the second statement. Because a planar end has zero flux, $X_\lambda$ is well-defined for all $\lambda > 0$, according to Proposition 3.2. Because it has finite total curvature,

Figure 3.1 The limit surfaces in the proof of Proposition 3.4

Left: The k-Enneper surface given by (3.9), with $k = 2$.
Right: The limit surface given by (3.10), with $k = 2$. There is one horizontal flat end and one immersed end with winding number equal to three.

the Weierstrass data extends to zero meromorphically. Without loss of generality, we may assume that the limit normal at the end is $(0, 0, -1)$, so $g(0) = 0$. Because the end is flat, the zero of $g$ at $z = 0$ has order $k \geq 2$ and $dh$ has a zero of order $k - 2$ at $z = 0$. After a change of variables, if necessary, we may assume that

$$g(z) = z^k \qquad dh = z^k(az^{-2} + f(z))dz,$$

with $a \neq 0$ and $f(z)$ holomorphic on some disk $D_r(0) = \{|z| < r\}$. As in the proof of the first part of the proposition, we change variables by $z = \lambda^{-\frac{1}{k}}\zeta$ and rescale—this time by $\lambda^{1-\frac{1}{k}}$—so that

$$g_\lambda(\zeta) = \zeta^k \qquad dh = \zeta^k(a\zeta^{-2} + \lambda^{-\frac{2}{k}}f(\lambda^{-\frac{1}{k}}\zeta))d\zeta$$

on $D'_{\lambda^{\frac{1}{k}}r}(0)$. The associated immersions converge, uniformly on compact subsets of $\mathrm{C}-\{0\} = \lim_{\lambda \to \infty} D'_{\lambda^{\frac{1}{k}}r}(0)$, to an immersion $\widehat{X}\colon \mathrm{C} - \{0\} \to \mathrm{R}^3$ given by the Weierstrass data

$$\widehat{g}(\zeta) = \zeta^k \qquad \widehat{dh} = a\zeta^{k-2}d\zeta. \qquad (3.10)$$

This immersion is complete and has a flat end at 0. However, it has a nonflat end at infinity, and that end is not embedded. (See Figure 3.1.) It follows that $X_\lambda$ is not one-to-one for $\lambda > 0$ sufficiently large. □

# 4 Construction of the deformation family with three ends

We wish to construct complete, embedded minimal surfaces of genus $k-1, k \geq 2$ with three catenoid ends and $k$ vertical planes of symmetry that intersect in a common vertical line. They will be the surfaces $M_{k,x}$ of Theorem 3.3. These surfaces should be deformations of the examples $M_k$ of Costa ($k = 2$) and Hoffman-Meeks [35, 36] ($k \geq 3$), given in Theorem 3.2, which have two catenoid ends and a middle flat end. There will be, for each genus, a one parameter of such deformations. We begin by locating the zeros and poles of $g = \sigma \circ N$, the stereographic projection of the Gauss map. We will then show that the conformal types of the examples are quite restricted.

Assume we have such a surface $S$ for some fixed $k \geq 2$. Recall from Theorem 2.2 that the conformal type of $S$ is that of a compact Riemann surface punctured in three points (one for each end). Each of these ends is fixed by each of the reflections. Thus the closure of the fixed point sets of these reflective isometries contain each of the puncture points. We express this by saying that: "Each of the $k$ symmetry planes passes through each of the three punctures."

Because the surface is properly embedded, it divides $\mathbb{R}^3$ into two regions. The vertical line of intersection of the symmetry planes may be taken to be the $x_3$-axis. The symmetry planes make equal angles of $\pi/k$, otherwise the surface would have more symmetry planes by Schwarz reflection. The top and bottom ends must be catenoid-type. (This follows from Proposition 2.5. See also the discussion after Theorem 2.4.) Since these ends have the same limit-normals (because there are an odd number of ends and their normals must alternate between $(0, 0, +1)$ and $(0, 0, -1)$), it follows that, for $x_3 > 0$ sufficiently large, the $x_3$-axis is in one component of $\mathbb{R}^3 - S$, and in the other component for $x_3 < 0$, sufficiently negative. Thus the $x_3$-axis passes through the surface an odd number of times.

We will show that it passes through *exactly once*. The reflections in vertical planes generate a cyclic rotation group of symmetries, of order $k$, about the $x_3$-axis. Because $S$ is embedded, at any point $p$ of intersection of $S$ and the $x_3$-axis, $S$ must have a vertical normal vector. The order of the zero or pole of $g$ at $p$ is one less than the number of curves in $S \cap T_p M$. If $S$ has a rotational symmetry of order $k$ this means that the zero or pole has order at least $k - 1$. On the other hand, it follows from Proposition 2.2 or Proposition 2.4 that

$$\text{degree } g = \text{genus } S + 2 = k + 1\,. \tag{4.1}$$

If there is more than one axis point, then the top-most and the bottom-most points have the same orientation as the catenoid ends. Orient $S$ so that $g \circ N$ has zeros at these points. Then the number of zeros of $g$ is at least $2(k-1) + 2$, which is strictly greater than $k+1$ for $k \geq 2$. This contradicts (4.1). Hence there is a unique axial point where $g$ has a zero of order at least $k - 1$. But since $g$ has a zero at each of the extreme catenoid ends, this makes $k + 1$ zeros. We can conclude from (4.1) that the axial point is a zero of $g$ of order

Figure 4.0.1 Successive wedges of $S$.

exactly $k - 1$, and there are no other zeros.

We can now locate the poles of $g$. If the middle end were flat, then at the level of the end, say $x_3 = c$, the plane $\{x_3 = c\}$ must intersect $S$ in $2kj$ divergent curves, $j \geq 1$. Using this, we can easily show that $g$ has a pole of order $kj + 1$. By (4.1) we can conclude that $j = 1$ and that there are no other poles.

On the other hand, if the middle end is of catenoid-type, the Gauss map has a simple pole there, leaving $k$ more poles to be found at finite points of $S$. The planes of symmetry divide $S$ into $2k$ symmetric pieces. If a pole of the Gauss map is found in the interior of one of these regions, there must be at least $2k$ of them. This is a clear contradiction. The sole remaining possibility is to have the poles occurring in the symmetry planes, one on every other of the $2k$ halfplanes.

We will now make a sketch of the portion of $S$ in a wedge. Each symmetry line from the saddle must diverge to an end with the same normal as at the saddle, i.e. the top or the bottom end. Because the order of the Gauss map at the saddle is $k + 1$, adjacent curvature lines must diverge to *different* ends.

The picture suggests that this piece of surface is simply connected. Clearly this surface piece cannot have positive genus, else $S$ would have genus in excess of $2k$, a contradiction.

To show it is simply connected, we will compute its Euler characteristic. By (4.1), its total curvature is $-4\pi(k+1)/2k$. By taking a sufficiently large piece of the surface we see that each end contributes $\pi + \pi/k$ to the total curvature of the boundary. The first term comes from two right angles at vertices, while the second term comes from the fact that the surface is asymptotic to a catenoid (or plane) and we are seeing a $\pi/k$-arc of a circle, which has essentially no normal curvature. All the other boundaries are geodesics. Adding in a contribution of $\pi - \pi/k$ at the single axis point, we have from the Gauss-Bonnet formula

$$\begin{aligned} 2\pi\chi &= \int K dA + \int k_g + \sum_1^4 \text{ vertex angles} \\ &= \frac{-4\pi(k+1)}{2k} + 0 + 3(\pi + \frac{\pi}{k}) + (\pi - \frac{\pi}{k}) \\ &= 2\pi \,. \end{aligned}$$

Thus the Euler characteristic of the piece of $S$ in the wedge is 1 and it is simply-connected. We will now show that the sketch is qualitatively correct. There is but one component of the boundary and if we start at the top catenoid end, (considered now in the compactification as a vertex), neither edge emanating from this vertex can go to the bottom catenoid end. This is because the surface is embedded and the top and bottom ends have the same orientation. Hence one edge runs from the top end to the middle end, while the other runs to the axis point. After a Euclidean motion, we may assume that if the middle end is a catenoid end, then its logarithmic growth is positive, i.e. it goes up. Thus the two edges emanating from the top catenoid end must look qualitatively like Figure 4.0.1.

We may choose the wedge on the left in this figure.

## 4.1 Hidden conformal symmetries

We wish to determine the underlying conformal structure of the surfaces described above. In the previous section, we established that the compactification of the example of genus $k-1$, $k \geq 2$, assumed to have $k$ vertical planes of symmetry, can be decomposed into $2k$ geodesic 4-gons, all of whose vertex angles equal $\pi/k$. All the 4-gons have the same vertices. The Riemann surface must be chosen to have conformal involutions corresponding to the reflections. If we choose a metric of constant curvature, these conformal involutions are isometries. In the *constant curvature model* our surface is made up of congruent geodesic 4-gons with all angles equal to $\pi/k$ (rectangles in the case $k = 2$, hyperbolic 4-gons when $k > 2$) and with common vertices. For such 4-gons we need to prove a standard fact from hyperbolic geometry. This will allow us to determine the Riemann surface structure explicitly.

In the case $k = 2$ of genus one, a Euclidean rectangle has two reflective symmetry lines and, conformally, there is precisely a one-parameter of rectangles. For the case $k > 2$ we will need the same information from hyperbolic geometry about $\alpha$-*quadrangles*. By an $\alpha$-quadrangle we mean a quadrilateral in the hyperbolic plane, all of whose interior angles are equal to $\alpha$.

**Proposition 4.1** *For every $\alpha \in (0, \frac{\pi}{2})$ there exists a 1-parameter of $\alpha$-quadrangles. Every $\alpha$-quadrangle has all the symmetries of a Euclidean rectangle and, for each $\alpha$, there is a unique $\alpha$-square, i.e. an $\alpha$-quadrangle with all the symmetries of a square. (In particular its sides all have the same length.)*

Figure 4.1 Euclidean rectangles and hyperbolic $\alpha$-quadrangles.

Proposition 4.1 implies that the underlying Riemann surfaces of the minimal surfaces of interest to us have conformal symmetries that are not given by isometries of the minimal surfaces. This is important enough to prove.

**Proof of Proposition 4.1.** In the Euclidean plane we can parametrize the conformal family of rectangles by the angle $\gamma \in (0, \frac{\pi}{2})$ between the diagonal and an edge. For a given $\gamma$, take a triangle with angles $\gamma, \frac{\pi}{2}$ and $\frac{\pi}{2} - \gamma$, and rotate it about the midpoint of the edge opposite the right angle. The two triangles together give the rectangle associated to $\gamma$. In the hyperbolic plane: For any angles $\alpha, \beta, \gamma$, the condition $\alpha + \beta + \gamma < \pi$ is necessary and sufficient for the existence of a triangle with these three angles, and, such a triangle is unique up to hyperbolic isometries. Given $\alpha \in (0, \frac{\pi}{2})$, for any $\gamma \in (0, \alpha)$, let $\beta := \alpha - \gamma$ and consider the $\alpha, \beta, \gamma$ - triangle. Rotate this triangle about the midpoint of the edge opposite $\alpha$; the two triangles together form an $\alpha$-quadrangle. For each fixed $\alpha$, different values of $\gamma \in (0, \alpha)$ give conformally different $\alpha$-quadrangles because the $\alpha, \beta, \gamma$-triangles are different.

Next, we wish to show that any $\alpha$-quadrangle has the symmetries of a rectangle. Take any edge, $e_1$, of the $\alpha$-quadrangle, and let $\sigma$ be the geodesic that is the perpendicular bisector of $e_1$. Reflection in $\sigma$ maps $e_1$ into itself and interchanges the geodesics that extend the edges, $e_0, e_2$, adjacent to $e_1$. If the edge $e_3$ opposite $e_1$ is symmetric with respect to $\sigma$, we are done. If not, $e_3$ and its $\sigma$-reflection $\widetilde{e}_3$ intersect $\sigma$ and meet the geodesic extensions of $e_0$ and $e_2$ at an angle equal to $\alpha$. This produces two triangles with a common vertex at $e_3 \cap \widetilde{e}_3 \in \sigma$ and the same angle $\phi > 0$ there, whose other two angles are $\alpha$ and $\pi - \alpha$, a contradiction since the angle sum is $\pi + \phi > \pi$. Hence $\widetilde{e}_3 = e_3$ and the $\alpha$-quadrangle has the required symmetry.

We wish to show that every $\alpha$-quadrangle lies in the family constructed above. (Note that our proof that an $\alpha$-quadrangle has reflectional symmetries did *not* assume that the quadrangle was in this family.) Given an $\alpha$-quadrangle, the two symmetry geodesics cross at a point, which we will refer to as the *center*. Observe that 180° rotation about the center is a symmetry of the $\alpha$-quadrangle; it is the composition of the two reflective symmetries.

The geodesic join from any vertex to the center extends by this rotation to a diagonal of the $\alpha$-quadrangle with the center as its midpoint. The diagonal divides the $\alpha$-quadrangle into two $\alpha, \gamma, \alpha - \gamma$-triangles for some $\gamma \in (0, \alpha)$. Hence any $\alpha$-quadrilateral lies in the family we constructed.

If our $\alpha$-quadrangle was constructed as above with $\gamma = \alpha/2$, note that the diagonal produced from the edge of the $\alpha, \alpha/2, \alpha/2$ triangle (used to construct it) bisects the angles at its endpoints. It follows that reflection in this diagonal is a symmetry of this special $\alpha$-quadrangle. This is the $\alpha$-*square*. □

In the next section, we will build Riemann surfaces of genus $k - 1$ from $2k$ congruent $\pi/k$-quadrangles. Like rectangular tori, they form noncompact 1-parameter families.

### 4.2 The birdcage model

From the previous two sections we can conclude that the conformal structure of our minimal surface of genus $k-1$ is given by $8k$ quadrilaterals in hyperbolic (Euclidean) space when $k > 2$ ($k = 2$), with three right angles, one angle of $\pi/k$, and edges appropriately identified. In this section, which strictly speaking is not necessary for the next sections, we introduce an embedded surface in 3-space that realizes the conformal symmetries of our minimal surfaces in a geometrically evident manner. It has the same conformal type and will make the assembly of the surface out of quadrilaterals very clear.

Start with $k$ symmetrically placed half-meridians on the unit sphere, $k \geq 2$. Let $S = S_k$ be a tubular neighborhood of these curves in $\mathrm{R}^3$. The thickness, $\varepsilon > 0$, of the tubes will correspond to a conformal parameter. The $k$ symmetry planes of the meridians and the equatorial plane slice $S$ into $4k$ geodesic 4-gons. Each 4-gon has two vertices on the polar axis with angle $\pi/k$ and two vertices on the equator with angle $\pi/2$. For each quadrilateral we have from Gauss-Bonnet that

$$\int K dA = 2\pi - \left[2(\frac{\pi}{2}) + \frac{2\pi(k-1)}{k}\right],$$

or

$$\int K dA = \pi - \frac{2\pi(k-1)}{k} = -\pi + \frac{2\pi}{k}, \qquad (4.2)$$

which implies that the total curvature of $S$ is $-(4k - 8)\pi = 2\pi(2 - 2(k - 1))$. This verifies that the genus of $S$ is $k-1$, which is evident from its construction. Note that the reflectional symmetries in the meridian planes generate a $k$-fold rotational symmetry, and this symmetry fixes the four polar points.

The birdcage surface has the additional isometry of reflection in the equatorial plane. There is another conformal automorphism, which amounts to inversion through the sphere of radius $\sqrt{1 - \varepsilon^2}$. These two additional automorphisms restricted to the geodesic 4-gons produce the automorphisms corresponding to the symmetries of a Euclidean rectangle. In general, the composition of the equatorial reflection with the other conformal automorphism

Figure 4.2.0 The birdcage model, $k = 3$, and the fundamental 4-gon.

Figure 4.2.1 The birdcage decomposition The decomposition of the birdcage surface into hyperbolic polygons. Illustration for the case $k = 3$

      Upper left: The birdcage with one basic 4-gon darkened. This 4-gon has three right angles and one angle of $\pi/k$ at an axis point.
      Upper right: Four basic 4-gons darkened. Each one has an angle $\pi/k$ at a different axis point. Taken together, they form a quadrilateral all of whose angles equal $\pi/k$.
      Lower left: The birdcage decomposed into $2k$ $\pi/k$ quadrilaterals.
      Lower right: The birdcage decomposed into 4 right $2k$-gons, each centered at an axis point.

Figure 4.3 The conformal mapping used to define $u$

The illustration here is for k=3. The dotted vertex is the branch point of $u$. Its image is determined by the mapping, and is not assigned in advance.

is an orientation-*preserving* conformal involution, the *hyperelliptic involution*. The quotient of the birdcage by this involution is a sphere; in the next section the quotient map is constructed by means of the Riemann Mapping Theorem (using a choice of identification of $S^2$ with $C \cup \infty$), and called the function $u$. When the genus is greater than 1, the branch points of $u$ are the so-called *Weierstrass points* of the surface.

**Remark 4.1** *The most symmetric of these Riemann surfaces are also the conformal models for one series of Lawson's minimal surfaces in $S^3$ ([47]). They can be found by extending Plateau solutions, since the symmetry-diagonals of the quadrilaterals, the ones from the $\pi/k$-angles, are great-circle arcs of length $\pi/2$ in $S^3$.*

The birdcage model serves as a convenient visualization of the decomposition of the Riemann surface into $8k$ 4-gons. We will now work directly with the hyperbolic (or Euclidean when $k = 2$) 4-gons having three right angles and one angle of $\pi/k$. Note that each of the 4-gons with all angles equal to $\pi/k$, decomposes into four 4-gons, each having three right angles and one angle $\pi/k$ at one of the polar points.

In the next section we construct functions on the surface using this decomposition and the Riemann Mapping Theorem. This will allow us to make the transition from this hyperbolic picture to an algebraic description of our surfaces.

### 4.3 Meromorphic functions constructed by conformal mappings

Notice that each axis point is the center of a $2k$-gon, all of whose angles are right angles. These $2k$-gons are made up of $2k$ of the 4-gons with three right angles, the vertex with angle $\pi/k$ lying at the polar point. Consider one of these small 4-gons.

We map the quadrilateral into the sector $0 \leq \theta \leq \pi/k$ on the unit disk considered as the closed hemisphere of $S^2$ by taking the central vertex with angle $\pi/k$ into 0, its adjacent vertices to 1 and $e^{i\pi/k}$. This determines the image of the opposite vertex, which is a point on the unit circle. We will denote its value by $e^{2i\alpha/k}$, $0 < \alpha < \pi/2$. We will call this map $u$ and extend it to $S$ by Schwarz reflection. Note that at the midpoint of each of the four $2k$-gons, $u$ has a simple pole or a simple zero, and $u$ maps each of the four $2k$-gons onto a hemisphere. Thus, $u$ is a degree-2 meromorphic function and we may think of $S$ as a branched covering of $S^2$ with $2k$ simple branch points and branch values $e^{i(2\alpha+j\pi)/k}$.

Figure 4.3.1

In particular, $S$ is *hyperelliptic*. In the case $k = 2$, $u$ is clearly a degree-2 elliptic function. (In this case it follows directly from the construction of $S$ that $S$ is a rectangular torus.)

We observe that the symmetric construction of $u$ allows us to choose $\alpha$ such that $0 < \alpha \leq \pi/2$. Given a choice of $\alpha$ between 0 and $\pi/4$, simply choose as fundamental domain an adjacent quadrilateral and adjust $u$ by multiplication by $e^{\pi i/k}$.

### 4.3.1 The function $T$ and its relationship to $u$

Fix the same quadrilateral on which $u$ is defined as in the previous section. We map this domain onto a quarter-disk, specifying that the vertex opposite the $\pi/k$-vertex goes to the origin and the other two $\pi/2$-vertices are mapped to the unit circle. This determines the mapping up to a rotation in C, i.e. up to a unitary multiplicative factor. We rotate so that the $\pi/k$-vertex is mapped to $-1$. We call this mapping $T$ and we label (as in the diagram above) the other two vertices $-e^{-i\beta}$ and $-ie^{-i\beta}$, where the first one is in the upper half-plane.

We will show that, in fact, $\beta = \alpha$, where $\alpha$ is the angle defined by the mapping $u$. To do this we will construct, from $T$ and $u$, the same conformal map to the lower semi-disk. First of all $(T)^2$ maps our quadrilateral to a semi-disk, which we rotate by $-e^{2i\beta}$ so that $-e^{-i\beta} \to -1, -ie^{-i\beta} \to +1$, of course $0 \to 0$, and $-1 \to -e^{2i\beta}$. In particular, $-e^{2i\beta}(T)^2$ takes the $\pi/k$-vertex to $-e^{+2i\beta}$. Turning now to $u$, note that $u^k$ maps the quadrangle to the upper semi-disk. Composition of $u^k$ with the linear fractional transformation $\nu(\zeta) = \frac{\zeta - e^{2i\alpha}}{1 - \zeta e^{2i\alpha}}$ not only takes the upper semi-disk to the lower semi-disk, but also has the property that $\nu(u^k) = -e^{2i\beta}(T)^2$ on the three vertices opposite the $\pi/k$-vertex. Thus $\nu(u^k) = -e^{2i\beta}(T)^2$ everywhere, in particular at the $\pi/k$-vertex itself; i.e.

$$-e^{+2i\alpha} = \nu(0) = -e^{2i\beta}(-1)^2 = -e^{2i\beta}.$$

Since $\alpha$ and $\beta$ were both chosen in the range $(0, \frac{\pi}{2})$, $\alpha = \beta$ and

$$\frac{u^k - e^{2i\alpha}}{1 - u^k e^{2i\alpha}} = -e^{2i\alpha}(T)^2.$$

The above equation actually should be considered as the relation between the coordinate maps $T$ and $u$, and is in fact an equation for our surface.

## 4.4 The function $z$ and the equation for the Riemann surface in terms of $z$ and $u$

We prefer to modify $T$ by a Möbius transformation so that: its zeros and poles are more closely related to what we will need in the Weierstrass representation; and it is related to $u$ by an equation with real coefficients. Accordingly let $M(\zeta) = \frac{i(\zeta+1)}{1-\zeta}$ and define

$$z = M \circ T.$$

Note that $T = \frac{z-i}{z+i}$ and that $M$ takes $-1$ to $0$, $0$ to $i$ and $1$ to $\infty$. Since $T$ takes on the values $-1, 0, -ie^{-i\alpha}, -e^{-i\alpha}$ at the vertices of the fundamental 4-gon, we know that $z$ takes on the values $0$ and $i$ at two of the vertices. A calculation shows that $z$ takes on the value $\cot(\frac{\pi}{4} + \frac{\alpha}{2})$ at the vertex where $T$ takes on the value $-ie^{-i\alpha}$.[1]

We note that $z$ extends to the $2k$-gon with a zero of order $k$ at the midpoint, and covers the hemisphere $k/2$ times. Extending to the rest of $S$ by Schwarz reflection produces a function of degree $k$ whose zeros and poles (all of order $k$) occur at the special polar points. We chose the branch values of $u$ in order to relate $u$ to $z$ in a particularly simple way. We actually have this relationship already because we know that

$$\frac{u^k - e^{2i\alpha}}{1 - u^k e^{2i\alpha}} = -e^{2i\alpha} T^2 = -e^{2i\alpha}\left(\frac{z-i}{z+i}\right)^2.$$

Manipulation of this expression will yield (4.3). We prefer to derive (4.3) by matching the zeros and poles of $z$ and $u$. Consider $Z := -z + z^{-1}$. From the identity

$$\cot \alpha - \tan \alpha = 2 \cot 2\alpha$$

---

[1] Let $M(\zeta) = \dfrac{i(\zeta+1)}{1-\zeta}$. Since $M(e^{i\theta}) = i \cdot \frac{e^{i\theta}+1}{1-e^{i\theta}} = -\frac{e^{i\theta/2}+e^{-i\theta/2}}{(e^{i\theta/2}-e^{-i\theta/2})/i} = -\frac{\cos\theta/2}{\sin\theta/2} = -\cot\theta/2$. It follows that

$$M(-e^{-i\alpha}) = M(e^{i(\pi-\alpha)}) = -\cot\left(\frac{\pi}{2} - \frac{\alpha}{2}\right) = -\tan\frac{\alpha}{2},$$

$$M(e^{i(-\frac{\pi}{2}-\alpha)}) = -\cot\left(-\frac{\pi}{4} - \frac{\alpha}{2}\right) = \cot\left(\frac{\pi}{4} + \frac{\alpha}{2}\right)$$

Figure 4.4 The functions $z$ and $u$.

and the values of $z$ and $u$ given in the previous diagrams, we can easily determine the values of the functions $z^{-1} - z$, and $u^{-k}$ at vertices. Notice that the poles and zeros of $Z + \cot\alpha - \tan\alpha$ match those of $u^{-k}$ and $(\cot\alpha + \tan\alpha)u^{-k}$ equals $Z$ at the points where $u^{-k} = \pm 1$. Hence
$$z^{-1} - z + \cot\alpha - \tan\alpha = (\cot\alpha + \tan\alpha)u^{-k}.$$

This produces an equation for the surface $S$, which we write in two equivalent ways:
$$\begin{aligned} z^{-1} - z &= (x + x^{-1})u^{-k} - (x - x^{-1}) \\ z^{-1}(z-x)(z+x^{-1}) &= -(x+x^{-1})u^{-k}, \quad \text{where } x = \cot\alpha. \end{aligned} \tag{4.3}$$

Here, $x$ may be considered a conformal parameter for our of Riemann surfaces, $0 < x \leq 1$.

Figure 4.4.1

When $k = 2$, the Riemann surface is a rectangular torus. The fundamental domain constitutes $1/16th$ of the surface. Values of $u$ are indicated on the surface. The mapping is specified by fixing the values – $0, 1$, and $i$ at three vertices. It is branched at the fourth vertex; – its value is determined by the mapping, which is in turn determined by the proportions of the rectangle.

Differentiate this equation to get

$$-(z + z^{-1})\frac{dz}{z} = k(x + x^{-1})u^{-k}\frac{du}{u}$$

We now observe that branch points of $u$ occur at points where $z = \pm i$, as we know by construction. By construction or by calculation from equation (4.3) we can conclude that the branch values of $u$ are of the form $e^{\pm 2i\alpha/k}\omega$, where $\omega$ is a $k^{th}$ root of unity. At points where $z \neq 0, \infty$, poles of $u$ are branch points of $z$, with branch values $x = \cot \alpha$ or $-x^{-1} = -\tan \alpha$.

**Remark 4.2** *In Section 3, we wrote the above equations in a different and, for historical reasons, more familiar form. Multiply both sides of (4.3) by $z^k$ and define $w := z/u$. Then we may write (4.3) as*

$$\left(\frac{z}{u}\right)^k = w^k = -X^{-1}z^{k-1}(z - x)(z + x^{-1}) \qquad X = (x + x^{-1}), \tag{4.4}$$

*and $x, 0 < x \leq 1$, parametrizes our family of Riemann surfaces. In case $k = 2$, we have*

$$w^2 = -X^{-1}z(z - x)(z + x^{-1}), \qquad x \geq 1 \tag{4.5}$$

*and we recognize the equation of a rectangular torus. The function $z$ is thus seen to be a geometrically normalized form of the Weierstrass $\wp$-function: $z = c_1(\wp - c_2)$, where $c_2$ is such that $z$ has a double zero at the center of one of the 4 gons.*

We also note that on the Riemann surfaces (4.4) one can define data for Chen-Gackstatter-type surfaces (see Section 2.2) with higher symmetry using

$$g = \rho \cdot \left(\frac{w}{z}\right)^k = \rho \cdot u^{-k},$$

$$dh = i(z + z^{-1} - 2)\frac{dz}{z}.$$

*(See Karcher [43].) This construction has been extended by Thayer [71].*

## 4.5 The Weierstrass data

We now determine the Weierstrass data for our minimal surface of genus $k-1, k \geq 2$, with three vertical catenoid ends and $k$ symmetrically placed vertical planes of reflective symmetry. We know that the underlying Riemann surface will be $S$, as presented in equation (4.3) or (4.4). Note that we have a 1-parameter family depending on $\pi/4 \leq \alpha < \pi/2$ (or equivalently $x^{-1} = \tan \alpha \geq 1$).

In terms of our decomposition of $S$, we will look at the region between two adjacent symmetry halfplanes. This corresponds to four basic quadrilaterals joined together to form a $\pi/k$-quadrangle with vertices consisting of the polar points of the birdcage where the vertex angle is $\pi/k$. The top and bottom ends are located at points where $u = \infty$. The middle end is placed where $z = \infty$, which is where $u = 0$.

We will specify the Gauss map here by noticing that the values of $g$ must be 0 or $\infty$ at the ends. (See Figure 4.5.1 and the table in Figure 4.5.0.) Since there is no branching of the Gauss map at a catenoid end, these are simple zeros and poles. By contrast, at the saddle we must have at least $k$ curves of intersection with the horizontal tangent plane, which implies that the branching of $g$ at the saddle has order at least $k-1$. Orient the surface so that $g = 0$ at the saddle. We know that the degree of $g$ must be equal to $k+1$, and since $g = 0$ at two other catenoid ends, we can conclude that $g$ has a zero of order $k-1$ at the saddle. Referring to Figure 4.5.0, we note that

$$g_1 := \frac{u^{-1}z}{mz+1} \tag{4.6}$$

| $(z,w) =$ | $(0,0)$ | $(x,0)$ | $(-x^{-1},0)$ | $\infty,\infty$ | $z = -m^{-1}$ |
|---|---|---|---|---|---|
|  |  | Bottom | Top | Middle | Finite |
|  | Saddle | End | End | End | Vertical Points |
| $g$ | $0^{k-1}$ | 0 | 0 | $\infty$ | $\infty$ |
| $u^{-1}$ | $\infty$ | 0 | 0 | $\infty$ | $-$ |
| $z$ | $0^k$ | $x$ | $-x^{-1}$ | $\infty^k$ | $-m^{-1}$ |
| $mz+1$ | $-$ | $-$ | $-$ | $\infty^k$ | 0 |
| $dh$ | $0^{k-1}$ | $\infty$ | $\infty$ | $\infty$ | 0 |
| $u^k$ | $0^k$ | $\infty^k$ | $\infty^k$ | $0^k$ | $-$ |
| $dz/z$ | $\infty$ | $0^{k-1}$ | $0^{k-1}$ | $\infty$ | $-$ |
| $du/u$ | $\infty$ | $\infty$ | $\infty$ | $\infty$ | $-$ |
| $\frac{mz+1}{z+z^{-1}}$ | $0^k$ | $-$ | $-$ | $-$ | 0 |

Figure 4.5.0

has the correct zero/pole structure at the saddle and the ends, while having a simple pole at an edge-point in between. Here, $m$ is a free parameter, placing the finite vertical normal at the points where $z = -m^{-1}$. Placement of the vertical normal as it is in the picture requires us to choose $-\infty < (-m^{-1}) < -1/x$, or

$$0 < m < x = \cot \alpha \,. \tag{4.7}$$

Otherwise, the finite vertical normal would be positioned elsewhere. Fixing $m$ determines the Gauss map up to a real multiplicative factor. (Multiplying by a unitary complex number produces a rotation of the surface.)

$$g = \rho g_1 \qquad \rho > 0 \,. \tag{4.8}$$

We are now in a position to determine the complex differential $dh$ of the height function. It must have a simple pole at each catenoid end, a simple zero at the (simple) finite pole of $g$ and a zero of order $k-1$ at the saddle. Direct calculation of the zeros and poles of $\dfrac{dz}{z}, u^k$ and $mz+1$ is straightforward. (See Figure 4.5.1.) We may conclude that

$$dh = cu^k(mz+1)\frac{dz}{z}, \tag{4.9}$$

where $c$ is a nonzero complex multiplicative constant. We wish to show that this constant must be chosen to be real. First, of all we note that the curves on the boundary of this fundamental quadrilateral correspond to planar geodesics (in vertical planes) along which the Gauss map must take values in a meridian circle. Thus $\dfrac{dg}{g}$ is real on the tangent vectors of these curves. On the other hand, the quadratic differential associated with the second fundamental form of a minimal surface is given as in (2.12) by

$$\langle S(V)), V \rangle = Re\left(\frac{dg(V)}{g} \cdot dh(V)\right). \tag{4.10}$$

Since these planar geodesics are also principal curves, we must have $\dfrac{dg}{g}(V) \cdot dh(V)$ real along these edges. (See 2.14.) But this implies that $dh$ is real along these edges. From our definitions of $u$ and $z$, it follows that $z$ and $du/u$ are real along these curves. Hence from (4.10), it follows that $c$ is real. We rescale so that this constant is equal to one. Using (4.3) we may write:

$$dh = u^k(mz+1)\frac{dz}{z} = \frac{-(x+x^{-1})(mz+1)dz}{(z-x)(z+x^{-1})}\,. \tag{4.11}$$

We will use both expressions for $dh$.

Summary

We summarize what we have achieved.

Figure 4.5.1 The portion of the surface in a wedge and its corresponding $\pi/k$-gon.

The Riemann surface and Weierstrass data are as in (4.4), (4.6)–(4.8), (4.11):

$$(x + x^{-1})u^{-k} = z^{-1} - z + (x - x^{-1}) = -z^{-1}(z - x)(z + x^{-1}),$$

$x = \cot\alpha, \pi/4 \leq \alpha < \pi/2$, minus the points $(z, u) = (x, \infty), (-x^{-1}, \infty)(\infty, 0)$.

The end at $(\infty, 0)$ is the middle end.

$$g = \rho g_1 = \rho \frac{u^{-1}z}{mz + 1} \quad \rho > 0, \quad 0 \leq m < x;$$

$$dh = u^k(mz + 1)\frac{dz}{z} = \frac{-(x + x^{-1})(mz + 1)dz}{(z - x)(z + x^{-1})}$$

(4.12)

This data will produce in the Weierstrass Representation

$$X(p) = (x_1, x_2, x_3)(p) = Re \int_{p_0}^{p} \Phi,$$

$$\Phi = \left(\frac{1}{2}(g^{-1} - g), \frac{i}{2}(g^{-1} + g), 1\right) dh,$$

a (possibly multivalued) minimal immersion into $R^3$. From the development in Sections 4.1 to 4.5, it is clear that if the surfaces $M_{k,x}$ of Theorem 3.3 (also described at the beginning of Section 4) exist, they must have a Weierstrass Representation in the form above, for appropriate choices of constants $k, x, m$ and $\rho$. All Weierstrass data of the above form produces the symmetry of the surfaces that we require. To see this, first note that the conformal automorphisms

$$f(z, u) := (\bar{z}, \bar{u})$$
$$\tau(z, u) := (z, \varepsilon u) \qquad \varepsilon = e^{-\frac{2\pi i}{k}}$$

fix $p_0 := (0, 0)$ and the ends $(x, \infty), (-x^{-1}, \infty)$ and $(\infty, 0)$. We choose to integrate from this $p_0$ in the expression for $X(p)$, and observe that we may write

$$x_1 - ix_2(p) = \int_{p_0}^{p} g^{-1}dh - \overline{\int_{p_0}^{p} gdh}$$

$$x_3(p) = Re \int_{p_0}^{p} dh.$$

Since $f^*(dh) = \overline{dh}$ and $\tau^*(dh) = dh$, it is evident that $x_3(\tau(p)) = x_3(f(p)) = x_3(p)$. Because $f^*g = \bar{g}$ and $\tau^*g = \varepsilon g$, it follows that

$$\begin{aligned} x_1 - ix_2(f(p)) &= \overline{x_1 - ix_2(p)} &= x_1 + ix_2(p) \\ x_1 - ix_2(\tau^j(p)) &= \bar{\varepsilon}^j(x_1 - ix_2(p)) &= 1 \leq j < k. \end{aligned}$$

Thus $f$ induces reflection in the $(x_1, x_3)$-plane, while $\tau^j$ induces rotation by $2\pi j/k$ about the vertical axis. (We could have used the Schwarz Reflection Principle to prove that $f$ induces a reflection.)

Figure 4.6 Logarithmic growth rates at the ends

## 4.6 The logarithmic growth rates

We proceed to determine the logarithmic growth rates of the catenoid ends. From Proposition 2.5*vi)*, we know that this growth rate is equal to minus the residue of $dh$ at the end. At the bottom end, where $z = x$, and at the top end, where $z = -x^{-1}$, $dz$ is regular, while at the middle end, where $z = \infty$, $dz/z$ has a simple pole with residue $-1$. From (4.12) we can read off the residue of $dh$ at $z = x$ and $z = -x^{-1}$:

$$\text{Residue}(dh)_{z=x} = \frac{-(x+x^{-1})(mx+1)}{x+x^{-1}} = -(mx+1)$$
$$\text{Residue}(dh)_{z=-x^{-1}} = \frac{-(x+x^{-1})(-mx^{-1}+1)}{(x^{-1}-x)} = -\frac{m}{x}+1\,. \tag{4.13}$$

Because the residues must sum to zero (and $dh$ has poles only at the three ends) we get for the middle end

$$\text{Residue}(dh)_{z=\infty} = m(x+x^{-1})\,.$$

We could also compute this directly by rewriting the expression for $dh$ in (4.12) in terms of $z$ and expressing it in terms of $z^{-1}$.

$$dh = \frac{(x-x^{-1})(m+z^{-1})}{z^{-2}+(x-x^{-1})z^{-1}-1}\,\frac{dz^{-1}}{z^{-1}}\,.$$

The logarithmic growth rates are the negatives of these residues: With $m \in [0, x)$ arbitrary, we have $-(mx+1) < 0 < m(x+x^{-1})$. However, the "middle" end can grow faster than the "top" end, in which case the surface would not be embedded. To avoid this we have the growth rate condition
$$m(x + x^{-1}) < -mx^{-1} + 1.  \tag{4.14}$$

## 4.7 The period and embeddedness problems for the surfaces $M_{k,x}$

We must show that the surfaces $M_{k,x}$ in Theorem 3.3 are in fact well defined by the mapping $X: \Sigma \to \mathrm{R}^3$ in the Weierstrass representation, (2.7) and (2.8) on the Riemann surface $\Sigma$ given by (4.12). That is, according to (2.9), we must show that

$$\operatorname{Period}_\alpha(\Phi) = Re \oint_\gamma \Phi = Re \oint_\gamma (\phi_1, \phi_2, \phi_3) = \vec{0} \tag{4.15}$$

for all closed curves $\gamma$ on the punctured Riemann surface. These formulae are restated in the summary at the end of Section 4.5, where we proved that the immersions must have the required Euclidean symmetry. The built-in symmetry of the mapping is crucial. Consider a fundamental quadrilateral as in Figure 4.7. We label the saddle by $S$ and the top, middle and bottom ends by $T, M$ and $B$, respectively. Opposite sides of the quadrilateral are mapped into parallel vertical planes by $X = Re \int \Phi$. We require that these planes coincide. After a rotation, we can assume that the curve $X(\widehat{ST})$ lies in the $x_2 = 0$ plane. For $X(\widehat{BM})$ to lie in the same plane, we need only check that

$$Re \oint_{\gamma_1} \phi_2 = 0, \tag{4.16}$$

where $\gamma$ is the symmetry line of the quadrilateral between edges $ST$ and $BM$.

The curves $X(\widehat{SB})$ and $X(\widehat{TM})$ lie in planes parallel to the plane $x_2 \cos \pi/k = x_1 \sin \pi/k$. They lie in the same plane if and only if

$$Re \oint_{\gamma_2} (\phi_2 \cos \pi/k - \phi_1 \sin \pi/k) = 0. \tag{4.17}$$

Notice that the symmetry forces $Re \int_\alpha \phi_3 = 0$ for every closed curve $\alpha$ on the punctured Riemann surface; we get no conditions from the third differential.

There is no geometric reason why the second period condition (4.17) should be more complicated than the first. In fact it is sufficient to look in detail at the first condition (4.16) and, at the end, handle the second condition by a parameter substitution. To see this, consider the following expressions for $\phi_1$ and $\phi_2$.

$$\begin{aligned} 2\rho\phi_1 &= (g_1^{-1} - \rho^2 g_1)dh = (u(\frac{1+mz}{z})^2 - \rho^2 u^{-1})u^k dz \\ 2\rho\phi_2 &= i(g_1^{-1} + \rho^2 g_1)dh = i(u(\frac{1+mz}{z})^2 + \rho^2 u^{-1})u^k dz\,. \end{aligned} \tag{4.18}$$

If we substitute, in the expression for $\phi_2$, $\widetilde{u} = e^{i\pi/k}u, \widetilde{z} = -z$ and $\widetilde{m} = -m$ and define $\widetilde{\phi}_2 = \phi_2(\widetilde{u}, \widetilde{z}, \widetilde{m}, \rho)$, we have $\widetilde{\phi}_2 = \phi_2 \cos \pi/k - \phi_1 \sin \pi/k$. Thus the second period condition (4.17) can be written in the form

$$Re \oint_{\gamma_2} \widetilde{\phi}_2 = 0. \tag{4.17'}$$

Figure 4.7 The period problem

Now note that $\tilde{u}$ and $\tilde{z}$ satisfy the same Riemann surface equation as $u$ and $z$, provided we substitute $\tilde{x} = 1/x$ for $x$:

$$(\tilde{x} + \tilde{x}^{-1})\tilde{u}^{-k} = (x + \frac{1}{x})(-1)u^{-k} = -(z^{-1} - z) - (x - x^{-1})$$
$$= (\tilde{z}^{-1} - \tilde{z}) + (\tilde{x} - \tilde{x}^{-1}).$$

Notice that $\tilde{x} = x^{-1} = \tan\alpha = \cot(\tilde{\alpha})$, where $\tilde{\alpha} = \pi/2 - \alpha$. Since we will explicitly express (4.16) in terms of $m$ and $\alpha$, we will be able to express (4.17) by the substitution $m \to \tilde{m} = -m$ and $\alpha \to \tilde{\alpha} = \pi/2 - \alpha$. Using (4.18) we may express (4.16) and (4.17) in terms of the integrals

$$Q_j = Q_j(m,\alpha) := \int_{\gamma_j} (\frac{1+mz}{z})^2 u^{k+1} dz,$$

$$C_j = C_j(\alpha) := \int_{\gamma_j} u^{k-1} dz,$$

$j = 1, 2$; namely,

$$Q_1(m, \alpha) = \rho^2 C_1(\alpha),$$
$$Q_2(\tilde{m}, \tilde{\alpha}) = \rho^2 C_2(\tilde{\alpha}) \qquad \tilde{m} = -m \;\; \tilde{\alpha} = \tfrac{\pi}{2} - \alpha. \tag{4.19}$$

For each $\alpha$, we want to find $(\rho, m)$ satisfying these conditions, which express the period problem. Note that $Q_j$ is quadratic in $m$, while, $C_j$ depends on $\alpha$ only. Each condition in (4.19) determines $\rho^2 = \rho^2(m, \alpha)$, and it is easy to check that $\rho^2 > 0$. The two conditions agree if and only if

$$\frac{Q_1(m,\alpha)}{Q_2(\tilde{m},\tilde{\alpha})} = \frac{C_1(\alpha)}{C_2(\tilde{\alpha})} \quad \text{(Compatibility condition for } m = m(\alpha)\text{)}. \tag{4.20}$$

Our strategy is to reduce this condition to a linear equation in $(m - m^{-1})$, whose obvious solvability shows that, for each $\alpha$, we can find $m(\alpha), \rho(\alpha)$ solving the period problem (4.19). We then need to show that the logarithmic growth rates are correctly ordered. That is, the growth rate at the middle end (the one corresponding to the puncture $(z, u) = (\infty, 0)$) lies between the growth rate at the other two ends. From (4.13) and (4.14), this is equivalent to

$$2\tan\alpha + \cot\alpha = 2x^{-1} + x < m^{-1} \quad \text{(Growth rate condition for } m = m(\alpha)\text{)}. \tag{4.21}$$

**Proposition 4.2** *Fix $k \geq 2$. For every $\pi/4 \leq \alpha < \pi/2$, there exists a unique non-negative $m(\alpha) < \cot\alpha$ that satisfies the compatibility condition (4.20). Moreover $m(\alpha)$ satisfies the growth condition (4.21) and is a continuous function of $\alpha$.*

We will give a complete proof of this Proposition in the case $k > 2$. For the torus case, $k = 2$, the arguments are a bit more difficult. (See [30].) As indicated above, existence of the surfaces $M_{k,x}$ in Theorem 3.3 follows from Proposition 4.2. When $\alpha = \frac{\pi}{4}, m(\alpha) = m(\frac{\pi}{4}) = 0$, and the middle end is flat; these are the surfaces $M_{k,1}$ of Theorem 3.3, which correspond to the surfaces $M_k$ of Theorem 3.2. These surfaces $M_k$ are embedded. Thus each $M_{k,x}$ lies in a continuous family of examples whose first member is embedded. The embeddedness of the surfaces in the $M_{k,x}$ is a consequence of the growth condition (4.21) and the following result from [30].

**Proposition 4.3** *Suppose $h_t: N \to R^3$, $0 \leq t \leq A$ is a continuous one-parameter family of complete minimal surfaces of finite total curvature. Suppose that for all $t$, the ends of $h_t(N) = N_t$ are all vertical, and that $h_0$ is an embedding but $h_A$ is not. Let $T = \sup\{t | h_t$ is an embedding.$\}$ Then $h_T$ is an embedding if $T < A$ and if one orders the ends by height, there are at least two ends with the same logarithmic growth.*

## 4.8 The details of the solution of the period problem, I. Simplification of the integrals.

Recall from (4.3) (or (4.12)) that the Riemann surface equation is

$$z^{-1} - z = (x + x^{-1})u^{-k} - (x - x^{-1}). \tag{4.22}$$

Differentiation gives us a relation between the differentials $dz$ and $du$:

$$(z^{-1} + z)\frac{dz}{z} = k(x + x^{-1})u^{-k}\frac{du}{u}. \tag{4.23}$$

Two differential forms have the same periods, or are homologous, $w \sim \eta$, if they differ by the differential of a function: $w = \eta + df$. For example:

$$g_1 dh = u^{k-1} dz \sim -(k-1)zu^{k-2} du.$$

This gives us explicit integral expressions for $C_1$ and $C_2$ defined in the previous section.

**Lemma 4.1**

$$C_1 = \operatorname{Re} \int_{\gamma_1} ig_1 dh = \frac{2}{k}(k-1)(x + x^{-1}) \int_0^\alpha \sqrt{\cos^2 \phi - \cos^2 \alpha} \, \cos((1 - \frac{2}{k})\phi) d\phi$$

$$C_2 = \frac{2}{k}(k-1)(x + x^{-1}) \int_0^{\widetilde{\alpha}} \sqrt{\cos^2 \phi - \cos^2 \widetilde{\alpha}} \, \cos((1 - \frac{2}{k})\phi) d\phi.$$

**Proof.** The path $\gamma_1$ is symmetric with respect to a branch point of $u$; therefore $u$ has the same values at opposite points and $du$ changes sign at the branch point. Moreover $z$ has the value $i$ at this point and is symmetric under $180°$ rotation, i.e. its values at opposite points are $z, -1/z$. Hence

$$\int_{\gamma_1} 2g_1 dh = -(k-1) \int_{\gamma_1} (z + z^{-1})u^{k-1}\frac{du}{u}.$$

We can eliminate $z + z^{-1}$ as follows:

$$\left(\frac{z + 1/z}{x + 1/x}\right)^2 = \left(\frac{1/z - z}{x + 1/x}\right)^2 + \frac{4}{(x + 1/x)^2}$$

$$= \left(u^{-k} - \frac{x - 1/x}{x + 1/x}\right)^2 + \frac{4}{(x + x^{-1})^2}.$$

Recall $x = \cot \alpha$. Hence $(x - x^{-1})(x + x^{-1})^{-1} = \cos 2\alpha$; with this

$$\left(\frac{z + 1/z}{x + 1/x}\right)^2 = u^{-k} \cdot (u^k + u^{-k} - 2\cos 2\alpha). \tag{4.24}$$

Furthermore, along the path $\gamma_1$ we can take $u$ as coordinate function on the Riemann surface; since $u$ is unitary along $\gamma_1$ we have

$$u = e^{i\phi}, u^k + u^{-k} = 2\cos\phi, \quad 0 \le \phi \le \frac{2\alpha}{k}$$

$$du/u = id\phi.$$

We insert this and (4.24) in the last integral and get

$$\operatorname{Re} \int_{\gamma_1} ig_1 dh = (k - 1)(x + x^{-1}) \int_0^{2\alpha/k} \sqrt{2\cos k\phi - 2\cos 2\alpha} \, \cos\left(\left(\frac{k}{2} - 1\right)\phi\right) d\phi.$$

We substitute $t = k\phi$ and use $\cos 2\alpha = 2\cos^2 \alpha - 1$ to derive the expression for $C_1$ in the lemma. Using the symmetric treatment discussed above, we get $C_2$ if we replace $\alpha$ by $\widetilde{\alpha} := \frac{\pi}{2} - \alpha$. □

Now we turn to $Q_1(m)$. We have at first

$$\frac{1}{g_1}dh = \left(m^2 + \frac{1}{z^2} + \frac{2m}{z}\right)u^{k+1}dz,$$

so that $\int_{\gamma_1} \frac{1}{g_1}dh$ is given in terms of three different period integrals on the Riemann surface. We have to reduce these integrals to the two in the previous lemma and to two others before we can discuss them. First, multiply the differential relation (4.23) by $u^{k+\ell}$ to get

$$u^{k+\ell}(1 + z^{-2})dz = k(x + x^{-1})u^\ell \cdot \frac{du}{u} \underset{[\ell \ne 0]}{\sim} 0.$$

In particular,

$$\frac{1}{z^2}u^{k+1}dz \sim -u^{k+1}dz.$$

Hence we achieved a first reduction to two integrals per path:

$$\frac{1}{g_1}dh \sim \left(m^2 - 1 + \frac{2m}{z}\right)u^{k+1}dz.$$

We deal with $u^{k+1}dz \sim -(k+1)u^{k+1}zdu/u$ as before by using the symmetry along $\gamma_1$:

$$\int_{\gamma_1} 2u^{k+1}dz = -\int_{\gamma_1}(k + 1)u^{k+1}\left(z + \frac{1}{z}\right)\frac{du}{u}.$$

Again $(z + 1/z)$ is eliminated using (4.24) to get almost the same integral as in the lemma:

$$2Re \int_{\gamma_1} iu^{k+1}dz = \frac{4(k+1)}{k}(x + x^{-1}) \int_0^\alpha \sqrt{\cos^2 \phi - \cos^2 \alpha} \cdot \cos((1 + \frac{2}{k})\phi)d\phi. \quad (4.25)$$

The remaining integral $\int u^{k+1}dz/z$ can be rewritten as a linear combination of the integrals we already have. We multiply the differential relation (4.23) to get

$$(z + \frac{1}{z})u^{k+1}dz = k(x + x^{-1})zdu$$

and we multiply the surface equation (4.22) to get

$$(\frac{1}{z} - z)u^{k+1}dz = (x + x^{-1})udz + (x^{-1} - x)u^{k+1}dz.$$

Before we add these two equations, we use $udz \sim -zdu$. Then

$$\frac{2}{z}u^{k+1}dz \sim (k-1)(x + x^{-1})zdu + (x^{-1} - x)u^{k+1}dz.$$

The second differential on the right has been handled in (4.25). The remaining one is treated with path-symmetry as before:

$$\int_{\gamma_1} 2zdu = \int_{\gamma_1} (z + \frac{1}{z})du,$$

and $(z + 1/z)$ is eliminated with (4.24). We obtain

$$2Re \int_\gamma izdu = -\frac{4}{k}(x + x^{-1}) \int_0^\alpha \sqrt{\cos^2 \phi - \cos^2 \alpha} \cos((-1 + 2/k)\phi)d\phi. \quad (4.26)$$

This reduction to the same integral as in the lemma allows us to finish the period discussion. (That this should happen seemed at first an undeserved piece of good luck. On reflection, it is simple to see why this happens, at least in the case of genus one: If one has two differential forms without residues and with linearly independent periods– as is the case in the situation at hand– then the periods of *any* differential form without residues are linear combinations of the periods of the first two.) First shorten the notation by defining

$$E_\mp(\gamma) := \int_0^\gamma \sqrt{\cos^2 \phi - \cos^2 \gamma} \cdot \cos((1 \mp 2/k)\phi)d\phi,$$
$$c := \frac{2}{k}(x + x^{-1}). \quad (4.27)$$

In these terms the lemma says

$$C_1 = c \cdot (k-1) \cdot E_-(\alpha), \quad C_2 = c \cdot (k-1) \cdot E_-(\tilde{\alpha}).$$

The simplifications concentrated in (4.25), (4.26) deal with all the terms in

$$\frac{1}{g_1}dh \sim (m^2 - 1 + \frac{2m}{z})u^{k+1}dz$$

to give

$$Q_1(m) := Re \int_{\gamma_1} -i\frac{dh}{g_1}$$
$$= c \cdot [(1-m^2)(k+1)E_+(\alpha) + (\cot\alpha - \tan\alpha) \cdot m \cdot (k+1)E_+(\alpha)$$
$$+ (\cot\alpha + \tan\alpha) \cdot m \cdot (k-1)E_-(\alpha)].$$

The substitution $\alpha \to \widetilde{\alpha} := \pi/2 - \alpha$, $m \to \widetilde{m} := -m$ gives

$$Q_2(m) = c \cdot [(1-m^2)(k+1)E_+(\widetilde{\alpha}) + (\cot\alpha - \tan\alpha) \cdot m \cdot (k+1)E_+(\widetilde{\alpha})$$
$$- (\cot\alpha + \tan\alpha) \cdot m \cdot (k-1)E_-(\widetilde{\alpha})].$$

The compatibility condition (4.20) is now fractional linear in $M := \frac{1}{m} - m$ and of the form

$$\frac{M + a + b \cdot q}{M + a - b \cdot \widetilde{q}} = \frac{q}{\widetilde{q}},$$

namely:

$$\frac{(\frac{1}{m} - m) + (\cot\alpha - \tan\alpha) + (\cot\alpha + \tan\alpha)\frac{k-1}{k+1} \cdot \frac{E_-(\alpha)}{E_+(\alpha)}}{(\frac{1}{m} - m) + (\cot\alpha - \tan\alpha) - (\cot\alpha + \tan\alpha)\frac{k-1}{k+1} \cdot \frac{E_-(\widetilde{\alpha})}{E_+(\widetilde{\alpha})}} = \frac{E_-(\alpha)/E_+(\alpha)}{E_-(\widetilde{\alpha})/E_+(\widetilde{\alpha})}$$

This simplifies to

$$((\frac{1}{m} - m) + (\cot\alpha - \tan\alpha)) \cdot \frac{k+1}{k-1} \cdot \left(\frac{E_+(\widetilde{\alpha})}{E_-(\widetilde{\alpha})} - \frac{E_+(\alpha)}{E_-(\alpha)}\right) = 2 \cdot (\cot\alpha + \tan\alpha). \quad (4.28)$$

## 4.9 The details of the solution of the period problem, II. The monotonicity lemma.

We will prove the following

**Lemma 4.2 (The Monotonicity Lemma)** *Fix $k \geq 2$. The function $\gamma \to \frac{E_+(\gamma)}{E_-(\gamma)}$ is strictly decreasing on $(0, \frac{\pi}{2})$, with*

$$\lim_{\gamma \to 0} \frac{E_+(\gamma)}{E_-(\gamma)} = 1, \quad \lim_{\gamma \to \frac{\pi}{2}} \frac{E_+(\gamma)}{E_-(\gamma)} = \frac{k-1}{k+1}.$$

*The integrals $E_\pm(\gamma)$ are defined in (4.27).*

Before presenting the proof, we will show how the Lemma is used to prove Proposition 4.2. Recall from Section 4.3 that the Riemann surfaces determined by the conformal parameter $\alpha$, with $\alpha > \pi/4$, are obtained from those with $\alpha < \pi/4$ by exchanging generators. The choice of $\alpha > \pi/4$ was made to be compatible with the location of the point where the normal is vertical (i.e. where $z = -m^{-1}$) on the planar geodesic symmetry curve *between* the end-punctures at the points where $z = -x^{-1} = -\tan\alpha$ and $z = \infty$. To see this, let $\alpha > \pi/4$. Then the Monotonicity Lemma implies

$$\frac{2}{k+1} = \frac{E_+}{E_-}(0) - \frac{E_+}{E_-}(\frac{\pi}{2}) \geq \frac{E_+(\tilde{\alpha})}{E_-(\tilde{\alpha})} - \frac{E_+(\alpha)}{E_-(\alpha)} > 0.$$

With this, it is clear that (4.28) is *always* solvable for $m^{-1} - m$ and that (for $\alpha > \pi/4$)

$$m^{-1} - m > (k-1)(\cot\alpha + \tan\alpha) + (\tan\alpha - \cot\alpha) > 0 \tag{4.29}$$

$$m^{-1} - \tan\alpha + (\cot\alpha - m) > (k-1)(\cot\alpha + \tan\alpha) > 0.$$

This shows that $M := m^{-1} - m > 0$, that $m = -\frac{M}{2} + \sqrt{1 + \frac{M^2}{4}} < 1/m$, and that $m < \cot\alpha$. Therefore the positive solution $m = m(\alpha)$ of (4.28) (which is our compatibility condition (4.20)), puts the vertical normal on the symmetry line defined by $-\infty < z < -\tan\alpha = -x^{-1}$, that runs between the punctures at $z = -\tan\alpha$ and $z = \infty$. Thus we have established the first part of Proposition 4.2 for any $k \geq 2$. For $k > 2$, the inequality (4.29) has as an immediate consequence the inequality

$$m^{-1} > 2\tan\alpha + \cot\alpha = 2x^{-1} + x,$$

which is the growth rate condition (4.21). This concludes the proof of Proposition 4.2, when $k > 2$.

**Remark 4.3** *The case $k = 2$, the genus=1 case, turned out to be more difficult for us than the higher genus cases. In addition to the Monotonicity Lemma, we proved, for $k = 2$:*

$$\frac{E_+(\tilde{\alpha})}{E_-(\tilde{\alpha})} - \frac{E_+(\alpha)}{E_-(\alpha)} \leq \frac{2}{3}\sin^2\alpha.$$

*If we insert this in the compatibility condition (4.28), and use the identity $\sin^{-2}\alpha = \cot^2\alpha + 1$, we get, instead of (4.29),*

$$m^{-1} - m > \frac{\cot\alpha + \tan\alpha}{\sin^{-2}\alpha} - \cot\alpha + \tan\alpha$$

$$= \cot^3\alpha + \tan\alpha[\cot^2\alpha + 2]$$

$$= \cot^3\alpha + \cot\alpha + 2\tan\alpha$$

$$> 2\tan\alpha + \cot\alpha = 2x^{-1} + x,$$

*which implies the growth rate condition (4.21) in the case $k = 2$.*

### 4.9.1 Proof of the Monotonicity Lemma

To establish the claim

$$0 \leq \alpha \leq \gamma \leq \frac{\pi}{2} \Rightarrow \frac{E_+}{E_-}(\alpha) \geq \frac{E_+}{E_-}(\gamma)$$

it is equivalent to prove
$$E_-(\alpha) \cdot (E_-(\gamma) - E_+(\gamma)) - E_-(\gamma) \cdot (E_-(\alpha) - E_+(\alpha)) \geq 0 \qquad (4.30)$$
and sufficient to prove
$$E_-(\alpha) \cdot (E_- - E_+)'(\alpha) - E_-'(\alpha) \cdot (E_- - E_+(\alpha)) \geq 0. \qquad (4.31)$$
From the definition of these integrals in (4.27) we have
$$E_-'(\alpha) = \sin\alpha \cos\alpha \int_0^\alpha \frac{1}{\sqrt{\cos^2\phi - \cos^2\alpha}} \cos((1 - \frac{2}{k})\phi) d\phi$$
and using $\cos(\alpha - \beta) - \cos(\alpha + \beta) = 2\sin\alpha \sin\beta$ we have
$$(E_- - E_+)(\alpha) = \int_0^\alpha \sqrt{\cos^2\phi - \cos^2\alpha} \cdot 2\sin\phi \cdot \sin(\frac{2}{k}\phi) d\phi,$$
$$(E_- - E_+)'(\alpha) = \sin\alpha \cos\alpha \int_0^\alpha \frac{2\sin\phi \cdot \sin(\frac{2}{k}\phi)}{\sqrt{\cos^2\phi - \cos^2\alpha}} d\phi.$$

Now we define decreasing functions $f_1, g_1$ and increasing functions $f_2, g_2$ with choices of constants $c_i > 0$ such that $\int_0^\alpha g_i(\phi) d\phi = 1$:
$$0 \leq f_1(\phi) := \sqrt{\cos^2\phi - \cos^2\alpha}$$
$$0 \leq g_1(\phi) := c_1 \cdot \cos((1 - \frac{2}{k})\phi) \qquad \text{(decreasing)}$$
$$0 \leq f_2(\phi) := 1/\sqrt{\cos^2\phi - \cos^2\alpha}$$
$$0 \leq g_2(\phi) := c_2 \cdot 2\sin\phi \cdot \sin(\frac{2}{k}\phi) \qquad \text{(increasing)}.$$

With this notation (4.31) is equivalent to
$$\int_0^\alpha f_1 \cdot g_1 d\phi \cdot \int_0^\alpha f_2 \cdot g_2 d\phi \geq \int_0^\alpha f_1 \cdot g_2 d\phi \cdot \int_0^\alpha f_2 \cdot g_1 d\phi$$
and this inequality is implied by
$$\int_0^\alpha f_1 \cdot (g_1 - g_2) d\phi \geq 0, \qquad \int_0^\alpha -f_2 \cdot (g_1 - g_2) d\phi \geq 0.$$

Since $g_1 - g_2$ is decreasing, with $\int_0^\alpha (g_1 - g_2) = 0$, there exists a $\phi^* \in (0, \alpha)$ such that
$$\phi \leq \phi^* \Rightarrow (g_1 - g_2)(\phi) \geq 0,$$
$$f_1(\phi) \geq f_1(\phi^*) \geq 0 \text{ and } -f_2(\phi^*) \leq -f_2(\phi) \leq 0$$
$$\phi^* \leq \phi \Rightarrow (g_1 - g_2)(\phi) \leq 0$$
$$0 \leq f_1(\phi) \leq f_1(\phi^*) \text{ and } -f_2(\phi) \leq -f_2(\phi^*) \leq 0.$$

Therefore
$$\int_0^\alpha f_1(\phi)(g_1 - g_2)(\phi) d\phi \geq \int_0^\alpha f_1(\phi^*)(g_1 - g_2)(\phi) d\phi = 0$$
$$\int_0^\alpha -f_2(\phi)(g_1 - g_2)(\phi) d\phi \geq \int_0^\alpha -f_2(\phi^*)(g_1 - g_2) = 0.$$
$\square$

Figure 5.0.1 Four-ended embedded minimal surfaces of genus $2k, k > 0$.
These are the surfaces described in 1).

## 5  The structure of the space of examples

We begin this section by giving a list of the constructions to date of complete, embedded minimal surfaces with finite total curvature. None of these examples have been as completely analyzed as those of Theorem 3.3. For some, a complete existence proof independent of computation has yet to be found. For others there are existence proofs that amount to computing, with error estimates, the degree of the period map. All of these examples have the property that they are part of deformation families of surfaces, which begin with all but two of the ends flat and continue through surfaces all of whose ends are catenoidal. Most of the families, in contrast to the surfaces $M_{k,x}$ of Theorem 3.3 eventually contain surfaces that are not embedded because the growth rates of the ends change their order and a lower end overtakes a higher one.

1) Four ended examples with $k \geq 2$ vertical planes of symmetry, one horizontal plane of symmetry, genus $2(k-1)$, two flat ends and two catenoid ends [74]. (See Figure 5.0.1.)

2) A one parameter deforming the surfaces above through surfaces with four catenoid ends. Eventually these surfaces cease to be embedded. [73, 4, 74, 75].

3) The surfaces in 2) with symmetric tunnels through their waist planes [73, 4]. These surfaces have $k \geq 2$ vertical planes of symmetry and genus $3(k-1)$. The begins with a surface with two flat ends and moves through surfaces with all four ends catenoidal. Eventually the surfaces are not embedded. (See Figure 5.0.2.)

4) Five ended examples. [5]. (See Figure 5.0.3.)

These examples have in common the fact that they were discovered first with the help of insights guided by and guiding computation. The software MESH, developed by James T.

Figure 5.0.2 Four-ended embedded minimal surfaces.

Left: The surface described in 3) with $k = 1$.
Right: Deformation of the four-ended surface on the left, now with all ends catenoidal. The surface is no longer embedded.

Figure 5.0.3 Five-ended embedded minimal surface of genus 3, discovered by Boix and Wohlgemuth. See 4).

Hoffman et. al. ([8, 41]) together with various numerical routines for killing periods (finding simultaneous zeros of several real-valued functions of several variables) has been crucial to this research. This was also the case for the three-ended examples of Theorem 3.3. One has to produce Weierstrass data by methods analogous to those presented in Section 4.1–4.4. The main obstacle to constructing more complex examples with methods like those in Section 4 is the difficulty in killing the periods. Even numerically, for $n$ large, finding the simultaneous zero of $n$ functions on a domain in $R^n$ is time-consuming and unstable for functions (the periods given by Weierstrass integrals) that themselves are approximated by numerical integration.

The state of knowledge at present is rather poor. It is not known whether these surfaces can be perturbed further to reduce their symmetry. There is at present no general implicit-function-theorem-type result for the period mapping. Some important theoretical progress has been made by Wohlgemuth [74, 75]. The existence of certain of these surfaces (the most symmetric ones with flat ends) can be proved by theoretical means involving finding desingularizing limits of unstable solutions to the Plateau Problem. (See, for example [34].) However, while this is done for the surfaces of Theorem 3.2, this procedure has never produced an example not previously discovered by the method outlined above.

## 5.1 The space of complete, embedded minimal surfaces of finite total curvature

Let $\mathcal{M}$ be the space of complete, embedded minimal surfaces of finite total curvature, normalized so that all the ends have vertical normals and so that the maximum absolute value of the Gauss curvature is equal to one and occurs at the origin. (Note that the plane is not in $\mathcal{M}$.)

**Definition 5.1** *The subspace $\mathcal{M}_{i,j} \subset \mathcal{M}, i \geq 0, j \geq 2$, consists of all surfaces in $\mathcal{M}$ of genus $i$ with $j$ ends. We define*

$$\mathcal{M}^*_{k,r} = \bigcup_{i \leq k,\, j \leq r} \mathcal{M}_{i,j}.$$

We consider two elements of $\mathcal{M}$ to be the same if they differ by a Euclidean motion. From (2.21) it follows that the total curvature of an example in $\mathcal{M}_{ij}$ is equal to $-4\pi(1-(i+j))$. This means that the subspace of $\mathcal{M}$ with total curvature at least $-4\pi n$ is a union of a finite number of $\mathcal{M}^*_{k,r}$. In this terminology, we may restate some of Theorem 3.1 and Theorem 3.4 as follows:

i) $\mathcal{M}^*_{o,r}$ consists of one surface, the catenoid [51];

ii) $\mathcal{M}^*_{k,2}$ consists of one surface, the catenoid [67];

iii) $\mathcal{M}_{1,3}$ consists of the surfaces $M_{2,x}$ in Theorem 3.3, which is naturally a half-open interval $[1, \infty)$ [16].

Note from this last item that the spaces $\mathcal{M}_{ij}$ are not necessarily compact. However it can be shown that in $\mathcal{M}$, which we consider with the compact-open topology, every divergent sequence of surfaces $M_{k,x}$ has a subsequence that converges to the catenoid [30]. Thus while $\mathcal{M}_{1,3}$ is not compact, $\mathcal{M}^*_{1,3}$ is naturally a closed interval and is compact. (Its end points are Costa's surface and the catenoid.) The following result generalizes this fact.

**Theorem 5.1 ([30])** $\mathcal{M}^*_{k,r}$ *is compact. Specifically, every sequence of surfaces in $\mathcal{M}^*_{p,k}$ possesses a subsequence that converges smoothly on compact subsets of $R^3$ to a surface in $\mathcal{M}^*_{p,k}$.*

This theorem is of limited utility because at the present time it is not known, in any general fashion, which of the spaces $\mathcal{M}_{i,j}$ are empty and which are not.

## 5.2 Some questions and conjectures

One could imagine that the Costa example is constructed in the following heuristic manner. Take a catenoid and intersect it by its waist plane. Remove a tubular neighborhood of the intersection circle and smoothly join up the four boundary circles by a genus one minimal surface. The higher-genus examples $M_k, k \geq 3$ can be thought of as improved approximations.

**Theorem 5.2 ([38])**

i) *As normalized in Theorem 3.2, the surfaces $M_k, k \geq 2$, possess a subsequence that converges as $k \to \infty$ to the union of the plane and the catenoid. The convergence is smooth away from the intersection circle.*

ii) *Considered as surfaces in $\mathcal{M}$ (i.e., normalized so that the maximum value of $K$ is equal to 1 and occurs at the origin), the sequence $\{M_k\}$ possesses a subsequence that converges to Scherk's singly-periodic surface.*

It has been conjectured by Hoffman and Meeks that for a complete, embedded, nonplanar minimal surface of finite total curvature,

$$k \geq r - 2, \tag{5.1}$$

where $k$ is the genus and $r$ is the number of ends. Theorem 3.1, Statements 2 and 3, results of Lopez-Ros and of Schoen, can be interpreted as saying that equality holds in (5.1) when either $k = 0$ or $r = 2$. The examples of Theorems 3.2 and 3.3, as well as all other known examples satisfy this inequality. Inequality 5.1 may be stated as a conjecture

**Conjecture 5.1**

$$\mathcal{M}_{k,r} \text{ is empty if } k < r - 2\,.$$

Figure 5.2.0

One of a family of complete minimal surfaces of genus zero with two parallel, embedded catenoid ends and one embedded flat end, at which the Gauss map has an order-two branch point. The surface is not embedded because the flat end is not parallel to the catenoid ends.

The heuristic idea behind the conjecture is that to increase the number of ends ($r$) by adding a flat end and then looking for a Costa style surface like the surface $M_k$, which is nearby and embedded, forces the genus ($k$) to increase. Related to this conjecture is the question of whether the order of the Gauss map at a flat end of a complete, embedded minimal surface can be equal to two. According to Remark 2.3*iii)* the order must be at least two, and in all known *embedded* examples of finite total curvature the order is at least three.

(See Remark 2.8.) However, there are complete embedded periodic minimal surfaces with flat ends where the Gauss map has degree two. The most famous of these is the example of Riemann [62, 63]. (See also [31], [37] and [39].) For *immersed* minimal surfaces, it is possible to have such an end. Consider the Weierstrass data

$$g(z) = \rho((z-r)(z+r))^{-1}, dh = (z^2 - r^2)(z^2 - 1)^{-2} dz, r \neq 1,$$

on $S^2 - \{\pm 1, \infty\}$. The surface has sufficient symmetry (vertical plane of reflection and one horizontal line) to make the period problem one-dimensional. The period problem is solvable by a residue calculation, which determines the choice of $r$ (for each choice of $\rho$ within a prescribed range). The ends at $\pm 1$ are catenoid ends and the end at $g(z) = \infty$, where $g(z)$ has a zero of order 2, is flat. However, while all the ends are separately embedded, they are not parallel; the surface is complete and immersed but not embedded.

Figure 5.2.1 The Horgan Surface

The surface depicted here with genus two and three ends probably does not exist. The periods appear to be zero, but one of them is just very small. Computationally, it appears that this period always has the same sign and is never zero, although it becomes vanishingly small. However, in the limit, the surface degenerates.

Recall from Proposition 2.1) and Remark 2.4*i)* that each embedded catenoid end of a complete embedded minimal surface of finite total curvature is asymptotic to the end of a specific vertical catenoid. There are always at least two such ends, according to Proposition 2.5*v)* and the discussion following Theorem 2.4. Each such end has an axis determined by the term of order $\mathcal{O}(\rho^{-1})$ in (2.21) and made precise in Definition 2.7. Because all of the known examples can be constructed by symmetry methods like those used in Section 4, it follows that all of the axes of the catenoid ends of these examples coincide. One basic question is whether or not there exist examples whose catenoid ends have *different* axes. Such an example could have, at most, a symmetry group of order four, generated by a reflection through a plane containing all of the axes (if such a plane exists) and an order two rotation about a horizontal line (if such a line exists). If the axes were not coplanar, only the rotation could exist. Such a rotation would imply that the number of catenoid ends was even.

As observed in Remark 2.8, if there are three catenoid ends then their axes are coplanar and the end whose axis is in the middle has logarithmic growth with sign opposite that of the other two.

It is known that any intrinsic isometry of a complete *embedded* minimal surface of finite total curvature in $R^3$ must extend to a symmetry of that surface; that is an isometry of $R^3$ that leaves the surface invariant [13, 36]. If the surface is not embedded, this is not necessarily the case. The Enneper surface and its generalizations $g(z) = z^k, dh = g(z)dz, z \in \mathbf{C}$, have intrinsic isometric rotations $z \to uz$, $|u| = 1$, but only $k$ of them extend to symmetries of the surface. ( See Remark 2.3ii).

An even more basic problem concerns the question of whether or not there exists a complete embedded minimal surface with no symmetries except the identity. That all the known examples have significant symmetry groups is an artifact of their construction and says nothing about the general situation. It would be extremely surprising if the hypotheses of finite total curvature together with embeddedness implied the existence of a nontrivial symmetry.

A related question concerns the surfaces $M_{k,x}$ in Theorem 3.3. Costa's Theorem 3.1.5 states that for the genus-one examples ($k = 2$)

$$\mathcal{M}_{1,3} = \{M_{2,x} \mid 0 < x \leq 1\}. \tag{5.2}$$

It would be very interesting to know to what extent this is true for higher genus:

$$\mathcal{M}_{k,3} \stackrel{?}{=} \{M_{k+1,x} \mid 0 < x \leq 1\}. \tag{5.3}$$

It is not even known if this is true with the added hypothesis that the symmetry group be of order at least $2(k+1)$, as is the case for the surfaces $M_{k+1,x}$. Costa's proof of (5.2) uses the full force of elliptic function theory, and hence is restricted to genus one. A geometric proof of (5.2) would be extremely useful to have, in order to begin to understand what happens in the higher-genus, three-ended case.

## 6 Finite total curvature versus finite topology

A surface has finite topology provided it is homeomorphic to a compact Riemann surface, from which a finite number of points have been removed. As a consequence of Theorem 2.2, we know that a complete immersed minimal surface whose total curvature is finite must have finite topology. In fact, it is *conformal* to a punctured compact Riemann surface. In this section, we will discuss to what extent finite topology implies finite total curvature for complete properly embedded minimal surfaces.

The helicoid is a complete, embedded, simply-connected (genus zero and one end) minimal surface. Because it is periodic and not flat, its total curvature is infinite. This simple classical example shows that finite topology does not imply finite total curvature. Osserman has posed the question of whether or not the helicoid is the only nonplanar, complete, embedded, minimal surface of genus-zero with one end. According to Theorem 3.1 there is no such surface of finite total curvature of any genus. One could ask, more generally, if there are any complete embedded minimal surfaces of finite topology and one end, besides the plane and the helicoid. For example, does there exist such a surface that is morphologically the helicoid with a finite number of handles? (We have recently shown that the answer is yes. See Section 6.3.)

An end of a complete surface of finite topology is necessarily annular. That is, it is homeomorphic to a punctured disk. If it is conformal to a punctured disk and the total curvature is infinite, then the Gauss map must have an essential singularity. All known embedded annular ends are asymptotic to the end of a helicoid. No embedded end is known that is conformal to an annulus $r_1 \leq |z| < r_2$.

## 6.1 Complete, properly-immersed, minimal surfaces with more than one end

Under the assumption that the surface has more than one end, much more is known.

**Theorem 6.1 (Hoffman and Meeks [35])** *Let $M$ be a complete properly embedded minimal surface in $R^3$ with more than one topological end. Then at most two of the annular ends have infinite total curvature.*

**Corollary 6.1 ([35])** *Let $M$ be a complete embedded minimal surface in $R^3$ with finite topology. Then all except two of the ends have finite total curvature and are asymptotic to planes or catenoids, all with parallel limit normals.*

This Corollary was strengthened by Fang and Meeks who proved

**Theorem 6.2 (Fang and Meeks [21])** *Let $M$ be a properly embedded minimal surface in $R^3$ with two annular ends that have infinite total curvature. Then these ends lie in disjoint closed halfspaces and all other annular ends are flat ends parallel to the boundary of the halfspaces.*

**Corollary 6.2** *Suppose $M$ is a complete properly embedded minimal surface of finite topology. If $M$ does not have finite total curvature, then there are three possibilities:*

*i) $M$ has a catenoid end and exactly one end with infinite total curvature that lies in a halfspace;*

*ii) $M$ has exactly one end (of infinite total curvature) and this end does not lie in a halfspace;*

*iii) $M$ has two ends of infinite total curvature that lie in disjoint closed halfspaces, and a finite number of flat ends in the slab between these halfspaces. There are no catenoid ends.*

In particular a complete properly embedded minimal surface of finite topology with more than one end can never have an end asymptotic to a helicoid end (because such an end never lies in a halfspace).

## 6.2 Complete embedded minimal surfaces of finite topology and more than one end

What sorts of infinite total curvature ends can appear on a complete embedded minimal surface with finite topology? Nitsche considered ends that were fibred by embedded Jordan curves in parallel planes. These are called *Nitsche ends* by Meeks and Rosenberg [53]. For concreteness we assume that such an end $A$ lies in the halfspace $\{x_3 \geq 0\}$ with boundary a simple closed curve in $x_3 = 0$, and that $\{x_3 = t\} \cap A$ is a simple closed curve, for all $t \geq 0$. Nitsche proved

**Proposition 6.1 ([57])** *A Nitsche end, all of whose level curves are star shaped, is a catenoidal end.*

In particular, it has finite total curvature. A Nitsche end necessarily has the conformal type of a punctured disk. Certainly, we may assume that the end is conformally $\{z \,|\, 0 \leq \alpha < |z| \leq 1\}$. The height function $h(z) = x_3(z)$ is harmonic and so is essentially determined by its boundary values. Suppose $\alpha > 0$. Then $h = 0$ on $|z| = 1$ and, $h = \infty$ on $|z| = \alpha$. Since $c_1(\ln|z|)$ is also harmonic, and $c_1 \ln|z|$ lies above $h$ on the boundary of $\alpha < |z| \leq 1$, it follows that
$$h(z) < c_1 \ln|z| \qquad \alpha < |z| < 1$$
for any real value of $c_1$. This means $h(z) \equiv -\infty$ on the interior of the annulus, which is a contradiction. Hence we can take the disk to be $\{z \,|\, 0 < z \leq 1\}$, and $h(z) = c_1 \ln|z| c_1 < 0$. Thus $dh$ is determined and using the Weierstrass representation (2.7), (2.8), we can determine the end completely by the knowledge of its Gauss map. For the Nitsche end to have infinite total curvature, $g$ must have an essential singularity. Thus finding a non-catenoidal Nitsche end is equivalent to finding an analytic function on the punctured unit disk with an essential singularity at the origin for which (2.7) (2.8) provide a well-defined embedding. Toubiana and Rosenberg have given examples which are well-defined, but not embedded [66]. We also remark that on $\mathbb{C}$, $g(z) = \exp(F(z))$, and $dh = dz$, where $F(z)$ is entire, gives a complete, immersed, simply-connected minimal surface whose Gauss map has an essential singularity at infinity.

The importance of understanding Nitsche ends is underscored by the following strong result of Meeks and Rosenberg.

**Theorem 6.3 ([53])** *Let $M$ be a complete embedded minimal surface of finite topology and more than one end. Then any end of infinite total curvature is a Nitsche end. In particular $M$ is conformal to a compact Riemann surface punctured in a finite number of points.*

As a result of this theorem, a proof of the Nitsche conjecture would imply that finite topology and finite total curvature are equivalent for complete embedded minimal surfaces with at least two ends.

Figure 6.3 The genus one helicoid, $\mathcal{H}e_1$

## 6.3  Complete embedded minimal surfaces of finite topology with one end

In 1992-3, Hoffman, Karcher and Wei [32, 33] constructed a complete embedded minimal surface of finite topology and infinite total curvature.

**Theorem 6.4** *There exists a complete embedded minimal surface, $\mathcal{H}e_1$, of genus-one with one end that is asymptotic to the end of a helicoid. The surface contains two lines, one vertical and corresponding to the axis of the helicoid, the other horizontal crossing the axial line. Rotation about these lines generates the full isometry group of $\mathcal{H}e_1$, which is isomorphic to $Z_2 \oplus Z_2$.*

This is the first properly embedded minimal surface with infinite total curvature, whose quotient by symmetries does not have finite total curvature.

# 7  Stability and the index of the Gauss map

Recall that a minimal surface $S \subset \mathrm{R}^3$ is by definition a critical point of the area functional for compactly supported perturbations. Let $X \colon \Sigma \to \mathrm{R}^3$ be a conformal parametrization of $S$ and $N \colon \Sigma \to S^2$ its Gauss map. If $D$ is a compact subdomain of $\Sigma$, and $X_t \colon D \to \mathrm{R}^3$ is a normal variation of $D$ with variation vector-field $fN$, $f \in C_0^\infty(D)$, the second derivative of area is given by

$$\frac{d^2}{dt^2} Area(X_t) = \int_D (|\nabla f|^2 + 2Kf^2) dA =: Q_D(f,f), \qquad (7.1)$$

where $K$ is the Gauss curvature. The second variation operator is defined to be

$$L =: -\Delta + 2K, \qquad (7.2)$$

in terms of which

$$Q_D(f,f) = \int_D fLf\, dA.$$

The sign of $\Delta$ is determined by our definition in (2.3). The index of $L$ on $D$, denoted $ind(D)$, is defined to be the dimension of the maximal subspace of $C_0^\infty(D)$ on which $Q_D$ is negative definite. A domain $X(D) \subset S$ is said to be *stable* if $ind(D) = 0$, i.e., if (7.1) is nonnegative for all $f \in C_0^\infty(D)$.

For a non-compact minimal surface $X(\Sigma) = S \subset \mathrm{R}^3$, the index $(\Sigma)$ is defined to be the supremum of $ind(D)$ taken over all compact $D \subset \Sigma$. The following theorem summarizes the results of several authors relating index to stability.

**Theorem 7.1** *Let $X \colon \Sigma \to R^3$ be a conformal minimal immersion of an orientable surface.*

i) *(Schwarz [68].) If the Gaussian image of $\Sigma$ lies in a hemisphere, $S$ is stable;*

ii) *(Barbosa and do Carmo [1].) If the area of the Gaussian image satisfies $A(N(\Sigma)) < 2\pi$, then $S$ is stable. (Here $A(N(\Sigma))$ is the area of the Gaussian image* disregarding multiplicities*);*

iii) *(Fischer-Colbrie and Schoen [23], do Carmo and Peng [18], Pogorelov [61].) If $S$ is complete, $S$ is stable if and only if $S$ is a plane;*

iv) *(Fischer-Colbrie [22], Gulliver [25], Gulliver and Lawson [26].) Index $(\Sigma) < \infty$ if an only if $\int_\Sigma |K| dA < \infty$.*

Without going too much into detail, we will give some of the background to the last two parts of the theorem above. The Gauss curvature $K$ can be expressed as $-|\nabla N|^2/2$ and therefore we may write $L = -(\Delta + |\nabla N|^2)$. Observe that $Q_D(f,f) = \int_D (|\nabla f|^2 - |\nabla N|^2 f^2) dA$ is invariant under conformal change of metric and therefore $ind(D)$ depends only on $N$.

We may therefore define, for a meromorphic mapping $\Phi\colon \Sigma \to S^2$, the index $ind_\Phi(D)$ of a compact subdomain $D \subset \Sigma$ to be the number of negative eigenvalues of $L_\phi = -(\Delta + |\nabla \phi|^2)$, with Dirichlet boundary conditions on $D$. By definition, $ind(D) = ind_N(D)$. In particular if $\overline{\Sigma}$ is a compact Riemann surface and $\Phi$ a meromorphic function on $\overline{\Sigma}$,

$$index(\Phi) := ind_\Phi(\overline{\Sigma}) \tag{7.3}$$

is an integer invariant of $\Phi\colon \overline{\Sigma} \to S^2$.

For a complete immersed minimal surface $S \subset \mathrm{R}^3$, the condition of finite total curvature is equivalent to the conditions that the underlying Riemann surface $\Sigma$ is a compact Riemann surface $\overline{\Sigma}$ with a finite number of points removed, and that the Gauss map $N$ extends to a meromorphic function on $\overline{\Sigma}$ (See Theorem 2.2). Thus for any compact subdomain $D \subset \Sigma$,

$$ind(D) = ind_N(D) \leq index(N),$$

In particular the condition of finite total curvature implies that

$$index(\Sigma) = \sup\{ind(D) \mid D \subset \Sigma, D\,compact\}$$

is finite. This is the easier implication of part *iii)* of Theorem 7.1. In fact these two indices are equal:

**Proposition 7.1 ([22])** *Let $X\colon \Sigma \to R^3$ be a complete immersed minimal surface of finite total curvature. Let $N\colon \Sigma \to S^2$ be its Gauss map. Then $ind(\Sigma) = index(N)$.*

As a consequence of this Proposition, statement iii) of Theorem 7.1 follows from the fact that the index of any nonconstant meromorphic function $\Phi\colon \overline{\Sigma} \to S^2$ is nonzero.

**Remark 7.1** *Let $ds_N^2$ be the metric produced by pulling back the standard metric on $S^2$ by the Gauss map $N\colon \Sigma \to S^2$. If $ds^2$ denotes the metric induced by $X\colon \Sigma \to R^3$, we have $ds_N^2 = \frac{1}{2}|\nabla N|^2 ds^2 = -K ds^2$. This metric is regular, away from the isolated zeros of $K$, and has constant curvature 1. The volume of $\overline{\Sigma}$ in this metric is therefore $4\pi\, degree\,(\Phi)$. Because the Dirichlet integrand $|\nabla f|^2 dA$ is invariant under conformal change of metric*

$$Q(f,f) = \int (|\nabla f|^2 + 2Kf^2) dA = \int (|\widetilde{\nabla} f|^2 - 2f^2) d\widetilde{A},$$

*where $\widetilde{\nabla} f$ and $d\widetilde{A}$ are computed in the metric $ds_N^2$.*

*The associated differential operator is $\widetilde{L} = -(\widetilde{\Delta} + 2)$. Eigenfunctions of $-\widetilde{\Delta}$ can be defined by the standard variational procedure [54, 72]. The spectrum of $-\widetilde{\Delta}$ is discrete, infinite and nonnegative and we may think of $index\,(\Phi)$ as the number of eigenvalues of $-\widetilde{\Delta}$ that are strictly less than two.*

It is natural to try to compute the index of stability of the minimal surfaces mentioned in this survey. We have seen that this is equivalent to computing the index of their extended Gauss maps. For the simplest case of the identity map from $S^2$ to itself, the index is just the number of eigenvalues of $\widetilde{\Delta} = \Delta$ (the spherical Laplacian) less than two. Since the first nonzero eigenvalue is equal to 2 and $\lambda_0 = 0$ is a simple eigenvalue, it follows from Remark 7.1 above that $index(Id) = 1$.

**Proposition 7.2** *The index of stability of Enneper's surface and the catenoid is equal to one.*

This follows immediately from the fact that the extended Gauss map of these genus-zero surfaces is a conformal diffeomorphism to the sphere. Osserman [58, 59] proved that the catenoid and Enneper's surface were the only complete minimal surfaces with this property. It is natural to ask whether or not index one characterizes these surfaces among all complete minimal surfaces of finite total curvature. The following result of Montiel-Ros ([54]) can be used to resolve this question.

**Theorem 7.2** *Suppose $\Phi\colon \overline{\Sigma} \to S^2$ is a nonconstant holomorphic map of degree $d$ defined on a compact Riemann surface of genus $k$.*

*i) If $k = 0$, $index(\Phi) = 1 \iff d = 1$.*

*ii) If $k = 1$, $index(\Phi) \geq 2$.*

*iii) If $k \geq 2$ then, for a generic conformal structure, $index(\Phi) \geq 2$.*

*iv) If $index(\Phi) = 1$, $d \leq 1 + \left[\frac{1+k}{2}\right]$.*

Statements *i)* and *ii)* are straightforward. Statements *iii)* and *iv)* follow from more general results in [54] (see Corollary 8 and Theorem 6) and the Brill-Noether Theorem.

Now suppose that $S$ is a complete minimal surface whose index of stability is equal to 1. From (2.9), we have for the Gauss map $N$ of a complete minimal surface $X\colon \Sigma \to \mathbf{R}^3$ of genus $k$
$$d = degree(N) \geq k + 1\,.$$

Combined with statement *iv)* of Theorem 7.2, this implies $k \leq 1$. Then statements *i)* and *ii)* imply that $k = 0$ and $d = 1$. As noted above, the catenoid and Enneper's surface are the only possibilities. Thus we have

**Corollary 7.1 (Corollary 9 of [54] [50])** *A complete minimal surface with stability index equal to one must be the catenoid or Enneper's surface.*

**Remark 7.2** *Montiel and Ros conjectured that it was possible to have a hyperelliptic Riemann surface whose two-fold covering map to $S^2$ had index equal to one, provided the branch values were sufficiently well distributed. This was proved by Souam [69].*

**Theorem 7.3** *Given $\varepsilon > 0$, choose $\mathcal{P} = \{p_1 \ldots p_{n_\varepsilon}\} \subset S^2$ such that the open $\varepsilon$-disks about the $p_i$ are disjoint, but maximal in the sense that their complement does not contain an $\varepsilon$-disk. Suppose $\Sigma$ is a hyperelliptic Riemann surface and $\phi\colon \Sigma \to S^2$ a two-fold covering of $S^2$ whose branch values are precisely $\mathcal{P}$. Then if $\varepsilon > 0$ is small enough, $index(\Phi) = 1$.*

The following results allows us to compute the index of stability of the Meeks-Jorge n-noids, the Chen-Gackstatter surface and its generalizations to higher symmetry (See Remark 4.2).

**Proposition 7.3 (Montiel and Ros [54])** *Suppose $\phi\colon \overline{\Sigma} \to S^2$ is a degree $d \geq 1$ holomorphic map. If all the branch value of $\phi$ lie on a great circle*

$$index(\Phi) = 2d - 1. \tag{7.4}$$

Choe [12], proved a weaker version of this proposition.

**Theorem 7.4 ([54])** *The index of stability of:*

   i) *the n-noid is $2n - 3$;*

   ii) *the Chen-Gackstatter surface is $3$.*

**Proof.** From Section 2.2 we know that the Gauss map of the Chen-Gackstatter surface has degree $d = 2$, and all of its branch values are real. The n-noid $n \geq 2$ can be produced with the Weierstrass data $\overline{\Sigma} = S^2$,

$$g = z^{n-1}$$
$$dh = (z^n + z^{-n} + 2)z^{-1}dz.$$

(When $n = 2$ we have the catenoid.) In both cases we can apply Proposition 7.3. □

The equivalence of finite total curvature and finite index suggests that there should be an explicit relationship between $d$, the degree of the Gauss map, and the index of stability. Tysk [72] was the first to show that the index could be estimated by the degree:

$$\frac{Index(N)}{degree(N)} < c \text{ where } c \sim 7.68.$$

Choe [12] conjectured that a weakened version of (7.3) was true in general for all complete minimal surfaces $M$ of finite total curvature $-4\pi d$:

**Conjecture 7.1**
$$index(M) \leq 2d - 1 \qquad (7.5)$$

The conjectured inequality (7.5) cannot be replaced in general by an equality. Nayatani (private communication) has pointed out to us that an example written down by Kusner [44], Rosenberg-Toubiana [64, 65], R. Bryant [6, 7] with 4 flat ends, genus 0, and Gauss map of degree $d = 3$, gives strict inequality. On any surface the support function $u = X \cdot N$ satisfies $Lu = 0$, where $L$ is defined in (7.2), and because the ends are flat, $u$ is bounded. Using this fact, it can be shown that the index of this surface is 4. Since $d = 3$, index = $4 < 5 = 2d - 1$.

Nayatani [55] computed the index of the Costa surface, $M_2$ of Theorem 3.2. This was done by computing the index of its Gauss map, which is the derivative of the Weierstrass $\wp$-function, geometrically normalized, on the square torus. The geometric normalization does not change the index.

He subsequently improved that result:

**Theorem 7.5 (Nayatani [56])** *Let $M_k$ be the Hoffman-Meeks surface of genus $k-1$, described in Theorem 3.2. For $2 \leq k \leq 38$,*

$$index(M_k) = 2k + 1\,.$$

The Gauss map of $M_k$ has degree $k+1$, so equality holds in (7.5) for these surfaces. For genus 0 surfaces Montiel and Ros [54] and Ejiri and Kotani [20] independently proved the following result.

**Theorem 7.6** *Formula (7.4) holds for a generic complete, genus zero, minimal surface of total curvature $-4\pi d$. In general, $index(M) \leq 2d - 1$ for such a surface.*

This result holds even for branched minimal surfaces and involves an analysis of bounded Jacobi fields.

We end this chapter with a discussion of the work of Choe [12] about the *vision number* of a minimal surface with respect to a vector field $V$ on $\mathrm{R}^3$. For a vector field $V$ on $S \subset \mathrm{R}^3$, Choe defines the *horizon* of $V$ to be

$$H(S, V) = \{s \in S \mid V \in T_s S\}\,,$$

and the vision number, to be

$$v(S, V) = \#\text{ components of } S - H(S, V)\,.$$

Suppose $V$ is a variation vector field associated with a of parallel translations, or rotations about a line. (Both are Killing vector fields on $\mathrm{R}^3$). Consider the restriction of such a vector field to $S$, and let $\widehat{V}$ be its projection onto the normal bundle of $S$. Then $\widehat{V}$ is a Jacobi field. That is, $\widehat{V} \cdot N$ is in the kernel of $L = -\Delta + 2K$. The same holds for the position vecter field $V_0(p) = \overline{p - p_0}$, for any fixed origin $p_0 \in \mathrm{R}^3$. This can be used to show that any connected subset of $S - H(S,V)$ is stable, for any one of these vector fields $V$. Choe uses this to prove

**Theorem 7.7** *Let $S \subset \mathrm{R}^3$ be a complete minimal surface of finite total curvature $-4\pi d$. Let $V_T, V_R$ denote the Killing vector fields generated by translation and rotation, respectively. Then*

   *i) $Index(S) \geq v(S, V_T) - 1$;*

   *ii) $d \geq v(S, V_T)/2$;*

*If all the ends are parallel*

   *iii) $Index(S) \geq v(S, V_R) - 1$;*

*If all the ends are embedded and parallel*

   *iv) $d \geq v(S, V_R)/2$.*

From statement *iii)* of the above theorem we can get a lower bound for the index of the surfaces $M_k$ of Theorem 3.2. Observe that the same is true for the $M_{k,\alpha}$ of Section 4. These surfaces have genus $k-1$, parallel embedded ends and $k$ vertical planes of reflective symmetry, which intersect in a vertical line, $\mathcal{L}$, orthogonal to the ends, and which divide the surface into $2k$ regions. Their Gauss maps have degree $d = k+1$. If $V_R$ is the Killing field generated by rotation about $\mathcal{L}$, then it is clear that $H(M_{k,\alpha}, V_R)$ contains the intersection of $M_{k,\alpha}$ with each plane of symmetry. Hence $v(M_{k,\alpha}, V_R) \geq 2k$. By Theorem 7.7 *ii)* we have

**Corollary 7.2** *The index of stability of the surfaces $M_k = M_{k,0}$ and $M_{k,\alpha}$ satisfies*

$$Index(M_{k,\alpha}) \geq 2k - 1 = 2(genus\ M_{k,\alpha}) + 1\,. \tag{7.6}$$

**Remark 7.3** *The symmetry of the surfaces $M_{k,\alpha}$ imply that $v(M_{k,\alpha}, V_R)$ is an integer multiple of $2k$. But by statement* iv) *of Theorem 7.7, $2(k+1) = 2d \geq v(M_{k,\alpha}, V_R)$. Hence $v(M_{k,\alpha}, V_R) = 2k$. Thus this particular argument cannot be improved to give a better lower bound. As we know from Theorem 7.5, $index(M_k) = 2k+1, 2 \leq k \leq 38$, so this estimate is not sharp in this range.*